\documentclass[12pt,a4paper,reqno]{amsart}

\usepackage[all,2cell,arrow,color]{xy}
\newdir{ >}{{}*!/-10pt/@{>}}
\usepackage[english]{babel}
\usepackage{mathtools}
\usepackage{stmaryrd}
\usepackage{tensor}
\numberwithin{equation}{section}
\usepackage{graphicx} 
\usepackage{rotating}
\usepackage[pagebackref,colorlinks,urlcolor=darkgreen,citecolor=darkgreen,linkcolor=darkgreen]{hyperref}
\usepackage{enumerate}
\usepackage{xspace}
\usepackage{xcolor}
\usepackage{floatpag}
\usepackage{MnSymbol}
\usepackage{bbm}
\usepackage[normalem]{ulem}

\definecolor{darkgreen}{rgb}{0,0.45,0} 

  \newtheorem{proposition}{Proposition}[section]

  \theoremstyle{definition}
  \newtheorem{definition}[proposition]{Definition}

\theoremstyle{remark}

  \newcounter{c}
  
  \newcommand{\etyk}[1]{\vspace{-7.4mm}$$\begin{equation}\Label{#1}
  \addtocounter{c}{1}}
  \renewcommand{\]}{\ifnum \value{c}=1 $$\else \end{equation}\fi}
\newcommand*{\textcal}[1]{%
  \textit{\fontfamily{qzc}\selectfont#1}%
}

\newenvironment{amssidewaysfigure}
  {\begin{sidewaysfigure}\vspace*{.5\textwidth}\begin{minipage}{\textheight}\centering}
  {\end{minipage}\end{sidewaysfigure}}

\newcommand{\longdownarrow}[1]{\mbox{\rotatebox{270}{$\Longrightarrow$} \raisebox{-12pt}{$#1$}}}
\newcommand{\Longdownarrow}[1]{\displaystyle #1}
\newcommand{\Mnd}{\mathbb M\mathsf{nd}}
\newcommand{\Sqr}{\mathbb S\mathsf{qr}}
\newcommand{\ldb}{\lsem}
\newcommand{\rdb}{\rsem}
\newcommand{\calone}{\scalebox{1.15}{$\textcal{1}$}}
\newcommand{\udot}{\raisebox{3pt}{$\cdot$}}
\newcommand{\tp}{\hspace{-1pt} , \hspace{-2pt}}
\newcommand{\TwoCat}{\mathsf{2{\text -}Cat}}

\begin{document}

\title{The Gray monoidal product of double categories}

\author{Gabriella B\"ohm} 
\address{Wigner Research Centre for Physics, H-1525 Budapest 114,
P.O.B.\ 49, Hungary}
\email{bohm.gabriella@wigner.mta.hu}
\date{Jan 2019}
%\subjclass{}
 
\begin{abstract}
The category of double categories and double functors is equipped with a symmetric closed monoidal structure. For any double category $\mathbb A$, the corresponding internal hom functor $\ldb\mathbb A,-\rdb$ sends a double category $\mathbb B$ to the double category whose 0-cells are the double functors $\mathbb A \to \mathbb B$, whose horizontal and vertical 1-cells are the horizontal and vertical pseudotransformations, respectively, and whose 2-cells are the modifications. Some well-known functors of practical significance are checked to be compatible with this monoidal structure.
\end{abstract}
  
\maketitle

%%%%%%%%%%%%%%%%%%     INTRODUCTION  %%%%%%%%%%%%%%%%%%%%%%%%

\section*{Introduction} \label{sec:intro}

The category $\TwoCat$ of 2-categories and 2-functors carries different monoidal structures. The simplest one is given by the Cartesian product. It is symmetric and closed. For any 2-category $\mathcal A$, the internal hom functor $\langle\mathcal A,-\rangle$ sends a 2-category $\mathcal B$ to the 2-category of 2-functors $\mathcal A\to \mathcal B$, 2-natural transformations, and modifications. This is, however, often too restrictive. For example, important examples of 2-categories which are intuitively monoidal, fail to be monoids for that \cite{BaezNeuchl,KapranovVoevodsky,SchommerPries}. A well established generalization is the so-called {\em Gray monoidal product} \cite{Gray}. It is also symmetric and closed and for any 2-category $\mathcal A$ the corresponding internal hom functor $[\mathcal A,-]$ sends a 2-category $\mathcal B$ to the 2-category of 2-functors $\mathcal A\to \mathcal B$, {\em pseudonatural} transformations, and modifications.
The Cartesian monoidal structure is more restrictive than the Gray one in the sense that the identity functor on $\TwoCat$ is a monoidal functor from the former to the latter one.

The category $\mathsf{DblCat}$ of double categories and double functors is also symmetric closed monoidal via the Cartesian product $\times$. For any double category $\mathbb A$, the corresponding internal hom functor $\langlebar\mathbb A,-\ranglebar$ sends a double category $\mathbb B$ to the double category whose 0-cells are the double functors $\mathbb A \to \mathbb B$, whose horizontal and vertical 1-cells are the horizontal and vertical transformations, respectively, and whose 2-cells are the modifications; see \cite{GrandisPare}. The analogue of the Gray monoidal product on  $\mathsf{DblCat}$, however, has apparently not yet been discussed in the literature. The current paper addresses this question.

For any double categories $\mathbb A$ and $\mathbb B$, there is a bigger double category $\ldb\mathbb A,\mathbb B\rdb$ in which the 0-cells are still the double functors $\mathbb A\to \mathbb B$. The horizontal and vertical 1-cells are, however, the horizontal and vertical {\em pseudo} (or strong) transformations of \cite{GrandisPare}. The 2-cells are their modifications. In Section \ref{sec:existence} we prove that for any double categories $\mathbb A$ and $\mathbb B$, there is a representing object $\mathbb B \otimes \mathbb A$ of the functor $\mathsf{DblCat}(\mathbb B,\ldb\mathbb A,\mathbb - \rdb):\mathsf {DblCat} \to \mathsf{Set}$.
Constructing the associativity and unit constraints, as well as the symmetry, in Section \ref{sec:coherence} we show that $\otimes$ equips $\mathsf {DblCat}$ with a symmetric monoidal structure. In order to support this choice of monoidal structure on $\mathsf {DblCat}$, in Section \ref{sec:examples} monoidality of the following functors is checked.
\begin{itemize}
\item The identity functor $(\mathsf {DblCat},\times) \to (\mathsf {DblCat},\otimes)$.
\item The functors $(\mathsf {DblCat},\otimes) \to (\TwoCat,\otimes)$ sending double categories to their horizontal -- or vertical -- 2-categories (for the Gray monoidal product $\otimes$ on $\TwoCat$).
\item The square (or quintet) construction functor $\Sqr:(\TwoCat,\otimes) \to (\mathsf {DblCat}, \otimes)$ due to Ehresmann \cite{Ehresmann}.
\item The functor $\Mnd:(\mathsf {DblCat},\otimes) \to (\mathsf {DblCat},\otimes)$, sending a double category to the double category of its monads by Fiore, Gambino and Kock \cite{FioreGambinoKock}.
\end{itemize}
We also give an explicit description of monoids in $(\mathsf {DblCat},\otimes)$ which generalize the strict monoidal double categories of \cite{BruniMeseguerMontanari}; that is, the monoids in $(\mathsf {DblCat},\times)$.

\subsection*{Acknowledgement} 
Financial support by the Hungarian National Research, Development and Innovation Office – NKFIH (grant K124138) is gratefully acknowledged.  

%%%%%%%%%%%%%%%%%%%%  SEC 1   %%%%%%%%%%%%%%%%%%%%%%%%%%%%

\section{Existence}
\label{sec:existence}

In this section we construct an adjunction $-\otimes \mathbb D \dashv \ldb \mathbb D,-\rdb$ of endofunctors on the category $\mathsf{DblCat}$ of double categories, for any double category $\mathbb D$. Our line of reasoning is similar to \cite[Proposition 3.10]{BourkeGurski}. The occurring double functor $\otimes:\mathsf{DblCat} \times \mathsf{DblCat} \to \mathsf{DblCat}$ is our candidate Gray monoidal product on $\mathsf{DblCat}$. Mac Lane's coherence conditions are checked in Section \ref{sec:coherence}.

\subsection{The category of double categories}
\label{sec:DblCat}

We begin with introducing the category $\mathsf{DblCat} $ of double categories and double functors, and recording some of its basic properties.

\begin{definition}
A {\em double category} is an internal category in the category $\mathsf{Cat}$ of categories and functors.
A {\em double functor} in an internal functor in $\mathsf{Cat}$.
Double categories are the objects, and double functors are the morphisms of the category $\mathsf{DblCat}$.
\end{definition}

So a double category consists of 0-cells, also called objects, (interpreted as the objects of the category of objects), vertical 1-cells (which are the morphisms of the category of objects), horizontal 1-cells (the objects of the category of morphisms) and 2-cells (the morphisms of the category of morphisms). 
They can be composed vertically (in the category of objects and the category of morphisms, respectively) and horizontally (via the composition functor of the double category).
As usual in the literature (see e.g. \cite{GrandisPare}), we denote 2-cells as squares surrounded by the appropriate horizontal and vertical source and target 1-cells. We denote by $1$ both horizontal and vertical identity 1-cells; and also identity 2-cells for the horizontal or vertical composition. 
Usually we neither make notational difference between the compositions of horizontal and vertical 1-cells; both are denoted by a dot (if not a diagram is rather drawn).

By \cite[Theorem 4.1]{FiorePaoliPronk} and its proof, $\mathsf{DblCat}$ is locally finitely presentable --- so in particular cocomplete --- and complete. Its terminal object $\mathbbm 1$ is the double category of a single object and only identity higher cells.

Consider the double category $\mathbb G$ which is freely generated by a single 2-cell. In more detail, $\mathbb G$ has four objects, we denote them by $X$, $Y$, $V$ and $Z$. There are identity horizontal and vertical identity 1-cells for each object as well as non-identity horizontal and vertical 1-cells
$$
\xymatrix{X \ar[r]^-t & Y} \qquad
\xymatrix{V \ar[r]^-b & Z} \qquad \textrm{and} \qquad
\raisebox{17pt}{$\xymatrix{X \ar[d]^-l \\  V} \qquad
\xymatrix{Y \ar[d]^-r \\  Z}$}.
$$
There are vertical identity 2-cells at each horizontal 1-cell, horizontal identity 2-cells at each vertical 1-cell, and a single non-identity 2-cell
$$
\xymatrix{X \ar[r]^-t \ar[d]_-l \ar@{}[rd]|-{\Longdownarrow \tau} & 
Y \ar[d]^-r \\
V \ar[r]_-b & 
Z.} 
$$
The functor $\mathsf{DblCat}(\mathbb G,-):\mathsf{DblCat} \to \mathsf{Set}$ sends a double category $\mathbb A$ to the set of double functors $\mathbb G \to \mathbb A$, which can be identified with the set of 2-cells in $\mathbb A$. A double functor $\mathsf F$ is sent to its 2-cell part, which is an isomorphism in $ \mathsf{Set}$ if and only if $\mathsf F$ is bijective on the 2-cells. Since this includes bijectivity also on the identity 2-cells of various kinds, it is equivalent to $\mathsf F$ being bijective on all kinds of cells; that is, its being an isomorphism in $\mathsf{DblCat}$.
By the so obtained conservativity of the functor $\mathsf{DblCat}(\mathbb G,-):\mathsf{DblCat} \to \mathsf{Set}$ we conclude that $\mathbb G$ is a strong generator of the finitely complete category $\mathsf{DblCat}$ with coproducts, see \cite[Proposition 4.5.10]{BorceuxI}.

\subsection{The double categories of double functors}
\label{sec:[A,B]}

Using similar constructions to those in \cite[Section 7]{GrandisPare}, any double categories $\mathbb A$ and $\mathbb B$ determine a double category $\ldb\mathbb A,\mathbb B\rdb$ as follows. 

The \underline{0-cells} are the double functors $\mathbb A \to \mathbb B$. 

The \underline{horizontal 1-cells} are the {\em horizontal pseudotransformations} (called {\em strong horizontal transformations} in \cite[Section 7.4]{GrandisPare}). A horizontal pseudotransformation $x:\mathsf F\to \mathsf G$ consists of the following data.
\begin{itemize}
\item For any 0-cell $A$ of  $\mathbb A$, a horizontal morphism in  $\mathbb B$ on the left; 
\item for any vertical 1-cell $f$ in $\mathbb A$, a 2-cell in $\mathbb B$ in the middle;
\item for any horizontal 1-cell $h$ in $\mathbb A$, a {\em vertically invertible} 2-cell in $\mathbb B$ on the right:
$$
\xymatrix@C=25pt{\mathsf FA \ar[r]^-{x_A} & \mathsf GA}
\qquad \qquad
\xymatrix@C=25pt@R=25pt{
\mathsf FA \ar[r]^-{x_A} \ar[d]_-{\mathsf Ff} \ar@{}[rd]|-{\Longdownarrow{x_f}} & 
\mathsf GA \ar[d]^-{\mathsf Gf} \\
\mathsf FB \ar[r]_-{x_B} &
\mathsf GB}
\qquad \qquad
\xymatrix@C=25pt@R=25pt{
\mathsf FA \ar[r]^-{\mathsf Fh} \ar@{=}[d] \ar@{}[rrd]|-{\Longdownarrow{x^h}} &
\mathsf FC \ar[r]^-{x_C} &
\mathsf GC \ar@{=}[d] \\
\mathsf FA \ar[r]_-{x_A} &
\mathsf GA \ar[r]_-{\mathsf Gh} &
\mathsf GC.}
$$
\end{itemize}
These ingredients are subject to the following axioms.
\begin{enumerate}[(i)]
\item {\em Vertical functoriality}, saying that for the identity vertical 1-cell $1$ on any object $A$ in $\mathbb A$, $x_1$ is equal to the vertical identity 2-cell on the left; and for any composable vertical 1-cells $f$ and $g$ in $\mathbb A$, the equality on the right holds:
$$
\xymatrix@C=25pt@R=62pt{
\mathsf FA \ar[r]^-{x_A} \ar@{=}[d] \ar@{}[rd]|-{\Longdownarrow 1} & 
\mathsf GA \ar@{=}[d] \\
\mathsf FA \ar[r]_-{x_A} &
\mathsf GA}
\qquad \qquad \qquad
\xymatrix{
\mathsf FA \ar[r]^-{x_A} \ar[d]_-{\mathsf Ff} \ar@{}[rd]|-{\Longdownarrow{x_f}} & 
\mathsf GA \ar[d]^-{\mathsf Gf} \\
\mathsf FB \ar[r]|-{\,x_B\,} \ar[d]_-{\mathsf Fg} \ar@{}[rd]|-{\Longdownarrow{x_g}} &
\mathsf GB \ar[d]^-{\mathsf Gg} \\
\mathsf FD \ar[r]_-{x_D} &
\mathsf GD} 
\raisebox{-38pt}{$=$}
\xymatrix@C=28pt@R=28pt{
\mathsf FA \ar[r]^-{x_A} \ar[dd]_-{\mathsf F(g.f)} \ar@{}[rdd]|-{\Longdownarrow{x_{g.f}}} & 
\mathsf GA \ar[dd]^-{\mathsf G(g.f)} \\
\\
\mathsf FD \ar[r]_-{x_D} &
\mathsf GD.}
$$
\item {\em Horizontal functoriality}, saying that for the identity horizontal 1-cell $1$ on any object $A$, $x^1$ is equal to the same vertical identity 2-cell on the left; and for any composable horizontal 1-cells $h$ and $k$ in $\mathbb A$, the equality on the right holds:
$$
\xymatrix@C=25pt@R=62pt{
\mathsf FA \ar[r]^-{x_A} \ar@{=}[d] \ar@{}[rd]|-{\Longdownarrow 1} & 
\mathsf GA \ar@{=}[d] \\
\mathsf FA \ar[r]_-{x_A} &
\mathsf GA}
\qquad \qquad
\xymatrix{
\mathsf FA \ar[r]^-{\mathsf Fh} \ar@{=}[d] \ar@{}[rd]|-{\Longdownarrow 1} &
\mathsf FC \ar[r]^-{\mathsf Fk} \ar@{=}[d] \ar@{}[rrd]|-{\Longdownarrow{x^k}} &
\mathsf FE \ar[r]^-{x_E} &
\mathsf GE \ar@{=}[d] \\
\mathsf FA \ar[r]_-{\mathsf Fh} \ar@{=}[d] \ar@{}[rrd]|-{\Longdownarrow{x^h}} &
\mathsf FC \ar[r]_-{x_C} &
\mathsf GC \ar@{=}[d] \ar[r]_-{\mathsf Gk} \ar@{}[rd]|-{\Longdownarrow 1} &
\mathsf GE \ar@{=}[d] \\
\mathsf FA \ar[r]_-{x_A} &
\mathsf GA \ar[r]_-{\mathsf Gh} &
\mathsf GC \ar[r]_-{\mathsf Gk} &
\mathsf GE}
\raisebox{-38pt}{$=$}
\xymatrix@C=28pt@R=28pt{
\mathsf FA \ar[r]^-{\mathsf F(k.h)} \ar@{=}[dd] \ar@{}[rrdd]|-{\Longdownarrow{x^{k.h}}} &
\mathsf FE \ar[r]^-{x_E} &
\mathsf GE \ar@{=}[dd] \\
\\
\mathsf FA \ar[r]_-{x_A} &
\mathsf GA \ar[r]_-{\mathsf G(k.h)} &
\mathsf GE.}
$$
\item {\em Naturality}, saying that for any 2-cell $\omega$ in $\mathbb A$,
$$
\xymatrix{
\mathsf FA \ar[r]^-{\mathsf Fh} \ar[d]_-{\mathsf Ff} \ar@{}[rd]|-{\Longdownarrow {\mathsf F\omega}} &
\mathsf FC \ar[d]|-{\mathsf Fg} \ar[r]^-{x_C} \ar@{}[rd]|-{\Longdownarrow {x_g}}  &
\mathsf GC \ar[d]^-{\mathsf Gg} \\
\mathsf FB \ar[r]_-{\mathsf Fk} \ar@{=}[d] \ar@{}[rrd]|-{\Longdownarrow{x^k}} &
\mathsf FD \ar[r]_-{x_D} &
\mathsf GD \ar@{=}[d] \\
\mathsf FB \ar[r]_-{x_B} &
\mathsf GB \ar[r]_-{\mathsf Gk} &
\mathsf GD}
\raisebox{-38pt}{$=$}
\xymatrix{
\mathsf FA \ar[r]^-{\mathsf Fh} \ar@{=}[d] \ar@{}[rrd]|-{\Longdownarrow{x^h}} &
\mathsf FC \ar[r]^-{x_C} &
\mathsf GC \ar@{=}[d] \\
\mathsf FA \ar[r]^-{x_A} \ar[d]_-{\mathsf Ff} \ar@{}[rd]|-{\Longdownarrow{x_f}} &
\mathsf GA \ar[r]^-{\mathsf Gh} \ar[d]|-{\mathsf Gf} \ar@{}[rd]|-{\Longdownarrow{\mathsf G\omega}} &
\mathsf GC \ar[d]^-{\mathsf Gg} \\
\mathsf FB \ar[r]_-{x_B} &
\mathsf GB \ar[r]_-{\mathsf Gk} &
\mathsf GD.}
$$
\end{enumerate}

The \underline{vertical 1-cells} are the {\em vertical pseudotransformations} (called {\em strong vertical transformations} in \cite[Section 7.4]{GrandisPare}). A vertical pseudotransformation $y:\mathsf F\to \mathsf H$ consists of the following data.
\begin{itemize}
\item For any 0-cell $A$ of  $\mathbb A$, a vertical 1-cell in  $\mathbb B$ on the left;
\item for any  horizontal 1-cell $h$ in $\mathbb A$, a 2-cell in $\mathbb B$ in the middle;
\item for any vertical 1-cell $f$ in $\mathbb A$, a horizontally invertible 2-cell in $\mathbb B$ on the right:
$$
\xymatrix@R=65pt{\mathsf FA \ar[d]^-{y_A} \\ \mathsf HA}
\qquad \qquad
\xymatrix@C=25pt@R=65pt{
\mathsf FA \ar[d]_-{y_A} \ar[r]^-{\mathsf Fh} \ar@{}[rd]|-{\Longdownarrow{y_h}} & 
\mathsf FC \ar[d]^-{y_C} \\
\mathsf HA \ar[r]_-{\mathsf Hh} &
\mathsf HC}
\qquad \qquad
\xymatrix@C=25pt@R=25pt{
\mathsf FA \ar[d]_-{y_A} \ar@{=}[r] \ar@{}[rd]|-{\Longdownarrow{y^f}} &
\mathsf FA \ar[d]^-{\mathsf Ff} \\
\mathsf HA \ar[d]_-{\mathsf Hf} &
\mathsf FB \ar[d]^-{y^B} \\
\mathsf HB \ar@{=}[r] &
\mathsf HB.}
$$
\end{itemize}
These ingredients are subject to the following axioms.
\begin{enumerate}[(i)]
\item {\em Horizontal functoriality}, saying that for the identity horizontal 1-cell $1$ on any object $A$ in $\mathbb A$, $y_1$ is equal to the horizontal identity 2-cell on the left; 
and for any composable horizontal 1-cells $h$ and $k$ in $\mathbb A$, the equality on the right holds:
$$
\xymatrix@C=25pt@R=23pt{
\mathsf FA \ar[d]_-{y_A} \ar@{=}[r] \ar@{}[rd]|-{\Longdownarrow 1} & 
\mathsf FA \ar[d]^-{y_A} \\
\mathsf HA \ar@{=}[r] &
\mathsf HA}
\qquad \qquad
\xymatrix{
\mathsf FA \ar[d]_-{y_A} \ar[r]^-{\mathsf Fh} \ar@{}[rd]|-{\Longdownarrow{y_h}} & 
\mathsf FC \ar[d]|-{y_C}  \ar[r]^-{\mathsf Fk} \ar@{}[rd]|-{\Longdownarrow{y_k}} &
\mathsf FE \ar[d]^-{y_E} \\
\mathsf HA \ar[r]_-{\mathsf Hh} &
\mathsf HC \ar[r]_-{\mathsf Hk} &
\mathsf HE} 
\raisebox{-18pt}{$=$}
\xymatrix{
\mathsf FA \ar[d]_-{y_A} \ar[rr]^-{\mathsf F(k.h)} \ar@{}[rrd]|-{\Longdownarrow{y_{k.h}}} && 
\mathsf FE \ar[d]^-{y_E} \\
\mathsf HA \ar[rr]_-{\mathsf H(k.h)} &&
\mathsf HE.}
$$
\item {\em Vertical functoriality}, saying that for the identity vertical 1-cell $1$ on any object $A$, $y^1$ is equal to the same horizontal identity 2-cell on the left; and for any composable vertical 1-cells $f$ and $g$ in $\mathbb A$, the equality on the right holds:
$$
\xymatrix@C=25pt@R=102pt{
\mathsf FA \ar[d]_-{y_A} \ar@{=}[r] \ar@{}[rd]|-{\Longdownarrow 1} & 
\mathsf FA \ar[d]^-{y_A} \\
\mathsf HA \ar@{=}[r] &
\mathsf HA}
\qquad \qquad
\xymatrix{
\mathsf FA \ar[d]_-{y_A} \ar@{=}[r] \ar@{}[rdd]|-{\Longdownarrow{y^f}} &
\mathsf FA \ar[d]_-{\mathsf Ff} \ar@{=}[r] \ar@{}[rd]|-{\Longdownarrow 1} &
\mathsf FA \ar[d]^-{\mathsf Ff}  \\
\mathsf HA \ar[d]_-{\mathsf Hf} &
\mathsf FB \ar[d]^-{y_B} \ar@{=}[r]  \ar@{}[rdd]|-{\Longdownarrow{y^g}} &
\mathsf FB \ar[d]^-{\mathsf Fg}  \\
\mathsf HB \ar[d]_-{\mathsf Hg} \ar@{=}[r] \ar@{}[rd]|-{\Longdownarrow 1} &
\mathsf HB \ar[d]^-{\mathsf Hg} &
\mathsf FD \ar[d]^-{y_D}  \\
\mathsf HD \ar@{=}[r] &
\mathsf HD \ar@{=}[r] &
\mathsf HD}
\raisebox{-58pt}{$\ =$}
\xymatrix@C=25pt@R=24pt{
\mathsf FA \ar[d]_-{y_A} \ar@{=}[r] \ar@{}[rddd]|-{\Longdownarrow{y^{g.f}}} &
\mathsf FA \ar[dd]^-{\mathsf F(g.f)}  \\
\mathsf HA \ar[dd]_-{\mathsf H(g.f)} \\
& \mathsf FD \ar[d]^-{y_D} \\
\mathsf HD \ar@{=}[r] &
\mathsf HD.}
$$
\item {\em Naturality}, saying that for any 2-cell $\omega$ in $\mathbb A$,
$$
\xymatrix{
\mathsf FA \ar[d]_-{y_A} \ar@{=}[r] \ar@{}[rdd]|-{\Longdownarrow {y^f}}  &
\mathsf FA \ar[r]^-{\mathsf Fh} \ar[d]_-{\mathsf Ff} \ar@{}[rd]|-{\Longdownarrow {\mathsf F\omega}} &
\mathsf FC \ar[d]^-{\mathsf Fg}  \\
\mathsf HA \ar[d]_-{\mathsf Hf} &
\mathsf FB \ar[r]^-{\mathsf Fk} \ar[d]_-{y_B} \ar@{}[rd]|-{\Longdownarrow{y_k}} &
\mathsf FD \ar[d]^-{y_D} \\
\mathsf HB \ar@{=}[r] &
\mathsf HB \ar[r]_-{\mathsf Hk} &
\mathsf HD}
\raisebox{-38pt}{$=$}
\xymatrix{
\mathsf FA \ar[r]^-{\mathsf Fh} \ar[d]_-{y_A} \ar@{}[rd]|-{\Longdownarrow{y_h}} &
\mathsf FC \ar[d]^-{y_C} \ar@{=}[r] \ar@{}[rdd]|-{\Longdownarrow{y^g}} &
\mathsf FA \ar[d]^-{\mathsf Fg} \\
\mathsf HA \ar[r]^-{\mathsf Hh} \ar[d]_-{\mathsf Hf} \ar@{}[rd]|-{\Longdownarrow{\mathsf H\omega}} &
\mathsf HC \ar[d]^-{\mathsf Hg} &
\mathsf FD \ar[d]^-{y_D} \\
\mathsf HB \ar[r]_-{\mathsf Hk} &
\mathsf HD \ar@{=}[r] &
\mathsf HD.}
$$
\end{enumerate}

The \underline{2-cells} are the {\em modifications}. A modification on the left is given by a collection of 2-cells in $\mathbb B$ on the right, for all 0-cells $A$ of $\mathbb A$:
$$
\xymatrix@C=25pt@R=25pt{
\mathsf F \ar[r]^-{x} \ar[d]_-y \ar@{}[rd]|-{\Longdownarrow \Theta} &
\mathsf G \ar[d]^-v \\
\mathsf H \ar[r]_-z &
\mathsf K}
\qquad \qquad \qquad
\xymatrix@C=25pt@R=25pt{
\mathsf FA \ar[r]^-{x_A} \ar[d]_-{y_A} \ar@{}[rd]|-{\Longdownarrow {\Theta_A}} &
\mathsf GA \ar[d]^-{v_A} \\
\mathsf HA \ar[r]_-{z_A} &
\mathsf KA}
$$
satisfying the following axioms.
\begin{enumerate}[(i)]
\item For any horizontal 1-cell $h$ in $\mathbb A$,
$$
\xymatrix{
\mathsf FA \ar[r]^-{\mathsf Fh} \ar[d]_-{y_A} \ar@{}[rd]|-{\Longdownarrow {y_h}} &
\mathsf FC \ar[r]^-{x_C} \ar[d]|-{y_C} \ar@{}[rd]|-{\Longdownarrow {\Theta_C}} &
\mathsf GC \ar[d]^-{v_C} \\
\mathsf HA \ar[r]_-{\mathsf Hh} \ar@{=}[d] \ar@{}[rrd]|-{\Longdownarrow {z^h}} &
\mathsf HC \ar[r]_-{z_C} &
\mathsf KC \ar@{=}[d] \\
\mathsf HA \ar[r]_-{z_A} &
\mathsf KA \ar[r]_-{\mathsf Kh}  &
\mathsf KC}
\raisebox{-38pt}{$=$}
\xymatrix{
\mathsf FA \ar[r]^-{\mathsf Fh} \ar@{=}[d] \ar@{}[rrd]|-{\Longdownarrow {x^h}} &
\mathsf FC \ar[r]^-{x_C} &
\mathsf GC \ar@{=}[d] \\
\mathsf FA \ar[r]^-{x_A} \ar[d]_-{y_A} \ar@{}[rd]|-{\Longdownarrow {\Theta_A}}  &
\mathsf GA \ar[r]^-{\mathsf Gh} \ar[d]|-{v_A} \ar@{}[rd]|-{\Longdownarrow {v_h}} &
\mathsf GC \ar[d]^-{v_C} \\
\mathsf HA \ar[r]_-{z_A} &
\mathsf KA \ar[r]_-{\mathsf Kh}  &
\mathsf KC.}
$$
\item For any vertical 1-cell $f$ in $\mathbb A$,
$$
\xymatrix{
\mathsf FA \ar[r]^-{x_A} \ar[d]_-{y_A} \ar@{}[rd]|-{\Longdownarrow {\Theta_A}} &
\mathsf GA \ar@{=}[r] \ar[d]^-{v_A} \ar@{}[rdd]|-{\Longdownarrow {v^f}} &
\mathsf GA \ar[d]^-{\mathsf Gf} \\
\mathsf HA \ar[r]^-{z_A} \ar[d]_-{\mathsf Hf} \ar@{}[rd]|-{\Longdownarrow {z_f}} &
\mathsf KA \ar[d]^-{\mathsf Kf} &
\mathsf GB \ar[d]^-{v_B} \\
\mathsf HB \ar[r]_-{z_B} &
\mathsf KB \ar@{=}[r]  &
\mathsf KB}
\raisebox{-38pt}{$=$}
\xymatrix{
\mathsf FA \ar[d]_-{y_A} \ar@{=}[r] \ar@{}[rdd]|-{\Longdownarrow {y^f}} &
\mathsf FA \ar[r]^-{x_A} \ar[d]_-{\mathsf Ff}  \ar@{}[rd]|-{\Longdownarrow {x_f}} &
\mathsf GA \ar[d]^-{\mathsf Gf} \\
\mathsf HA \ar[d]_-{\mathsf Hf} &
\mathsf FB \ar[r]^-{x_B} \ar[d]_-{y_B} \ar@{}[rd]|-{\Longdownarrow {\Theta_B}} &
\mathsf GB \ar[d]^-{v_B} \\
\mathsf HB \ar@{=}[r] &
\mathsf HB \ar[r]_-{z_B}  &
\mathsf KB.}
$$
\end{enumerate}

The identity horizontal pseudotransformation has the components
$$
\xymatrix{
\mathsf  F A \ar@{=}[r] \ar[d]_-{\mathsf F f} \ar@{}[rd]|-{\Longdownarrow 1} &
\mathsf F A \ar[d]^-{\mathsf F f} \\
\mathsf F B \ar@{=}[r] & 
\mathsf B}
\qquad \qquad
\xymatrix{
\mathsf F A \ar[r]^-{\mathsf F h} \ar@{=}[d] \ar@{}[rrd]|-{\Longdownarrow 1} &
\mathsf F C \ar@{=}[r] &
\mathsf F C \ar@{=}[d] \\
\mathsf F A \ar@{=}[r] &
\mathsf F A \ar[r]_-{\mathsf F h} &
\mathsf F C;}
$$
while the composite of some horizontal pseudotransformations 
$\xymatrix@C=15pt{\mathsf F \ar[r]^-x & \mathsf G \ar[r]^-z & \mathsf H}$ 
has the components
$$
\xymatrix@R=28pt{
\mathsf F A\ar[r]^-{x_A} \ar[dd]_-{\mathsf Ff} \ar@{}[rdd]|-{\Longdownarrow{x_f}} &
\mathsf GA \ar[r]^-{z_A} \ar[dd]|-{\mathsf Gf} \ar@{}[rdd]|-{\Longdownarrow{z_f}} &
\mathsf HA \ar[dd]^-{\mathsf H f} \\
\\
\mathsf F B \ar[r]_-{x_B} &
\mathsf GB \ar[r]_-{z_B} &
\mathsf H B}
\qquad \qquad
\xymatrix{
\mathsf F A \ar[r]^-{\mathsf Fh} \ar@{=}[d] \ar@{}[rrd]|-{\Longdownarrow {x^h}} &
\mathsf F C \ar[r]^-{x_C} &
\mathsf G C \ar[r]^-{z_C} \ar@{=}[d] \ar@{}[rd]|-{\Longdownarrow 1} &
\mathsf H C \ar@{=}[d] \\
\mathsf F A\ar[r]^-{x_A} \ar@{=}[d] \ar@{}[rd]|-{\Longdownarrow 1} &
\mathsf G A\ar[r]^-{\mathsf Gh} \ar@{=}[d]  \ar@{}[rrd]|-{\Longdownarrow {z^h}} &
\mathsf G C \ar[r]_-{z_C} &
\mathsf H C \ar@{=}[d] \\
\mathsf F A \ar[r]_-{x_A} &
\mathsf GA \ar[r]_-{z_A} &
\mathsf H A \ar[r]_-{\mathsf Hh} &
\mathsf H C}
$$
for any horizontal 1-cell $h$ and vertical 1-cell $f$ in $\mathbb A$.
Symmetric formulae apply to the vertical pseudotransformations.
The components of the horizontal composite of modifications are the horizontal composites of their components, and the components of the vertical composite of modifications are the vertical composites of their components.

Throughout, we identify any double category $\mathbb A$ with the isomorphic double category $\ldb \mathbbm 1,\mathbb A \rdb$.

\subsection{The functor $\ldb -,- \rdb:\mathsf{DblCat}^{\mathsf{op}} \times \mathsf{DblCat}\to \mathsf{DblCat}$}
\label{sec:[-,-]}

In this section we interpret the map, sending a pair of double categories $\mathbb A$ and $\mathbb B$ to the double category $\ldb \mathbb A,\mathbb B \rdb$ of Section \ref{sec:[A,B]}, as the object map of a functor in the title. So we need to construct its morphism map, sending a pair of double functors $\mathsf F:\mathbb A'\to \mathbb A$ and $\mathsf G:\mathbb B\to \mathbb B'$ to a double functor $\ldb \mathsf F,\mathsf G \rdb:\ldb \mathbb A,\mathbb B \rdb \to \ldb \mathbb A',\mathbb B' \rdb$.

Its \underline{0-cell} part sends a double functor 
$\xymatrix@C=13pt{
\mathbb A \ar[r]^-{\mathsf H} &
\mathbb B}$ 
to the composite 
$\xymatrix@C=13pt{
\mathbb A' \ar[r]^-{\mathsf F} &
\mathbb A \ar[r]^-{\mathsf H} &
\mathbb B \ar[r]^-{\mathsf G} &
\mathbb B'.}$

The \underline{horizontal 1-cell} part sends a horizontal pseudotransformation 
$\xymatrix@C=15pt{
\mathsf H \ar[r]^-x &
\mathsf H'}$
to the horizontal pseudotransformation with the components
$$
\xymatrix{
\mathsf{GHF} A \ar[r]^-{\mathsf  G x_{\mathsf FA}} \ar[d]_-{\mathsf{GHF} f}
\ar@{}[rd]|-{\Longdownarrow{\mathsf  G x_{\mathsf Ff}}} &
\mathsf{GH'F} A \ar[d]^-{\mathsf{GH'F} f} \\
\mathsf{GHF} B \ar[r]_-{\mathsf  G x_{\mathsf FB}} &
\mathsf{GH'F}B}
\qquad \qquad
\xymatrix{
\mathsf{GHF} A \ar[r]^-{\mathsf{GHF} h} \ar@{=}[d] 
\ar@{}[rrd]|-{\Longdownarrow{\mathsf  G x^{\mathsf Fh}}} &
\mathsf{GHF}  C \ar[r]^-{\mathsf  G x_{\mathsf FC}} &
\mathsf{GH'F}  C \ar@{=}[d] \\
\mathsf{GHF} A \ar[r]_-{\mathsf  G x_{\mathsf FA}} &
\mathsf{GH'F} A \ar[r]_-{\mathsf{GH'F} h} &
\mathsf{GH'F}  C}
$$
for any horizontal 1-cell $h$ and vertical 1-cell $f$ in $\mathbb A'$.

Symmetrically, the \underline{vertical 1-cell} part sends a vertical pseudotransformation $y$ to the vertical pseudotransformation with the components
$$
\xymatrix@R=27pt{
\mathsf{GHF} A \ar[r]^-{\mathsf{GHF} h} \ar[dd]_-{\mathsf  G y_{\mathsf FA}} 
\ar@{}[rdd]|-{\Longdownarrow{\mathsf  G y_{\mathsf Fh}}} &
\mathsf{GHF}  C \ar[dd]^-{\mathsf  G y_{\mathsf FC}} \\
\\
\mathsf{GH^{\prime \prime}F} A \ar[r]_-{\mathsf{GH^{\prime \prime}F} h} &
\mathsf{GH^{\prime \prime}F}  C}
\qquad \qquad
\xymatrix{
\mathsf{GHF} A \ar[d]_-{\mathsf  G y_{\mathsf FA}} \ar@{=}[r] 
\ar@{}[rd]|-{\Longdownarrow{\mathsf  G y^{\mathsf Ff}}} &
\mathsf{GHF}  A \ar[d]^-{\mathsf{GHF} f}  \\
\mathsf{GH^{\prime \prime}F} A \ar[d]_-{\mathsf{GH^{\prime \prime}F} f} &
\mathsf{GHF} B \ar[d]^-{\mathsf  G y_{\mathsf FB}} \\
\mathsf{GH^{\prime \prime}F}  B \ar@{=}[r]  &
\mathsf{GH^{\prime \prime}F}  B }
$$
for any horizontal 1-cell $h$ and vertical 1-cell $f$ in $\mathbb A'$.

Finally, the \underline{2-cell} part sends a modification in the first diagram to the modification with components in the second diagram:
$$
\xymatrix{
\mathsf H \ar[r]^-{x} \ar[d]_-y \ar@{}[rd]|-{\Longdownarrow \Theta} &
\mathsf H' \ar[d]^-v \\
\mathsf H^{\prime \prime} \ar[r]_-z &
\mathsf H^{\prime \prime \prime }}
\qquad
\xymatrix{
\mathsf {GHF}A \ar[r]^-{\mathsf Gx_{\mathsf FA}} \ar[d]_-{\mathsf Gy_{\mathsf FA}} 
\ar@{}[rd]|-{\Longdownarrow {\mathsf G\Theta_{\mathsf FA}}} &
\mathsf{GH'F}A \ar[d]^-{\mathsf Gv_{\mathsf FA}} \\
\mathsf {GH^{\prime \prime} F}A \ar[r]_-{\mathsf Gz_{\mathsf FA}} &
\mathsf {GH^{\prime \prime \prime }F}A.}
$$

\subsection{The extranatural transformation $\mathfrak l$ }
\label{sec:l}

In this section we construct an extranatural transformation 
$$
\xymatrix@C=60pt{
\mathsf{DblCat}^{\mathsf{op}} \times \mathsf{DblCat} \times 
\mathsf{DblCat} \times \mathsf{DblCat}^{\mathsf{op}}
\ar[r]^-{\ldb -,- \rdb \times ! \times !} 
\ar[d]_(.07){\xymatrix@C=3pt@R=17pt{
\ar[rrd] &
\ar[rrrrd] &&&
\ar[lllld] &
\ar[lld] \\
&&&&&}\hspace{2.3cm}}
\ar@{}[rd]|-{\longdownarrow{\displaystyle {\mathfrak l}}} &
\mathsf{DblCat} \\
\mathsf{DblCat} \times \mathsf{DblCat}^{\mathsf{op}} \times 
\mathsf{DblCat}^{\mathsf{op}} \times \mathsf{DblCat}
\ar[r]_-{\ldb \mathbb -,- \rdb^{\mathsf{op}}\times \ldb \mathbb -,- \rdb} &
\mathsf{DblCat}^{\mathsf{op}} \times \mathsf{DblCat} 
\ar[u]_-{\ldb -,- \rdb} }
$$
where $!$ in the top row denotes the unique functor to the terminal category, and in the left column the depicted symmetry natural isomorphism --- that is, the appropriate flip map --- occurs. We denote by lower indices that $\mathfrak l$ is {\em ordinary natural} in the first two arguments, and an upper index reminds us that it is {\em extranatural} in the last two factors.
At any object of the form $\mathbb A,\mathbb B,\mathbb D,\mathbb D$, it is given by the following double functor 
$\mathfrak l^{\mathbb D}_{\mathbb A,\mathbb B}:\ldb \mathbb A,\mathbb B \rdb\to 
\ldb \ldb \mathbb D,\mathbb A \rdb,\ldb \mathbb D, \mathbb B \rdb  \rdb$.

The \underline {0-cell} part sends a double functor $\mathsf G:\mathbb A \to \mathbb B$ to the double functor $\ldb1,\mathsf G \rdb:\ldb \mathbb D,\mathbb A \rdb \to \ldb \mathbb D, \mathbb B \rdb$ of Section \ref{sec:[-,-]}.

The \underline {horizontal 1-cell} part sends a horizontal pseudotransformation 
$\xymatrix@C=12pt{\mathsf G \ar[r]^-x & \mathsf {G'}}$
to the horizontal pseudotransformation $\ldb 1,\mathsf G \rdb \to \ldb 1,\mathsf G' \rdb$ with the following components, for any horizontal pseudotransformation $p$ and vertical pseudotransformation $q$ between double functors $\mathsf H ,\mathsf H':\mathbb D \to \mathbb A$.
\begin{itemize}
\item The horizontal pseudotransformation with components
$$
\xymatrix{
\mathsf{GH}A \ar[r]^-{x_{\mathsf H A}} \ar[d]_-{\mathsf{GH}f} 
\ar@{}[rd]|-{\Longdownarrow {x_{\mathsf H f}}} &
\mathsf{G'H}A \ar[d]^-{\mathsf{G'H}f}  \\
\mathsf{GH}B \ar[r]_-{x_{\mathsf H B}} &
\mathsf{G'H}B}
\qquad \qquad
\xymatrix{
\mathsf{GH}A \ar[r]^-{\mathsf{GH} h} \ar@{=}[d] \ar@{}[rrd]|-{\Longdownarrow {x^{\mathsf H h}}} &
\mathsf{GH}C \ar[r]^-{x_{\mathsf H C}} &
\mathsf{G'H}C \ar@{=}[d] \\
\mathsf{GH}A \ar[r]_-{x_{\mathsf H A}} &
\mathsf{G'H}A \ar[r]_-{\mathsf{G'H} h} &
\mathsf{G'H}C}
$$
for any horizontal 1-cell $h$ and vertical 1-cell $f$ in $\mathbb D$.
\item The modification with components
$$
\xymatrix{
\mathsf{GH}A \ar[r]^-{x_{\mathsf H A}} \ar[d]_-{\mathsf G q_A}
\ar@{}[rd]|-{\Longdownarrow {x_{q_A}}} &
\mathsf{G'H}A \ar[d]^-{\mathsf G' q_A} \\
\mathsf{GH'}A \ar[r]_-{x_{\mathsf H' A}}  &
\mathsf{G'H'}A}
$$
for any 0-cell $A$ in $\mathbb D$.
\item The vertically invertible modification with components
$$
\xymatrix{
\mathsf{GH}A \ar[r]^-{\mathsf G p_A} \ar@{=}[d] \ar@{}[rrd]|-{\Longdownarrow {x^{p_A}}} &
\mathsf{GH'}A \ar[r]^-{x_{\mathsf H' A}}  &
\mathsf{G'H'}A \ar@{=}[d] \\
\mathsf{GH}A \ar[r]_-{x_{\mathsf H A}}  &
\mathsf{G'H}A \ar[r]_-{\mathsf G' p_A} &
\mathsf{G'H'}A}
$$
for any 0-cell $A$ in $\mathbb D$.
\end{itemize}

Symmetrically, the \underline {vertical 1-cell} part sends a vertical pseudotransformation $y$ from $\mathsf G$ to $\mathsf G'$ to the vertical pseudotransformation from $\ldb 1,\mathsf G\rdb$ to $\ldb 1,\mathsf G'\rdb$ with the following components, for any horizontal pseudotransformation $p$ and vertical pseudotransformation $q$ between double functors $\mathsf H ,\mathsf H':\mathbb D \to \mathbb A$.
\begin{itemize}
\item The vertical pseudotransformation with components
$$
\xymatrix@R=27pt{
\mathsf{GH}A \ar[r]^-{\mathsf{GH} h} \ar[dd]_-{y_{\mathsf H A}} 
\ar@{}[rd]|-{\Longdownarrow {y_{\mathsf H h}}} &
\mathsf{GH}C \ar[dd]^-{y_{\mathsf H C}} \\
&\\
\mathsf{G'H}A \ar[r]_-{\mathsf{G'H} h} &
\mathsf{G'H}C}
\qquad \qquad
\xymatrix{
\mathsf{GH}A  \ar[d]_-{y_{\mathsf H A}} \ar@{=}[r] \ar@{}[rd]|-{\Longdownarrow {y^{\mathsf H f}}}  &
\mathsf{GH}A \ar[d]^-{\mathsf{GH} f} \\
\mathsf{G'H}A \ar[d]_-{\mathsf{G'H} f} &
\mathsf{GH} B \ar[d]^-{y_{\mathsf H B}}  \\
\mathsf{G'H}B \ar@{=}[r] &
\mathsf{G'H}B}
$$
for any horizontal 1-cell $h$ and vertical 1-cell $f$ in $\mathbb D$.
\item The modification with components
$$
\xymatrix{
\mathsf{GH}A  \ar[d]_-{y_{\mathsf H A}} \ar[r]^-{\mathsf G p_A}  
\ar@{}[rd]|-{\Longdownarrow {y_{p_A}}} &
\mathsf{GH'}A  \ar[d]^-{y_{\mathsf {H'} A}} \\
\mathsf{G'H}A \ar[r]_-{\mathsf{G'} p_A} &
\mathsf{G'H}A}
$$
for any 0-cell $A$ in $\mathbb D$.
\item The horizontally invertible modification with components
$$
\xymatrix{
\mathsf{GH}A \ar@{=}[r]  \ar[d]_-{y_{\mathsf H A}} 
\ar@{}[rd]|-{\Longdownarrow {y^{q_A}}}  &
\mathsf{GH}A \ar[d]^-{\mathsf{G} q_A} \\
\mathsf{G'H}A \ar[d]_-{\mathsf{G'} q_A} &
\mathsf{GH'}A \ar[d]^-{y_{\mathsf {H'} A}} \\
\mathsf{G'H'}A \ar@{=}[r] &
\mathsf{G'H'}A}
$$
for any 0-cell $A$ in $\mathbb D$.
\end{itemize}

The \underline{2-cell} part sends a modification in the first diagram to the modification with components in the second diagram, for any 0-cell $A$ in $\mathbb D$:
$$
\xymatrix{
\mathsf G \ar[r]^-x \ar[d]_-y \ar@{}[rd]|-{\Longdownarrow \Gamma} &
\mathsf{G'} \ar[d]^-v \\
\mathsf{G^{\prime \prime}} \ar[r]_-z &
\mathsf{G^{\prime \prime \prime }}}
\qquad \qquad
\xymatrix{
\mathsf {GH} A \ar[r]^-{x_{\mathsf H  A}} \ar[d]_-{y_{\mathsf H  A}} 
\ar@{}[rd]|-{\Longdownarrow {\Gamma_{\mathsf H A}}} &
\mathsf{G'H} A \ar[d]^-{v_{\mathsf H A}} \\
\mathsf{G^{\prime \prime}H} A \ar[r]_-{z_{\mathsf H A}} &
\mathsf{G^{\prime \prime \prime}H}  A.}
$$

For any double categories $\mathbb A$, $\mathbb B$, $\mathbb C$ and $\mathbb D$, direct computation verifies the commutativity of
\begin{equation} \label{eq:l_comm}
\xymatrix@C=18pt{
\ldb \mathbb A,\mathbb B \rdb
\ar[rr]^-{\mathfrak l^{\mathbb D}_{\mathbb A,\mathbb B}}
\ar[d]_-{\mathfrak l^{\mathbb C}_{\mathbb A,\mathbb B}} &&
\ldb \ldb \mathbb D,\mathbb A \rdb,\ldb \mathbb D,\mathbb B \rdb \rdb
\ar[d]^-{\ldb 1,\mathfrak l^{\mathbb C}_{\mathbb D,\mathbb B} \rdb} \\
\ldb \ldb \mathbb C,\mathbb A \rdb,\ldb \mathbb C,\mathbb B \rdb \rdb
\ar[r]_-{\raisebox{-10pt}{${}_{
\mathfrak l^{\ldb \mathbb C,\mathbb D \rdb}_{\ldb \mathbb C,\mathbb A \rdb,\ldb \mathbb C,\mathbb B \rdb}}$}} &
\ldb
\ldb \ldb \mathbb C,\mathbb D \rdb,\ldb \mathbb C,\mathbb A \rdb \rdb,
\ldb \ldb \mathbb C,\mathbb D \rdb,\ldb \mathbb C,\mathbb B \rdb \rdb
\rdb
\ar[r]_-{\raisebox{-10pt}{${}_{
\ldb
\mathfrak l^{\mathbb C}_{\mathbb D ,\mathbb A },1
\rdb}$}} &
\ldb 
\ldb \mathbb D,\mathbb A \rdb,
\ldb \ldb \mathbb C,\mathbb D \rdb,\ldb \mathbb C,\mathbb B \rdb \rdb
\rdb}
\end{equation}
and the equality of
\begin{equation} \label{eq:l_id}
\xymatrix{
\ldb \mathbb A,\mathbb B \rdb
\ar[r]^-{\mathfrak l^{\mathbb A}_{\mathbb A,\mathbb B}} &
\ldb \ldb \mathbb A,\mathbb A \rdb,\ldb \mathbb A,\mathbb B \rdb \rdb
\ar[r]^-{\ldb 1_{\mathbb A},1 \rdb} &
\ldb\mathbbm 1, \ldb \mathbb A,\mathbb B \rdb  \rdb \cong
\ldb \mathbb A,\mathbb B \rdb}
\end{equation}
to the identity double functor, where $1_{\mathbb A}:\mathbbm 1 \to \ldb \mathbb A,\mathbb A \rdb$ is the double functor sending the single object of $\mathbbm 1$ to the identity double functor $1_{\mathbb A}:\mathbb A \to \mathbb A$.

\subsection{The extranatural transformation $\mathfrak r$}
\label{sec:r}

In this section we construct another extranatural transformation
$$
\xymatrix@C=60pt{
\mathsf{DblCat}^{\mathsf{op}} \times \mathsf{DblCat} \times 
\mathsf{DblCat} \times \mathsf{DblCat}^{\mathsf{op}}
\ar[r]^-{\ldb-,-\rdb \times !\times !} 
\ar[d]_(.07){
\xymatrix@C=3pt@R=17pt{
\ar[rrrrd] &&
\ar[lld] &
\ar[rrd] &&
\ar[lllld] \\
&&&&&}\hspace{2.3cm}}
\ar@{}[rd]|-{\longdownarrow {\mathfrak r}} &
\mathsf{DblCat} \\
\mathsf{DblCat} \times  \mathsf{DblCat}^{\mathsf{op}} \times 
\mathsf{DblCat}^{\mathsf{op}} \times \mathsf{DblCat}
\ar[r]_-{\ldb-,-\rdb^{\mathsf{op}} \times \ldb-,-\rdb} &
\mathsf{DblCat}^{\mathsf{op}} \times \mathsf{DblCat}
\ar[u]_-{\ldb-,-\rdb}}
$$
where again, $!$ in the top row denotes the unique functor to the terminal category, and in the left column the depicted symmetry natural isomorphism --- that is, the appropriate flip map --- occurs.
As in Section \ref{sec:l}, we denote by lower indices that $\mathfrak r$ is {\em ordinary natural} in the first two arguments, and an upper index reminds us that it is {\em extranatural} in the last two factors.
At any object of the form $\mathbb A,\mathbb B,\mathbb D,\mathbb D$, it is given by the following double functor 
$\mathfrak r^{\mathbb D}_{\mathbb A,\mathbb B}:\ldb \mathbb A,\mathbb B \rdb\to 
\ldb \ldb \mathbb B,\mathbb D \rdb,\ldb \mathbb A,\mathbb D \rdb \rdb$.

The \underline {0-cell} part sends a double functor $\mathsf F:\mathbb A \to \mathbb B$ to the double functor $\ldb \mathsf F,1\rdb:\ldb \mathbb B,\mathbb D \rdb \to \ldb \mathbb A,\mathbb D \rdb$ in Section \ref{sec:[-,-]}.

The \underline {horizontal 1-cell} part sends a horizontal pseudotransformation $\xymatrix@C=12pt{\mathsf F \ar[r]^-x & \mathsf {F'}}$
to the horizontal pseudotransformation $\ldb \mathsf F,1\rdb \to \ldb \mathsf F',1\rdb$ with the following components, for any horizontal pseudotransformation $p$ and vertical pseudotransformation $q$ between double functors $\mathsf H ,\mathsf H':\mathbb B \to \mathbb D$.
\begin{itemize}
\item The horizontal pseudotransformation with components
$$
\xymatrix{
\mathsf {HF}A \ar[r]^-{\mathsf Hx_A} \ar[d]_-{\mathsf {HF} f} 
\ar@{}[rd]|-{\Longdownarrow {\mathsf Hx_f}} &
\mathsf {HF'}A \ar[d]^-{\mathsf {HF'} f} \\
\mathsf {HF}B \ar[r]_-{\mathsf Hx_B} &
\mathsf {HF'}B}
\qquad \qquad
\xymatrix{
\mathsf {HF}A \ar[r]^-{\mathsf {HF} h} \ar@{=}[d] \ar@{}[rrd]|-{\Longdownarrow{\mathsf Hx^h}} &
\mathsf {HF} C \ar[r]^-{\mathsf Hx_C} &
\mathsf {HF'} C \ar@{=}[d]  \\
\mathsf {HF}A \ar[r]_-{\mathsf Hx_A} &
\mathsf {HF'}A \ar[r]_-{\mathsf {HF'} h} &
\mathsf {HF'} C}
$$
for horizontal 1-cells $h$ and vertical 1-cells $f$ in $\mathbb A$.
\item The modification with the components
$$
\xymatrix{
\mathsf {HF}A \ar[r]^-{\mathsf Hx_A} \ar[d]_-{q_{\mathsf F A}} \ar@{}[rd]|-{\Longdownarrow{q_{x_A}}} &
\mathsf {HF'}A \ar[d]^-{q_{\mathsf F' A}} \\
\mathsf {H'F}A \ar[r]_-{\mathsf H'x_A} &
\mathsf {H'F'}A}
$$
for any 0-cell $A$ in $\mathbb A$.
\item The vertically invertible modification with the components
$$
\xymatrix{
\mathsf{HF} A \ar[r]^-{p_{\mathsf F A}} \ar@{=}[d] \ar@{}[rrd]|-{\Longdownarrow{(p^{x_A})^{-1}}} &
\mathsf{H'F} A \ar[r]^-{\mathsf{H'} x_A} &
\mathsf{H'F'} A \ar@{=}[d] \\
\mathsf{HF} A  \ar[r]_-{\mathsf H x_A} &
\mathsf{HF'} A \ar[r]_-{p_{\mathsf{F'} A}} &
\mathsf{H'F'} A}
$$
for any 0-cell $A$ in $\mathbb A$.
\end{itemize}

Symmetrically, the \underline {vertical 1-cell} part sends a vertical pseudotransformation $y$ from $\mathsf F$ to $\mathsf F'$ to the vertical pseudotransformation $\ldb \mathsf F,1\rdb$ to $\ldb \mathsf F',1\rdb$ with the following components, for any horizontal pseudotransformation $p$ and vertical pseudotransformation $q$ between double functors $\mathsf H ,\mathsf H':\mathbb B \to \mathbb D$.
\begin{itemize}
\item The vertical pseudotransformation with components
$$
\xymatrix@R=27pt{
\mathsf {HF}A \ar[dd]_-{\mathsf Hy_A} \ar[r]^-{\mathsf {HF} h} 
\ar@{}[rd]|-{\Longdownarrow {\mathsf Hy_h}} &
\mathsf {HF}C \ar[dd]^-{\mathsf Hy_C} \\
& \\
\mathsf {HF'}A \ar[r]_-{\mathsf {HF'} h} &
\mathsf {HF'}C}
\qquad \qquad
\xymatrix{
\mathsf {HF}A \ar[d]_-{\mathsf Hy_A} \ar@{=}[r] 
\ar@{}[rd]|-{\Longdownarrow{\mathsf Hy^f}} &
\mathsf {HF}A \ar[d]^-{\mathsf {HF} f} \\
\mathsf {HF'} A \ar[d]_-{\mathsf {HF'} f} &
\mathsf {HF}B \ar[d]^-{\mathsf Hy_B} \\
\mathsf {HF'} B \ar@{=}[r]  &
\mathsf {HF'} B}
$$
for horizontal 1-cells $h$ and vertical 1-cells $f$ in $\mathbb A$.
\item The modification with the components
$$
\xymatrix{
\mathsf {HF}A \ar[d]_-{\mathsf Hy_A} \ar[r]^-{p_{\mathsf F A}} \ar@{}[rd]|-{\Longdownarrow{p_{y_A}}} &
\mathsf {H'F}A \ar[d]^-{\mathsf H'y_A} \\
\mathsf {HF'}A \ar[r]_-{p_{\mathsf F' A}} &
\mathsf {H'F'}A}
$$
for any 0-cell $A$ in $\mathbb A$.
\item The horizontally invertible modification with the components
$$
\xymatrix{
\mathsf{HF} A \ar[d]_-{\mathsf{H} y_A} \ar@{=}[r]  
\ar@{}[rd]|-{\Longdownarrow{(q^{y_A})^{-1}}} &
\mathsf{HF} A \ar[d]^-{q_{\mathsf F A}} \\
\mathsf{HF'} A  \ar[d]_-{q_{\mathsf{F'} A}} &
\mathsf{H'F} A \ar[d]^-{\mathsf H' y_A} \\
\mathsf{H'F'} A \ar@{=}[r] &
\mathsf{H'F'} A}
$$
for any 0-cell $A$ in $\mathbb A$.
\end{itemize}

The \underline{2-cell} part sends a modification in the first diagram to the modification with components in the second diagram, for any 0-cell $A$ in $\mathbb A$:
$$
\xymatrix{
\mathsf F \ar[r]^-x \ar[d]_-y \ar@{}[rd]|-{\Longdownarrow \Phi} &
\mathsf{F'} \ar[d]^-v \\
\mathsf{F^{\prime \prime}} \ar[r]_-z &
\mathsf{F^{\prime \prime \prime }}}
\qquad \qquad
\xymatrix{
\mathsf {HF} A \ar[r]^-{\mathsf H x_A} \ar[d]_-{\mathsf H y_A} 
\ar@{}[rd]|-{\Longdownarrow {\mathsf H \Phi_A}} &
\mathsf{HF'} A \ar[d]^-{\mathsf H v_A} \\
\mathsf{HF^{\prime \prime}} A \ar[r]_-{\mathsf H z_A} &
\mathsf{HF^{\prime \prime \prime}}  A.}
$$

Direct computation verifies the commutativity of 
$$
\xymatrix@C=45pt{
\ldb \mathbb A,\mathbb B \rdb 
\ar[r]^-{\ldb 1,\mathfrak r^{\mathbb D}_{\mathbbm 1,\mathbb B} \rdb}
\ar[d]_-{\mathfrak r^{\mathbb D}_{\mathbb A,\mathbb B}} &
\ldb \mathbb A, \ldb \ldb \mathbb B,\mathbb D \rdb,\mathbb D \rdb \rdb \\
\ldb \ldb \mathbb B,\mathbb D \rdb,\ldb \mathbb A,\mathbb D \rdb \rdb
\ar[r]_-{\mathfrak r^{\mathbb D}_{\ldb \mathbb B,\mathbb D \rdb,\ldb \mathbb A,\mathbb D \rdb }}  &
\ldb 
\ldb \ldb \mathbb A,\mathbb D \rdb,\mathbb D \rdb, 
\ldb \ldb \mathbb B,\mathbb D \rdb,\mathbb D \rdb 
\rdb
\ar[u]_-{\ldb  \mathfrak r^{\mathbb D}_{\mathbbm 1,\mathbb A},1\rdb}}
$$
and the equality of 
$$
\xymatrix@C=45pt{
\ldb \mathbb A,\mathbb D \rdb 
\ar[r]^-{\mathfrak r^{\mathbb D}_{\mathbbm 1,\mathbb \ldb \mathbb A,\mathbb D \rdb}} &
\ldb \ldb \ldb \mathbb A,\mathbb D \rdb,\mathbb D \rdb ,\mathbb D \rdb
\ar[r]^-{\ldb  \mathfrak r^{\mathbb D}_{\mathbbm 1,\mathbb A},1\rdb} &
\ldb \mathbb A,\mathbb D \rdb }
$$
to the identity double functor, for any double categories $\mathbb A$, $\mathbb B$ and $\mathbb D$.
It follows from these properties that the double functors
\begin{equation} \label{eq:f}
\mathfrak f^{\mathbb D}_{\mathbb A, \mathbb B}=\Bigl(
\xymatrix@C=45pt{
\ldb \mathbb A,\ldb \mathbb B,\mathbb D \rdb  \rdb 
\ar[r]^-{\mathfrak r^{\mathbb D}_{\mathbb A,\ldb \mathbb B,\mathbb D \rdb}} &
\ldb \ldb \ldb \mathbb B,\mathbb D \rdb,\mathbb D \rdb ,  
\ldb \mathbb A,\mathbb D \rdb \rdb 
\ar[r]^-{\ldb \mathfrak r^{\mathbb D}_{\mathbbm 1,\mathbb B},1\rdb } &
\ldb \mathbb B,\ldb \mathbb A,\mathbb D \rdb  \rdb }\Bigr)
\end{equation}
constitute an idempotent natural transformation
$$
\xymatrix@C=40pt{
\mathsf{DblCat}^{\mathsf{op}} \times \mathsf{DblCat}^{\mathsf{op}}
\ar[r]^-{1\times \ldb -,\mathbb D \rdb}
\ar[d]_-{\mathsf {flip}} 
\ar@{}[rrd]|-{\longdownarrow  {\mathfrak f^{\mathbb D}}} &
\mathsf{DblCat}^{\mathsf{op}} \times \mathsf{DblCat} 
\ar[r]^-{\ldb -,- \rdb} &
\mathsf{DblCat}  \\
\mathsf{DblCat}^{\mathsf{op}} \times \mathsf{DblCat}^{\mathsf{op}}
\ar[rr]_-{1\times \ldb -,\mathbb D \rdb} &&
\mathsf{DblCat}^{\mathsf{op}} \times \mathsf{DblCat} 
\ar[u]_-{\ldb -,- \rdb}}
$$
for any double category $\mathbb D$.
This $\mathfrak f$ is natural in its lower indices by the naturality of $\mathfrak r$. 
It is natural in the upper index as well (here the upper index no longer refers to extranaturality). In order to see that, both naturality and extranaturality of $\mathfrak r$ are needed.

The extranatural transformation $\mathfrak r$ in this section and $\mathfrak l$ in Section \ref{sec:l} together render commutative the following diagram, for any double categories $\mathbb A$, $\mathbb B$ and $\mathbb D$.
\begin{equation} \label{eq:l-r_pentagon}
\scalebox{1}{$
\xymatrix@C=34pt{
% 1.1
\ldb \mathbb A,\mathbb B\rdb 
\ar[rr]^-{\mathfrak l^{\ldb \mathbb D,\mathbb A\rdb}_{\mathbb A,\mathbb B}} 
\ar[d]_-{\mathfrak l^{\mathbb D}_{\mathbb A,\mathbb B}} &&
% 1.3
\ldb \ldb \mathbb D,\mathbb A\rdb,\mathbb A \rdb,\ldb \ldb \mathbb D,\mathbb A\rdb\mathbb B\rdb 
\ar[d]^-{\ldb \mathfrak r^{\mathbb A}_{\mathbbm1,\mathbb D},1\rdb}  \\
% 2.1
\ldb \ldb \mathbb D,\mathbb A\rdb,\ldb \mathbb D,\mathbb B\rdb \rdb 
\ar[r]^-{\mathfrak r^{\mathbb B}_{\ldb \mathbb D,\mathbb A\rdb,\ldb \mathbb D,\mathbb B\rdb}}
\ar@/_1.3pc/[rr]_-{\mathfrak f^{\mathbb B}_{\ldb \mathbb D,\mathbb A\rdb,\mathbb D}} &
% 2.2
\ldb\ldb \ldb \mathbb D,\mathbb B\rdb,\mathbb B \rdb,
\ldb \ldb \mathbb D,\mathbb A\rdb,\mathbb B \rdb \rdb 
\ar[r]^-{\ldb \mathfrak r^{\mathbb B}_{\mathbbm1,\mathbb D},1\rdb} &
% 2.3
\ldb \mathbb D,\ldb\ldb \mathbb D,\mathbb A\rdb,\mathbb B \rdb \rdb
}$}
\end{equation}

\subsection{Representability of the functor $\mathsf{DblCat}(\mathbb G,\ldb \mathbb G,-\rdb):\mathsf{DblCat} \to \mathsf{Set}$}
\label{sec:represent}

In this section we investigate the functor in the title, for the double category $\mathbb G$ of Section \ref{sec:DblCat}.

By the description of $\mathbb G$ in Section \ref{sec:DblCat}, for any double category $\mathbb A$ the double functors $\mathbb G \to \ldb \mathbb G,\mathbb A\rdb$ correspond bijectively to the 2-cells of $\ldb \mathbb G,\mathbb A\rdb$. 
The 0-cells at the corners of such a 2-cell are double functors denoted as $(A,-):\mathbb G  \to \mathbb A$, for all 0-cells $A\in \{X,Y,Z,V\}$ of $\mathbb G$. 
The top and bottom horizontal 1-cells are horizontal pseudotransformations labelled by the horizontal 1-cells $h\in \{t,b\}$ in $\mathbb G$. We denote their components by
$$
\xymatrix{
(A,P) \ar[r]^-{(h,P)} \ar[d]_-{(A,v)} \ar@{}[rd]|-{\Longdownarrow{(h,v)}} &
(C,P) \ar[d]^-{(C,v)} \\
(A,Q) \ar[r]_-{(h,Q)} &
(C,Q)} 
\qquad \qquad
\xymatrix{
(A,P) \ar[r]^-{(A,n)} \ar@{=}[d] \ar@{}[rrd]|-{
\Longdownarrow{(h,n)}} &
(A,R) \ar[r]^-{(h,R)} &
(C,R) \ar@{=}[d]  \\
(A,P) \ar[r]_-{(h,P)} &
(C,P) \ar[r]_-{(C,n)} &
(C,R)}
$$
for any horizontal 1-cell $n$ and vertical 1-cell $v$ in $\mathbb G$. They satisfy the naturality condition
\begin{equation} \label{eq:nat_h}
\xymatrix{
(A,X) \ar[r]^-{(A,t)} \ar[d]_-{(A,l)} \ar@{}[rd]|-{\Longdownarrow{(A,\tau)}} &
(A,Y) \ar[r]^-{(h,Y)} \ar[d]|-{(A,r)} \ar@{}[rd]|-{\Longdownarrow{(h,r)}} &
(C,Y) \ar[d]^-{(C,r)} \\
(A,Z) \ar[r]_-{(A,b)} \ar@{=}[d] &
(A,V) \ar[r]_-{(h,V)} \ar@{}[d]|-{\Longdownarrow{(h,b)}} &
(C,V) \ar@{=}[d] \\
(A,Z) \ar[r]_-{(h,Z)} &
(C,Z) \ar[r]_-{(C,b)} &
(C,V)}
\raisebox{-38pt}{$=$}
\xymatrix{
(A,X) \ar[r]^-{(A,t)} \ar@{=}[d] &
(A,Y) \ar[r]^-{(h,Y)} \ar@{}[d]|-{\Longdownarrow{(h,t)}} &
(C,Y) \ar@{=}[d] \\
(A,X) \ar[r]^-{(h,X)} \ar[d]_-{(A,l)} \ar@{}[rd]|-{\Longdownarrow{(h,l)}} &
(C,X) \ar[r]^-{(C,t)} \ar[d]|-{(C,l)} \ar@{}[rd]|-{\Longdownarrow{(C,\tau)}} &
(C,Y) \ar[d]^-{(C,r)} \\
(A,Z) \ar[r]_-{(h,Z)} &
(C,Z) \ar[r]_-{(C,b)} &
(C,V).}
\end{equation}
Symmetrically, the left and right vertical 1-cells are vertical pseudotransformations labelled by the vertical 1-cells $w \in \{ l,r \}$ of $\mathbb G$; with components denoted by
$$
\xymatrix@R=30pt{
(A,P) \ar[r]^-{(A,n)} \ar[dd]_-{(w,P)} \ar@{}[rd]|-{\Longdownarrow{(w,n)}} &
(A,R) \ar[dd]^-{(w,R)} \\
&\\
(B,P) \ar[r]_-{(B,n)} &
(B,R)}
\qquad \qquad
\xymatrix{
(A,P) \ar@{=}[r] \ar[d]_-{(w,P)} \ar@{}[rd]|-{
\Longdownarrow{(w,v)}} &
(A,P) \ar[d]^-{(A,v)} \\
(B,P) \ar[d]_-{(B,v)} &
(A,Q) \ar[d]^-{(w,Q)} \\
(B,Q) \ar@{=}[r] &
(B,Q)}
$$
for any horizontal 1-cell $n$ and vertical 1-cell $v$ in $\mathbb G$. They satisfy the naturality condition
\begin{equation} \label{eq:nat_v}
\xymatrix{
(A,X) \ar@{=}[r] \ar[d]_-{(w,X)} \ar@{}[rd]|-{\Longdownarrow{(w,l)}} &
(A,X) \ar[r]^-{(A,t)} \ar[d]|-{(A,l)} \ar@{}[rd]|-{\Longdownarrow{(A,\tau)}} &
(A,Y) \ar[d]^-{(A,r)} \\
(B,X) \ar[d]_-{(B,l)} &
(A,Z) \ar[r]^-{(A,b)} \ar[d]_-{(w,Z)} \ar@{}[rd]|-{\Longdownarrow{(w,b)}} &
(A,V) \ar[d]^-{(w,V)} \\
(B,Z) \ar@{=}[r] &
(B,Z) \ar[r]_-{(B,b)} &
(B,V)}
\raisebox{-38pt}{$=$}
\xymatrix{
(A,X) \ar[r]^-{(A,t)} \ar[d]_-{(w,X)} \ar@{}[rd]|-{\Longdownarrow{(w,t)}} &
(A,Y) \ar@{=}[r] \ar[d]|-{(w,Y)} \ar@{}[rd]|-{\Longdownarrow{(w,r)}} &
(A,Y) \ar[d]^-{(A,r)} \\
(B,X) \ar[r]^-{(B,t)} \ar[d]_-{(B,l)} \ar@{}[rd]|-{\Longdownarrow{(B,\tau)}} &
(B,Y) \ar[d]^-{(B,r)} &
(A,V) \ar[d]^-{(w,V)} \\
(B,Z) \ar[r]_-{(B,b)} &
(B,V) \ar@{=}[r] &
(B,V). }
\end{equation}
Finally, the 2-cell itself is a modification with components denoted by
$$
\xymatrix{
(X,A) \ar[r]^-{(t,A)} \ar[d]_-{(l,A)} \ar@{}[rd]|-{\Longdownarrow {(\tau,A)}} &
(Y,A) \ar[d]^-{(r,A)} \\
(Z,A) \ar[r]_-{(b,A)} &
(V,A)}
$$
for all 0-cells $A\in \{X,Y,Z,V\}$ of $\mathbb G$. They satisfy the horizontal compatibility conditions 
\begin{equation} \label{eq:mod_h}
\xymatrix{
(X,A) \ar[r]^-{(X,h)} \ar[d]_-{(l,A)} \ar@{}[rd]|-{\Longdownarrow {(l,h)}} &
(X,C) \ar[r]^-{(t,C)} \ar[d]|-{(l,C)} \ar@{}[rd]|-{\Longdownarrow {(\tau,C)}} &
(Y,C) \ar[d]^-{(r,C)} \\
(Z,A) \ar[r]_-{(Z,h)} \ar@{=}[d] \ar@{}[rrd]|-{\Longdownarrow {(b,h)}} &
(Z,C) \ar[r]_-{(b,C)} &
(V,C) \ar@{=}[d] \\
(Z,A) \ar[r]_-{(b,A)} &
(V,A) \ar[r]_-{(V,h)} &
(V,C)}
\raisebox{-38pt}{$=$}
\xymatrix{
(X,A) \ar[r]^-{(X,h)} \ar@{=}[d]  \ar@{}[rrd]|-{\Longdownarrow {(t,h)}} &
(X,C) \ar[r]^-{(t,C)} &
(Y,C) \ar@{=}[d] \\
(X,A) \ar[r]^-{(t,A)} \ar[d]_-{(l,A)} \ar@{}[rd]|-{\Longdownarrow {(\tau,A)}} &
(Y,A) \ar[r]^-{(Y,h)} \ar[d]|-{(r,A)} \ar@{}[rd]|-{\Longdownarrow {(r,h)}} &
(Y,C) \ar[d]^-{(r,C)} \\
(Z,A) \ar[r]_-{(b,A)} &
(V,A) \ar[r]_-{(V,h)} &
(V,C)}
\end{equation}
for all horizontal 1-cells $h\in \{t,b\}$ in $\mathbb G$; and the vertical compatibility conditions 
\begin{equation} \label{eq:mod_v}
\xymatrix{
(X,A) \ar[r]^-{(t,A)} \ar[d]_-{(l,A)} \ar@{}[rd]|-{\Longdownarrow {(\tau,A)}} &
(Y,A) \ar@{=}[r] \ar[d]|-{(r,A)} \ar@{}[rd]|-{\Longdownarrow {(r,w)}} &
(Y,A) \ar[d]^-{(Y,w)} \\
(Z,A) \ar[r]^-{(b,A)} \ar[d]_-{(Z,w)} \ar@{}[rd]|-{\Longdownarrow {(b,w)}} &
(V,A) \ar[d]^-{(V,w)} &
(Y,B) \ar[d]^-{(r,B)} \\
(Z,B) \ar[r]_-{(b,B)} &
(V,B) \ar@{=}[r]  &
(V,B)}
\raisebox{-38pt}{$=$}
\xymatrix{
(X,A) \ar@{=}[r] \ar[d]_-{(l,A)} \ar@{}[rd]|-{\Longdownarrow {(l,w)}} &
(X,A) \ar[r]^-{(t,A)} \ar[d]|-{(X,w)} \ar@{}[rd]|-{\Longdownarrow {(t,w)}} &
(Y,A) \ar[d]^-{(Y,w)} \\
(Z,A) \ar[d]_-{(Z,w)} &
(X,B) \ar[d]_-{(l,B)} \ar[r]^-{(t,B)} \ar@{}[rd]|-{\Longdownarrow {(\tau,B)}} &
(Y,B) \ar[d]^-{(r,B)} \\
(Z,B) \ar@{=}[r]  &
(Z,B) \ar[r]_-{(b,B)} &
(V,B) }
\end{equation}
for all vertical 1-cells $w\in \{l,r\}$ in $\mathbb G$.

From all that we can read off that the functor $\mathsf{DblCat}(\mathbb G,\ldb \mathbb G,-\rdb):\mathsf{DblCat} \to \mathsf{Set}$ is represented by the following double category. 
\begin{itemize}
\item
The \underline{0-cells} are pairs $(A,B)$ of 0-cells in $\mathbb G$.
\item
There are two kinds of non-identity  \underline{horizontal 1-cells} $(A,h)$ and $(h,A)$; both are ordered pairs of a 0-cell $A$, and a horizontal 1-cell $h\in\{t,b\}$ in $\mathbb G$.
\item
Symmetrically, there are two kinds of non-identity \underline{vertical 1-cells} $(A,v)$ and $(v,A)$; both are ordered pairs of a 0-cell $A$, and a vertical 1-cell $v\in\{l,r\}$ in $\mathbb G$.
\item
There are non-identity \underline{2-cells} of ordered pairs
\begin{itemize}
\item[{-}]
$(A,\tau)$ and $(\tau, A)$, for all 0-cells $A$ in $\mathbb G$
\item[{-}]
$(h,v)$ and $(v,h)$ for horizontal 1-cells $h\in\{t,b\}$ and vertical 1-cells $v\in\{l,r\}$ in $\mathbb G$
\item[{-}]
vertically invertible 2-cells $(h,n)$ for horizontal 1-cells $h,n\in\{t,b\}$ in $\mathbb G$ 
\item[{-}]
horizontally invertible 2-cells $(v,w)$ for vertical 1-cells $v,w\in\{l,r\}$ in $\mathbb G$.
\end{itemize}
All further cells are generated by their compositions modulo the associativity and unitality conditions, the middle four interchange law, and the identities \eqref{eq:nat_h}, \eqref{eq:nat_v}, \eqref{eq:mod_h} and \eqref{eq:mod_v}.
\end{itemize}

\subsection{Representability of the functors $\mathsf{DblCat}(\mathbb A,\ldb \mathbb B,-\rdb):\mathsf{DblCat} \to \mathsf{Set}$}
\label{sec:left_adjoint}

In this section we investigate the functors in the title, for any double categories $\mathbb A$ and $\mathbb B$.

Consider a functor $\mathsf U:\mathsf C \to \mathsf C'$ between locally presentable categories. Bourke and Gurski's \cite[Lemma 3.9]{BourkeGurski} says that the functor $\mathsf C'(X,\mathsf U(-)):\mathsf C \to \mathsf{Set}$ is representable for all objects $X$ --- that is, $\mathsf U$ possesses a left adjoint --- if and only if $\mathsf C'(G_i,\mathsf U(-)):\mathsf C \to \mathsf{Set}$ is representable for all $G_i$ in a strong generator of $\mathsf C'$.

The category $\mathsf{DblCat}$ is locally presentable by \cite[Theorem 4.1]{FiorePaoliPronk}. The double category $\mathbb G$ in Section \ref{sec:represent} is a strong generator of $\mathsf{DblCat}$; see Section \ref{sec:DblCat}. So from the representability result of Section \ref{sec:represent} we conclude by Bourke and Gurski's lemma that the functor 
$\mathsf{DblCat}(\mathbb B,\ldb \mathbb G,-\rdb):\mathsf{DblCat} \to \mathsf{Set}$ is representable for any double category $\mathbb B$.

The 0-cell part of the iso double functor \eqref{eq:f} yields a bijection
$\mathsf{DblCat}(\mathbb B,\ldb \mathbb G,\mathbb A \rdb) \cong 
\mathsf{DblCat}(\mathbb G,\ldb \mathbb B,\mathbb A \rdb)$ 
for any double categories $\mathbb A$ and $\mathbb B$; which is natural in $\mathbb A$ by the extranaturality and naturality of $\mathfrak r$.
Therefore also the functor $\mathsf{DblCat}(\mathbb G,\ldb \mathbb B,- \rdb):\mathsf{DblCat} \to \mathsf{Set}$ is representable for any double category $\mathbb B$.
Applying again Bourke and Gurski's lemma, we obtain the representability of 
$\mathsf{DblCat}(\mathbb A,\ldb \mathbb B,-\rdb)$ 
for any double categories $\mathbb A$ and $\mathbb B$.

In other words, the functor $\ldb \mathbb B,-\rdb: \mathsf{DblCat} \to \mathsf{DblCat}$ possesses a left adjoint  for any double category $\mathbb B$ which we denote by $-\otimes \mathbb B$. By the functoriality in $\mathbb  B$, it gives raise to a double functor $\otimes:\mathsf{DblCat} \times \mathsf{DblCat} \to \mathsf{DblCat}$.
For any double categories $\mathbb A$ and $\mathbb B$, an explicit description of $\mathbb A \otimes \mathbb B$ can be given analogously to Section \ref{sec:represent}.

%%%%%%%%%%%%%%%%%       SEC 2    %%%%%%%%%%%%%%%%%%%%%%%%%%%%%%

\section{Coherence}
\label{sec:coherence}

This section is devoted to the proof that the double functor $\otimes:\mathsf{DblCat} \times \mathsf{DblCat} \to \mathsf{DblCat}$ of Section \ref{sec:left_adjoint} renders $\mathsf{DblCat}$ a symmetric monoidal category (which is then closed with the internal hom functors $\ldb \mathbb B,-\rdb$, for all double categories $\mathbb B$).

In what follows, the unit of the adjunction $-\otimes \mathbb B \dashv \ldb \mathbb B,-\rdb$ will be denoted by $\eta^{\mathbb B}_{\mathbb A}: \mathbb A\to \ldb \mathbb B,\mathbb A \otimes \mathbb B\rdb$, and the counit will be denoted by $\epsilon^{\mathbb B}_{\mathbb A}: \ldb \mathbb B,\mathbb A \rdb \otimes \mathbb B \to \mathbb A$, for all double categories $\mathbb A$ and $\mathbb B$.

\subsection{The associativity natural isomorphism}
\label{sec:alpha}

For any double category $\mathbb C$, consider the natural transformation
$$
\xymatrix{
\mathsf{DblCat}^{\mathsf{op}} \times \mathsf{DblCat}^{\mathsf{op}}
\ar[r]^-\otimes
\ar[d]_-{1\times \ldb -,\mathbb C \rdb }
\ar@{}[rd]|-{\longdownarrow {\mathfrak a^{\mathbb C}}} &
\mathsf{DblCat}^{\mathsf{op}} 
\ar[d]^-{ \ldb -,\mathbb C \rdb} \\
\mathsf{DblCat}^{\mathsf{op}} \times \mathsf{DblCat}
\ar[r]_-{\ldb -,\mathbb - \rdb} &
\mathsf{DblCat}}
$$
with the components
$$
\mathfrak a^{\mathbb C}_{\mathbb A,\mathbb B}:=\Bigl(
\xymatrix{
\ldb \mathbb A\otimes \mathbb B, \mathbb C\rdb
\ar[r]^-{\mathfrak l^{\mathbb B}_{\mathbb A\otimes  \mathbb B, \mathbb C}} &
\ldb 
\ldb \mathbb B,\mathbb A\otimes \mathbb B \rdb, 
\ldb \mathbb B,\mathbb C \rdb
\rdb
\ar[r]^-{\ldb  \eta^{\mathbb B}_{\mathbb A},1 \rdb} &
\ldb  \mathbb A,\ldb \mathbb B,\mathbb C \rdb \rdb}
\Bigr) 
$$
at any double categories $\mathbb A$ and $\mathbb B$.
It is natural in $\mathbb A$ and $\mathbb C$ by the naturality of $\mathfrak l$ and $\eta$ (the upper index $\mathbb C$ of $\mathfrak a$ no longer refers to extranaturality). It is natural in $\mathbb B$ as well which follows by the extranaturality of $\mathfrak l$ and $\eta$ together with the naturality of $\mathfrak l$.

The 0-cell parts of the iso double functors in \eqref{eq:f} yield bijections in the columns of the commutative diagram
\begin{equation} \label{eq:a_inv}
\xymatrix@C=15pt{
\mathsf{DblCat}(\mathbb D,\ldb \mathbb A \otimes\mathbb B,\mathbb C \rdb)
\ar[rr]^-{\mathsf{DblCat}(\mathbb D,\mathfrak a^{\mathbb C})}
\ar[d]_-{\mathfrak f^{\mathbb C}_0} &&
\mathsf{DblCat}(\mathbb D,
\ldb \mathbb A ,\ldb \mathbb B, \mathbb C \rdb \rdb ) \\
\mathsf{DblCat}(\mathbb A \otimes \mathbb B,\ldb \mathbb D,\mathbb C \rdb )
\ar[r]_-{\raisebox{-10pt}{${}_\cong$}} &
\mathsf{DblCat}(\mathbb A,
\ldb \mathbb B,\ldb \mathbb D, \mathbb C \rdb \rdb ) 
\ar[r]_-{\raisebox{-10pt}{${}_{\mathsf{DblCat}(\mathbb A,\mathfrak f^{\mathbb C})}$}} &
\mathsf{DblCat}(\mathbb A,
\ldb \mathbb D,\ldb \mathbb B, \mathbb C \rdb \rdb ) .
\ar[u]_-{\mathfrak f^{\ldb \mathbb B, \mathbb C \rdb}_0}}
\end{equation}
Since all of the occurring maps but the top row are known to be bijections, we conclude that so is the top row. Whence by Yoneda's lemma $\mathfrak a^{\mathbb C}$ is a natural isomorphism.
Using the adjunction isomorphisms in the first and last steps, we obtain a natural isomorphism
$$
\xymatrix@R=20pt{
\mathsf{DblCat}(\mathbb A \otimes (\mathbb B \otimes \mathbb C),\mathbb D) \ar[r]^-\cong &
\mathsf{DblCat}(\mathbb A , \ldb \mathbb B \otimes \mathbb C,\mathbb D\rdb)
\ar[d]^-{\mathsf{DblCat}(\mathbb A , \mathfrak a^{\mathbb D})} \\
& \mathsf{DblCat}(\mathbb A , \ldb \mathbb B ,\ldb \mathbb C,\mathbb D\rdb \rdb) \ar[r]_-\cong &
\mathsf{DblCat}((\mathbb A \otimes \mathbb B) \otimes \mathbb C,\mathbb D).}
$$
By Yoneda's lemma again, it determines a natural isomorphism $\alpha_{\mathbb A,\mathbb B,\mathbb C}:(\mathbb A \otimes \mathbb B) \otimes  \mathbb C \to \mathbb A \otimes (\mathbb B \otimes \mathbb C)$ which is our candidate associativity natural isomorphism.

\subsection{The pentagon condition}
\label{sec:pentagon}

By Yoneda's lemma, Mac Lane's pentagon condition on the natural isomorphism $\alpha$ of Section \ref{sec:alpha} is equivalent to the commutativity of the exterior of the diagram of Figure \ref{fig:pent}; hence also to the commutativity of the diagram of Figure \ref{fig:pent_eq}.

The left column in Figure \ref{fig:pent_eq} is equal to $\ldb \alpha_{\mathbb A,\mathbb B,\mathbb C},1\rdb$.
The triangles marked by $(\ast)$ commute by the naturality of $\mathfrak l$ and a triangle condition on the adjunction $-\otimes \mathbb C \dashv \ldb \mathbb C,-\rdb$, yielding the commutative diagram
\begin{equation} \label{eq:eps_a_l}
\xymatrix{
\ldb \mathbb P,\mathbb K\rdb 
\ar[r]^-{\mathfrak l^{\mathbb C}}
\ar[d]_-{\ldb \epsilon^{\mathbb C},1\rdb} &
\ldb
       \ldb \mathbb C,\mathbb P \rdb,
       \ldb \mathbb C,\mathbb K \rdb
\rdb
\ar[d]_-{\ldb \ldb 1,\epsilon^{\mathbb C}\rdb,1\rdb}
\ar@/^1.7pc/@{=}[rd] \\
\ldb
      \ldb \mathbb C,\mathbb P \rdb \otimes \mathbb C,
      \mathbb K
\rdb
\ar[r]^-{\mathfrak l^{\mathbb C}}
\ar@/_1.3pc/[rr]_-{\mathfrak a^{\mathbb K}} &
\ldb
        \ldb
               \mathbb C,
               \ldb \mathbb C,\mathbb P \rdb \otimes \mathbb C
        \rdb,
        \ldb \mathbb C,\mathbb K \rdb
\rdb
\ar[r]^-{\ldb \eta^{\mathbb C},1\rdb} &
\ldb
       \ldb \mathbb C,\mathbb P \rdb,
       \ldb \mathbb C,\mathbb K \rdb
\rdb}
\end{equation}
for any double categories $\mathbb C$, $\mathbb P$ and $\mathbb K$.
The region marked by $(\ast\ast)$ commutes by \eqref{eq:l_comm} and extranaturality of $\mathfrak l$, yielding the commutative diagram
$$
\scalebox{.94}{$
\xymatrix@C=2pt@R=35pt{
% 1.1
\ldb \mathbb P\tp \mathbb K \rdb
\ar[rr]^-{\mathfrak l^{\mathbb B \otimes \mathbb C}}
\ar[d]_-{\mathfrak l^{\mathbb C}}
\ar@{}[rrd]|-{\eqref{eq:l_comm}} &&
% 1.3
\ldb \hspace{-1pt}
        \ldb \mathbb B \!  \otimes  \! \mathbb C\tp \mathbb P \rdb \tp
        \ldb \mathbb B \! \otimes \! \mathbb C\tp \mathbb K \rdb
\hspace{-1pt} \rdb
\ar[d]_-{\ldb 1,\mathfrak l^{\mathbb C}\rdb}
\ar@/^6.5pc/[dd]^(.24){\ldb 1,\mathfrak a^{\mathbb K}\rdb} \\
% 2.1
\ldb \hspace{-1pt}
        \ldb \mathbb C\tp \mathbb P \rdb \tp
        \ldb \mathbb C\tp \mathbb K \rdb
\hspace{-1pt} \rdb
\ar[r]^-{\raisebox{10pt}{${}_{\mathfrak l^{\ldb \mathbb C,\mathbb B  \otimes \mathbb C\rdb}}$}}
\ar[d]_-{\mathfrak l^{\mathbb B}}
&
% 2.2
 \ldb \hspace{-1pt}
    \ldb \hspace{-1pt}
       \ldb \mathbb C\tp\mathbb B \! \otimes \! \mathbb C\rdb\tp
       \ldb \mathbb C\tp \mathbb P \rdb
    \hspace{-1pt} \rdb \tp
    \ldb \hspace{-1pt}
        \ldb \mathbb C\tp\mathbb B \! \otimes \!  \mathbb C\rdb\tp
        \ldb \mathbb C\tp \mathbb K \rdb
    \hspace{-1pt} \rdb
\hspace{-1pt} \rdb
\ar[r]^-{\raisebox{10pt}{${}_{\ldb \mathfrak l^{\mathbb C},1 \rdb}$}}
\ar[d]^-{\ldb 1, \ldb \eta^{\mathbb C},1 \rdb \rdb} &
% 2.3
\ldb \hspace{-1pt}
       \ldb \mathbb B \! \otimes \! \mathbb C\tp \mathbb P \rdb \tp
       \ldb \hspace{-1pt}
        \ldb \mathbb C\tp \mathbb B \! \otimes \! \mathbb C\rdb\tp
        \ldb \mathbb C\tp \mathbb K \rdb
    \hspace{-1pt} \rdb
\hspace{-1pt} \rdb
\ar[d]_-{\ldb 1, \ldb \eta^{\mathbb C},1 \rdb \rdb}  \\
% 3.1
\ldb \hspace{-1pt}
      \ldb \mathbb B\tp \ldb \mathbb C\tp\mathbb P \rdb \hspace{-1pt} \rdb\tp
      \ldb \mathbb B\tp \ldb \mathbb C\tp \mathbb K \rdb \hspace{-1pt} \rdb
\hspace{-1pt} \rdb
\ar[r]^-{\raisebox{10pt}{${}_{\ldb \ldb \eta^{\mathbb C},1 \rdb , 1\rdb}$}}
\ar@/_1.3pc/[rr]_-{\ldb \mathfrak a^{\mathbb P}, 1\rdb} &
% 3.2
\ldb \hspace{-1pt}
      \ldb \hspace{-1pt}
            \ldb \mathbb C\tp \mathbb B \! \otimes \! \mathbb C\rdb\tp 
            \ldb \mathbb C\tp \mathbb P \rdb 
      \hspace{-1pt} \rdb\tp
      \ldb 
            \mathbb B\tp
            \ldb \mathbb C\tp\mathbb K \rdb 
       \hspace{-1pt} \rdb
\hspace{-1pt} \rdb
\ar[r]^-{\raisebox{10pt}{${}_{\ldb \mathfrak l^{\mathbb C},1 \rdb}$}} &
% 3.3
\ldb \hspace{-1pt}
      \ldb 
            \mathbb B  \! \otimes \! \mathbb C\tp
            \mathbb P 
      \rdb\tp
      \ldb 
            \mathbb B\tp
            \ldb \mathbb C\tp\mathbb K \rdb 
       \hspace{-1pt} \rdb
\hspace{-1pt} \rdb}$}
$$
for any double categories $\mathbb B$, $\mathbb C$, $\mathbb P$ and $\mathbb K$.

\subsection{The unitality natural isomorphisms}
\label{sec:unitors}

For any double categories $\mathbb A$ and $\mathbb K$ there are natural isomorphisms
\begin{eqnarray*}
\mathsf{DblCat}(\mathbb A,\mathbb K)&\cong&
\mathsf{DblCat}(\mathbb A,\ldb \mathbbm 1,\mathbb K\rdb)\cong
\mathsf{DblCat}(\mathbb A \otimes  \mathbbm 1 ,\mathbb K) \\
\mathsf{DblCat}(\mathbb A,\mathbb K)&\cong&
\mathsf{DblCat}(\mathbb A,\ldb \mathbbm 1,\mathbb K\rdb) 
\stackrel{\mathfrak f^{\mathbb K}_0} \longrightarrow
\mathsf{DblCat}(\mathbbm 1,\ldb \mathbb A,\mathbb K\rdb) \cong
\mathsf{DblCat}(\mathbbm 1 \otimes  \mathbb A ,\mathbb K)
\end{eqnarray*}
where $\mathfrak f^{\mathbb K}_0$ denotes the 0-cell part of the iso double functor \eqref{eq:f}.
By Yoneda's lemma, they induce respective natural isomorphisms $\varrho$ and $\lambda$ with the components 
$$
\varrho_{\mathbb A}=\Bigl(
\xymatrix@C=12pt{
\mathbb A  \otimes  \mathbbm 1 \ar[r]^-\cong &
\ldb \mathbbm 1, \mathbb A \rdb  \otimes  \mathbbm 1 \ar[r]^-{\epsilon^{\mathbbm 1}} &
\mathbb A}\Bigr)
\quad \textrm{and}\quad
\lambda_{\mathbb A}=\Bigl(
\xymatrix@C=12pt{
\mathbbm 1 \otimes  \mathbb A \ar[rr]^-{1_A \otimes 1} &&
\ldb \mathbb A,\mathbb A\rdb \otimes  \mathbb A \ar[r]^-{\epsilon^{\mathbb A}} &
\mathbb A}\Bigr)
$$
at any double category $\mathbb A$ (where $1_A:\mathbbm 1 \to \ldb \mathbb A,\mathbb A\rdb$ is the double functor which sends the single object of $\mathbbm 1$ to the identity double functor $1_A:\mathbb A \to \mathbb A$). They are our candidate unitality natural isomorphisms.

\subsection{The triangle conditions}
\label{sec:triangle}

By Yoneda's lemma, Mac Lane's triangle condition on the natural isomorphisms $\alpha$ of Section \ref{sec:alpha}  and $\lambda,\varrho$ of Section \ref{sec:unitors} is equivalent to the commutativity of the exterior of Figure \ref{fig:triang}
for any double categories $\mathbb A$, $\mathbb B$ and $\mathbb K$;
hence also to the commutativity of 
\begin{equation} \label{eq:triang}
\xymatrix@C=55pt{
\ldb \mathbb B,\mathbb K\rdb
\ar@/^1.5pc/[rr]^-{\ldb \lambda,1\rdb}
\ar[r]_-{\ldb \epsilon^{\mathbb B},1\rdb}
\ar@/_1.5pc/[rd]_-{\mathfrak l^{\mathbb B}}
\ar@{}[rd]|-{\qquad \eqref{eq:eps_a_l}} &
\ldb \ldb \mathbb B,\mathbb B\rdb \otimes \mathbb B,\mathbb K\rdb
\ar[d]^-{\mathfrak a^{\mathbb K}}
\ar[r]_-{\ldb 1_{\mathbb B} \otimes 1,1\rdb} &
\ldb \mathbbm 1 \otimes \mathbb B,\mathbb K\rdb
\ar[d]^-{\mathfrak a^{\mathbb K}} \\
& 
\ldb \ldb \mathbb B,\mathbb B\rdb ,\ldb \mathbb B,\mathbb K\rdb \rdb
\ar[r]_-{\ldb 1_{\mathbb B} ,1\rdb} &
\ldb \mathbbm 1 ,\ldb \mathbb B,\mathbb K\rdb \rdb}
\end{equation}
whose left-bottom path is \eqref{eq:l_id}; that is, the canonical (usually omitted) isomorphism.

\subsection{The symmetry}
\label{sec:symm}

The natural isomorphism
$$
\xymatrix@C=15pt{
\mathsf{DblCat}( \mathbb B  \otimes \mathbb A,\mathbb K )
\ar[r]^-\cong &
\mathsf{DblCat}(  \mathbb B ,\ldb \mathbb A,\mathbb K \rdb )
\ar[r]^-{\mathfrak f^{\mathbb K}_0} &
\mathsf{DblCat}(  \mathbb A ,\ldb \mathbb B,\mathbb K \rdb)
\ar[r]^-\cong &
\mathsf{DblCat}(  \mathbb A  \otimes  \mathbb B,\mathbb K)}
$$
constructed from the 0-cell part of $\mathfrak f$ in \eqref{eq:f} induces a natural isomorphism
$\varphi:\otimes  \to \otimes .\mathsf{flip}$ with the components
$$
\varphi_{\mathbb A ,\mathbb B}=\Bigl(\xymatrix@C=12pt{
\mathbb A \otimes  \mathbb B
\ar[rr]^-{\mathfrak r^{\mathbb B \otimes  \mathbb A}_{\mathbbm 1,\mathbb A}
\otimes  1} &&
\ldb \ldb \mathbb A,\mathbb B \otimes \mathbb A \rdb,\mathbb B \otimes \mathbb A \rdb
\otimes \mathbb B
\ar[rr]^-{\ldb \eta^{\mathbb A},1\rdb \otimes 1} &&
\ldb \mathbb B,\mathbb B \otimes \mathbb A \rdb \otimes \mathbb B
\ar[r]^-{\epsilon^{\mathbb B}} &
\mathbb B \otimes  \mathbb A}\Bigr)
$$
at any double categories $\mathbb A$ and $\mathbb B$. It is our candidate symmetry. 

\subsection{The hexagon condition}
\label{sec:hexagon}
By Yoneda's lemma, the hexagon condition on the natural isomorphisms $\alpha$ of Section \ref{sec:alpha}  and $\varphi$ of Section \ref{sec:symm} is equivalent to the commutativity of the exterior of Figure \ref{fig:hex}
for any double categories $\mathbb A$, $\mathbb B$, $\mathbb C$ and $\mathbb K$. Hence it follows by the commutativity of the diagram of Figure \ref{fig:hex_eq}
whose left-bottom path is equal to $\mathfrak a^{\mathbb K}$. In order to see that the left column of Figure \ref{fig:hex_eq} is equal to $\mathfrak l^{\mathbb A}$, apply twice \eqref{eq:l-r_pentagon} to obtain the commutative diagram
\begin{equation} \label{eq:l-r_square}
\xymatrix@C=27pt{
% 1.1
\ldb \mathbb A\tp \mathbb B \rdb 
\ar[rr]^-{\mathfrak l
^{\ldb \hspace{-1pt} 
           \ldb \mathbb C, \mathbb A \rdb ,
           \mathbb A  
       \rdb }
_{\mathbb A, \mathbb B }}
\ar[ddd]_-{\mathfrak l^{\mathbb C}_{\mathbb A, \mathbb B }}
\ar[rd]^-{\mathfrak l^{\ldb \mathbb C, \mathbb A \rdb}_{\mathbb A, \mathbb B }} &
\ar@{}[rd]|-{\eqref{eq:l-r_pentagon}}  &
% 1.3
\ldb \hspace{-1pt} 
   \ldb \hspace{-1pt} 
      \ldb \hspace{-1pt} 
           \ldb \mathbb C\tp \mathbb A \rdb \tp 
           \mathbb A  
       \rdb \tp
       \mathbb A 
   \rdb \tp 
   \ldb \hspace{-1pt} 
     \ldb \hspace{-1pt} 
       \ldb \mathbb C\tp\mathbb A \rdb \tp
       \mathbb A
     \rdb\tp 
     \mathbb B 
   \rdb \hspace{-1pt} 
\rdb 
\ar[d]^-{\ldb \mathfrak r^{\mathbb A}_{\mathbbm 1, \ldb \mathbb C, \mathbb A\rdb},1\rdb}   \\
\ar@{}[rd]|-{\eqref{eq:l-r_pentagon}}
% 2.2
& \ldb \hspace{-1pt} \ldb \hspace{-1pt} \ldb \mathbb C \tp \mathbb A \rdb \tp   \mathbb A \rdb \tp \ldb \hspace{-1pt}  \ldb \mathbb C \tp \mathbb A\rdb \tp \mathbb B \rdb \hspace{-1pt}  \rdb 
\ar[r]^-{\mathfrak f^{\mathbb B}
_{\ldb \hspace{-1pt} \ldb \mathbb C, \mathbb A \rdb,   \mathbb A \rdb, 
\ldb  \mathbb C, \mathbb A\rdb}} 
\ar[d]^-{\ldb \mathfrak r^{\mathbb A}_{\mathbbm 1, \mathbb C},1\rdb} &
% 2.3
\ldb \hspace{-1pt}  \ldb \mathbb C \tp \mathbb A \rdb \tp 
\ldb \hspace{-1pt} 
     \ldb \hspace{-1pt} 
       \ldb \mathbb C\tp\mathbb A \rdb \tp
       \mathbb A
     \rdb\tp 
     \mathbb B 
   \rdb \hspace{-1pt} 
\rdb 
\ar[dd]^-{\ldb 1,\ldb \mathfrak r^{\mathbb A}_{\mathbbm 1, \mathbb C},1\rdb\rdb} \\
% 3.2
& \ldb \mathbb C \tp \ldb \hspace{-1pt}  \ldb \mathbb C \tp \mathbb A\rdb \tp \mathbb B \rdb \rdb 
\ar[rd]^-{\mathfrak f^{\mathbb B}_{\mathbb C,\ldb \mathbb C,\mathbb A \rdb}} \\
% 4.1
\ldb \hspace{-1pt} \ldb \mathbb C \tp \mathbb A \rdb \tp \ldb \mathbb C \tp \mathbb B \rdb \hspace{-1pt} \rdb 
\ar[ru]^-{\mathfrak f^{\mathbb B}_{\ldb \mathbb C, \mathbb A \rdb,\mathbb C}}
\ar@{=}[rr] &&
% 4.3
\ldb \hspace{-1pt} \ldb \mathbb C \tp \mathbb A \rdb \tp \ldb \mathbb C \tp \mathbb B \rdb \hspace{-1pt} \rdb }
\end{equation}
for any double categories $\mathbb A$, $\mathbb B$ and $\mathbb C$.
The region of Figure \ref{fig:hex_eq} marked by $(\ast)$ commutes by the extranaturality of $\mathfrak l$.

%%%%%%%%%%%%%%%%%%%%%%%%  SEC 3 %%%%%%%%%%%%%%%%%%%%%%%%%%

\section{Examples}
\label{sec:examples}

Although our notions of (horizontal and vertical) pseudotransformations and of corresponding modification in Section \ref{sec:[A,B]} may look quite natural, admittedly no higher principle fixes their choice. Therefore there is no {\em a'priori good} resulting Gray monoidal product of double categories. In this final section we support our construction by relating it to existing structures. Namely, we verify the monoidality of some well-known functors between  our monoidal category $(\mathsf{DblCat},\otimes)$ and some other monoidal categories which occur in the literature quite frequently.

\subsection{Monoidal functors between closed monoidal categories}
\label{sec:functors}

In any closed monoidal category we may take the mate 
$$
\mathfrak a^C_{A,B}:=
\xymatrix@C=1pt@R=1pt{
\Bigl(
[A\! \otimes \!B\tp C] \ar[r]^-{\eta^A} &
[A\tp [A\!\otimes \! B\tp C] \! \otimes \! A] \ar[rr]^-{[1,\eta^B]} &&
[A\tp [B\tp ([A\!\otimes \!B\tp C] \!\otimes \!A) \!\otimes \! B]\hspace{-1pt}]
\ar[rrrr]^-{[1,[1,\alpha]]} &&&& \\
& [A\tp [B\tp [A\!\otimes\! B\tp C] \!\otimes \! (A \!\otimes \!B)\hspace{-1pt}\hspace{-1pt}]\hspace{-1pt}]
\ar[rr]^-{[1,[1,\epsilon^{A\otimes B}]]} &&
[A\tp [B\tp C]\hspace{-1pt}] \Bigr)}
$$
of the associativity isomorphism $\alpha$ under the adjunctions $-\otimes X \dashv [X,-]$ for $X$ being the objects $A$, $B$ and $A\otimes B$; and its mate
\begin{equation} \label{eq:gen_l}
\mathfrak l^C_{A,B} :=\Bigl((\xymatrix{
[A,B] \ar[r]^-{[\epsilon^C,1]} &
[[C,A]\otimes C,B] \ar[r]^-{\mathfrak a^B} &
[[C,A],[C,B]]
}\Bigr).
\end{equation}

Consider now a functor $\mathsf H$ between closed monoidal categories. Some natural transformation $\mathsf H_2:\mathsf H(-) \otimes \mathsf H(-) \to \mathsf H(-\otimes -)$ (for both monoidal products denoted by $\otimes$) and a morphism $\mathsf H_0:I \to \mathsf H I$ (for both monoidal units denoted by $I$) render $\mathsf H$ monoidal if and only if the mate
$$
\chi_{A,B}:=\Bigl(\!  \xymatrix@C=10pt{
\mathsf H [A,B] \ar[r]^-{\eta^{\mathsf H A}} &
[\mathsf H A,\mathsf H [A,B] \otimes \mathsf H A] \ar[rr]^-{[1,\mathsf H_2]} &&
[\mathsf H A,\mathsf H ([A,B] \otimes A)]  \ar[rr]^-{[1,\mathsf H \epsilon^A]} &&
[\mathsf H A,\mathsf H B]
} \!  \Bigr)
$$
of $\mathsf H_2$ under the adjunctions $-\otimes A \dashv [A,-]$ and  $-\otimes \mathsf H A \dashv [\mathsf H A,-]$ makes the following diagrams commute.

The left and right unitality conditions translate to the commutativity of the respective diagrams in 
\begin{equation} \label{eq:unit}
\xymatrix@C=12pt@R=49pt{
I \ar[r]^-{\eta^{\mathsf H A}} \ar[d]_-{\mathsf H_0} &
[\mathsf H A,I \otimes \mathsf H A] \ar[r]^-{[1,\lambda]} &
[\mathsf H A,\mathsf H A] \\
\mathsf H I \ar[r]_-{\mathsf H \eta^A} &
\mathsf H [A,I \otimes A] \ar[r]_-{\mathsf H [1,\lambda]} &
\mathsf H [A,A]\ar[u]_-\chi}
\quad
\xymatrix@C=3pt@R=15pt{
\mathsf H A \ar[r]^-{\eta^I} \ar[d]_-{\mathsf H \eta^I} &
[I,\mathsf H A \otimes I] \ar[rr]^-{[1,\varrho]} &&
[I,\mathsf H A] \\
\mathsf H [I,A\otimes I] \ar[d]_-{\mathsf H [1,\varrho]} \\
\mathsf H [I, A] \ar[rrr]_-\chi &&&
[\mathsf H I,\mathsf H A]\ar[uu]_-{[\mathsf H_0,1]}}
\end{equation}
for all objects $A$ of the domain category.
The associativity condition translates to the commutativity of 
\begin{equation} \label{eq:assoc}
\xymatrix{
\mathsf H[A,B] \ar[r]^-\chi \ar[d]_-{\mathsf H\mathfrak l^C} &
[\mathsf H A,\mathsf H B] \ar[r]^-{\mathfrak l^{\mathsf H C}} &
[[\mathsf H C,\mathsf H A] ,[\mathsf H C,\mathsf H B] ]
\ar[d]^-{[\chi,1]} \\
\mathsf H [[ C,A] ,[C,B]] \ar[r]_-\chi &
[\mathsf H [ C,A] ,\mathsf H [C,B]] \ar[r]_-{[1,\chi]} &
[\mathsf H [ C,A] ,[\mathsf H C,\mathsf H B]] }
\end{equation}
for all objects $A$, $B$ and $C$ of the domain category. 

\subsection{The closed monoidal category $(\mathsf {DblCat},\otimes)$}
\label{sec:ox}

For the closed monoidal category $(\mathsf {DblCat},\otimes)$ of Section \ref{sec:existence} and Section \ref{sec:coherence}, 
$$
\xymatrix{
\mathbbm 1 \ar[r]^-{\eta^{\mathbb A}} &
\ldb \mathbb A,\mathbbm 1 \otimes \mathbb A\rdb \ar[r]^-{\ldb 1,\lambda\rdb} &
\ldb \mathbb A,\mathbb A\rdb}
$$
is equal to $1_{\mathbb A}$; that is, the double functor sending the single object of $\mathbbm 1$ to the identity double functor $1_{\mathbb A}:\mathbb A \to \mathbb A$ for any double category $\mathbb A$. 
The double functor
$$
\xymatrix{
\mathbb A \ar[r]^-{\eta^{\mathbbm 1}} &
\ldb \mathbbm 1,\mathbb  A \otimes \mathbbm 1\rdb \ar[r]^-{\ldb 1,\varrho \rdb} &
\ldb \mathbbm 1,\mathbb A\rdb}
$$
is the canonical isomorphism for any double category $\mathbb A$. 
The double functor of \eqref{eq:gen_l} is equal to that in Section \ref{sec:l}.

\subsection{Monoids in $(\mathsf {DblCat},\otimes)$}
\label{sec:monoid}

Monoidal 2-categories can be defined at different levels of generality. The most restrictive one in the literature is a monoid in the category of 2-categories and 2-functors with respect to the Cartesian monoidal structure. This is known as a strict monoidal 2-category.
The most general one is a single object tricategory \cite{GordonPowerStreet}; known as a monoidal bicategory. In between them are the so-called Gray monoids; these are again monoids in the category of 2-categories and 2-functors, but in this case with respect to the Gray monoidal structure \cite{Gray}.  Their importance stems from the coherence result of \cite{GordonPowerStreet}, proving that any monoidal bicategory is equivalent to a Gray monoid (as a tricategory).

Analogously, a strict monoidal double category \cite{BruniMeseguerMontanari} is a monoid in the category of double categories and double functors with respect to the Cartesian monoidal structure. In \cite{Shulman,GrandisPare:Intercategories_fw} it was generalized to a pseudomonoid in the 2-category of (pseudo) double categories and {\em pseudo} double functors and, say, vertical transformations. However, no double category analogues of Gray monoids and of monoidal bicategories seem to be available in the literature.
While the considerations in this paper do not promise any insight how to define most general monoidal (pseudo) double categories, monoids in $(\mathsf {DblCat},\otimes)$ are natural candidates for the double category analogue of Gray monoid. In this section we give their explicit characterization, similar to the characterization of Gray monoids in \cite[Lemma 4]{BaezNeuchl}.

A monoid in $(\mathsf {DblCat},\otimes)$ is equivalently a monoidal functor from the terminal double category $\mathbbm 1$ (with the trivial monoidal structure) to $(\mathsf {DblCat},\otimes)$. It can be described in terms of the data in Section \ref{sec:functors}. Namely, a  monoid structure on a double category $\mathbb A$ translates to double functors $I:\mathbbm 1 \to \mathbb A$ and $M:\mathbb A \to \ldb \mathbb A, \mathbb A\rdb $ which render commutative the diagrams of \eqref{eq:unit} and \eqref{eq:assoc}. Spelling out the details, this amounts to the following data.
\begin{itemize}
\item A distinguished 0-cell $I$.
\item For any 0-cells $X$ and $Y$, a 0-cell $X \ast Y$.
\item For any 0-cell $Y$ and any horizontal 1-cell on the left below, horizontal 1-cells on the right:
$$
\xymatrix@C=15pt{X \ar[r]^-h & X'}
\qquad 
\qquad
\xymatrix@C=20pt{X \ast Y \ar[r]^-{h  \ast Y} & X' \ast Y}
\ \textrm{and} \ 
\xymatrix@C=20pt{Y \ast X \ar[r]^-{Y \ast h} & Y \ast X'.}
$$
\item For any 0-cell $Y$ and any vertical 1-cell on the left below, vertical 1-cells on the right:
$$
\xymatrix@R=18pt{X \ar[d]^-v \\ X'}
\qquad 
\qquad
\qquad 
\qquad
\xymatrix@R=18pt{X \ast Y \ar[d]^-{v  \ast Y} \\ X' \ast Y}
\quad \raisebox{-15pt}{$\textrm{and}$} \quad
\xymatrix@R=18pt{Y \ast X \ar[d]^-{Y \ast v} \\ \ Y \ast X'.}
$$
\item For any 0-cell $Y$ and any 2-cell on the left below, 2-cells on the right:
$$
\xymatrix@C=15pt{
X \ar[r]^-h \ar[d]_-v \ar@{}[rd]|-{\Longdownarrow \omega} &
X' \ar[d]^-w \\ 
X^{\prime \prime} \ar[r]_-{k} &
X^{\prime \prime \prime}}
\qquad 
\qquad
\xymatrix@C=15pt{
X \ast Y \ar[r]^-{h \ast Y} \ar[d]_-{v \ast Y} 
\ar@{}[rd]|-{\Longdownarrow {\omega\ast Y }} &
X' \ast Y \ar[d]^-{w \ast Y} \\ 
X^{\prime \prime} \ast Y \ar[r]_-{k\ast Y } &
X^{\prime \prime \prime}\ast Y }
\quad \raisebox{-19pt}{$\textrm{and}$}\quad
\xymatrix@C=15pt{
Y \ast X \ar[r]^-{Y \ast h} \ar[d]_-{Y \ast v} 
\ar@{}[rd]|-{\Longdownarrow {Y \ast \omega}} &
Y \ast X' \ar[d]^-{Y \ast w} \\ 
Y \ast X^{\prime \prime} \ar[r]_-{Y \ast k} &
Y \ast X^{\prime \prime \prime}.}
$$
\item For any horizontal 1-cell $h$ and any vertical 1-cell $q$, 2-cells
$$
\xymatrix@C=15pt{
X \ast Y \ar[r]^-{h \ast Y} \ar[d]_-{X \ast q} 
\ar@{}[rd]|-{\Longdownarrow {h\ast q}} &
X' \ast Y \ar[d]^-{X' \ast q} \\ 
X\ast Y' \ar[r]_-{h\ast Y' } &
X'\ast Y' }
\quad \textrm{and} \quad
\xymatrix@C=15pt{
Y \ast X \ar[r]^-{Y \ast h} \ar[d]_-{q \ast X} 
\ar@{}[rd]|-{\Longdownarrow {q\ast h}} &
Y \ast X' \ar[d]^-{q \ast X'} \\ 
Y' \ast X \ar[r]_-{Y' \ast h} &
Y' \ast X'.}
$$
\item For any horizontal 1-cells $h$ and $p$, a vertically invertible 2-cell
$$
\xymatrix{
X \ast Y \ar[r]^-{X \ast p} \ar@{=}[d] \ar@{}[rrd]|-{\Longdownarrow {h\ast p}} &
X \ast Y' \ar[r]^-{h \ast Y'} &
X' \ast Y' \ar@{=}[d] \\
X \ast Y \ar[r]_-{h\ast Y} & 
X' \ast Y \ar[r]_-{X' \ast p} &
X' \ast Y'.}
$$
\item For any vertical 1-cells $v$ and $q$, a horizontally invertible 2-cell
$$
\xymatrix{
X \ast Y \ar[d]_-{v \ast Y} \ar@{=}[r] \ar@{}[rd]|-{\Longdownarrow {v\ast q}} &
X \ast Y \ar[d]^-{X \ast q} \\
X' \ast Y \ar[d]_-{X'\ast q} &
X\ast Y' \ar[d]^-{v\ast Y'} \\
X' \ast Y' \ar@{=}[r] &
X' \ast Y'.}
$$
\end{itemize}
One can memorize this as the  rule  that a pair of an $n$ dimensional and an $m$ dimensional cell is sent by the operation $\ast$ to an $n+m\leq 2$ dimensional cell.
These data are subject to the following conditions.
\begin{enumerate}[(i)]
\item For any 0-cell $X$, $X\ast -$ and $-\ast X$ are double functors $\mathbb A \to \mathbb A$.
\item $I\ast - =1_{\mathbb A}=-\ast I$.
\item For any 0-cells $X$ and $Y$, the following equalities of double functors hold.
$$
X\ast (Y\ast -) = (X \ast Y)\ast -\quad \quad
X\ast (-\ast Y) = (X \ast -)\ast Y\quad \quad
-\ast (X\ast Y) = (- \ast X)\ast Y
$$ 
\item For any 0-cell $X$, any horizontal 1-cells $h,p$ and any vertical 1-cells $v,q$, the following equalities of 2-cells hold. \\
$
\begin{array}{lll}
h\ast (X \ast q) =(h\ast X) \ast q \ \, &
(X\ast h)\ast q = X \ast (h \ast q) \ \, &
%\textrm{and} & 
h\ast (q\ast X) =(h\ast q) \ast X \\
v\ast (X\ast p)=(v\ast X) \ast p &
(X \ast v) \ast p = X \ast (v \ast p) &
%\textrm{and} &
v\ast (p\ast X)=(v\ast p) \ast X \\
h \ast (X\ast p) = (h \ast X) \ast p &
(X \ast h) \ast p = X \ast (h \ast p) &
%\textrm{and} &  
h \ast (p\ast X) = (h \ast p) \ast X \\
v \ast (X\ast q) = (v \ast X) \ast q &
(X \ast v) \ast q = X \ast (v \ast q) &
%\textrm{and} & 
v \ast (q\ast X) = (v \ast q) \ast X 
\end{array}
$
\item For any 0-cell $X$, we denote by $1_X$ the horizontal identity 1-cell; and by $1^X$ the vertical identity 1-cell on $X$.
For any horizontal 1-cell $h$ and  any vertical 1-cell $v$, the following equalities of 2-cells hold.
\\
$h\ast 1^X =\raisebox{17pt}{$\xymatrix{
\ar[r]^-{h\ast X} \ar@{}[rd]|-{\Longdownarrow 1} \ar@{=}[d] &
\ar@{=}[d] \\
\ar[r]_-{h\ast X} &}$}=h\ast 1_X$
and
$v\ast 1_X =\raisebox{17pt}{$\xymatrix{
\ar[d]_-{v\ast X} \ar@{}[rd]|-{\Longdownarrow 1} \ar@{=}[r] &
\ar[d]^-{v\ast X} \\
 \ar@{=}[r] &}$}=v\ast 1^X$
\\
$1^X\ast h =\raisebox{17pt}{$\xymatrix{
\ar[r]^-{X \ast h} \ar@{}[rd]|-{\Longdownarrow 1} \ar@{=}[d] &
\ar@{=}[d] \\
\ar[r]_-{X\ast h} &}$}=1_X\ast h$
and
$1_X\ast v =\raisebox{17pt}{$\xymatrix{
\ar[d]_-{X \ast v} \ar@{}[rd]|-{\Longdownarrow 1} \ar@{=}[r] &
\ar[d]^-{X \ast v} \\
\ar@{=}[r]  &}$}=1^X\ast v$
\item We denote the vertical composition (of vertical 1-cells and of 2-cells) by an upper dot, and the horizontal composition (of horizontal 1-cells and of 2-cells) by a lower dot.
For any composable pairs of horizontal 1-cells $h,h'$ and $p,p'$ and 
for any composable pairs of vertical 1-cells $v,v'$ and $q,q'$, the following hold. \\
$h\ast(q'\udot q)=(h\ast q')\udot (h\ast q)$ and 
$(h'.h)\ast q=(h'\ast q).(h\ast q)$ \\ 
$v\ast (p'.p)=(v\ast p').(v\ast p)$ and 
$(v'\udot v)\ast p=(v'\ast p) \udot (v\ast p)$ \\
$h\ast (p'.p)$ and $(h'.h)\ast p$
are equal to the respective 2-cells
$$
\xymatrix@C=18pt{
X\!\ast \! Y \ar[r]^-{X\ast p} \ar@{=}[d] \ar@{}[rd]|-{\Longdownarrow 1} &
X\!\ast \! Y' \ar[r]^-{X\ast  p'} \ar@{=}[d] \ar@{}[rrd]|-{\Longdownarrow {h\ast p'}} &
X \!\ast \! Y^{\prime\prime} \ar[r]^-{h\ast  Y^{\prime\prime} } &
X' \!\ast \! Y^{\prime\prime} \ar@{=}[d] \\
X\!\ast \! Y \ar[r]^-{X\ast  p} \ar@{=}[d] \ar@{}[rrd]|-{\Longdownarrow {h \ast p}} &
X \!\ast \! Y' \ar[r]^-{h\ast  Y'} &
X'\!\ast \! Y' \ar[r]^-{X' \ast  p'} \ar@{=}[d] \ar@{}[rd]|-{\Longdownarrow 1} &
X' \!\ast \! Y^{\prime\prime} \ar@{=}[d] \\
X \!\ast \! Y \ar[r]_-{h\ast  Y} &
X' \!\ast \! Y \ar[r]_-{X'\ast  p} &
X'\!\ast \! Y' \ar[r]_-{X' \ast  p'} &
X' \!\ast \! Y^{\prime\prime} } \quad
\xymatrix@C=18pt{
X\!\ast \! Y \ar[r]^-{X\ast  p} \ar@{=}[d] \ar@{}[rrd]|-{\Longdownarrow {h\ast p}}  &
X\!\ast \! Y' \ar[r]^-{h\ast  Y'} &
X' \!\ast \! Y' \ar[r]^-{h'\ast  Y'} \ar@{=}[d] \ar@{}[rd]|-{\Longdownarrow 1} &
X^{\prime\prime} \!\ast \! Y' \ar@{=}[d] \\
X\!\ast \! Y \ar[r]^-{h\ast  Y} \ar@{=}[d] \ar@{}[rd]|-{\Longdownarrow 1} &
X' \!\ast \! Y \ar[r]^-{X'\ast  p} \ar@{=}[d] \ar@{}[rrd]|-{\Longdownarrow {h' \ast p}}  &
X'\!\ast \! Y' \ar[r]^-{h' \ast  Y'} &
X^{\prime\prime}  \!\ast \! Y'\ar@{=}[d] \\
X \!\ast \! Y \ar[r]_-{h\ast  Y} &
X' \!\ast \! Y \ar[r]_-{h'\ast  Y} &
X^{\prime\prime} \!\ast \! Y \ar[r]_-{X^{\prime\prime} \ast  p} &
X^{\prime\prime}  \!\ast \!    Y' } 
$$
$v\ast (q'\udot q)$ and $(v'\udot v)\ast q$ are equal to the respective 2-cells
$$
\xymatrix@C=15pt{
X \ast Y \ar@{=}[r] \ar[d]_-{v\ast Y} \ar@{}[rd]|-{\Longdownarrow {v\ast q}} &
X \ast Y \ar@{=}[r] \ar[d]^-{X\ast q}  \ar@{}[rd]|-{\Longdownarrow 1} &
X \ast Y \ar[d]^-{X\ast q} \\
X' \ast Y \ar[d]_-{X'\ast q} &
X \ast Y' \ar@{=}[r] \ar[d]_-{v\ast Y'} \ar@{}[rd]|-{\Longdownarrow {v\ast q'}} &
X \ast Y' \ar[d]^-{X\ast q'} \\
X' \ast Y' \ar@{=}[r]  \ar[d]_-{X'\ast q'} \ar@{}[rd]|-{\Longdownarrow 1} &
X' \ast Y'  \ar[d]^-{X'\ast q'} &
X \ast Y^{\prime\prime} \ar[d]^-{v\ast Y^{\prime\prime}} \\
X' \ast Y^{\prime\prime} \ar@{=}[r] &
X' \ast Y^{\prime\prime} \ar@{=}[r] &
X' \ast Y^{\prime\prime}.}
\qquad
\xymatrix@C=15pt{
X \ast Y \ar@{=}[r] \ar[d]_-{v\ast Y}  \ar@{}[rd]|-{\Longdownarrow 1} &
X \ast Y \ar@{=}[r] \ar[d]_-{v\ast Y} \ar@{}[rd]|-{\Longdownarrow {v\ast q}} &
X \ast Y \ar[d]^-{X\ast q} \\
X' \ast Y \ar@{=}[r] \ar[d]_-{v'\ast Y} \ar@{}[rd]|-{\Longdownarrow {v'\ast q}} &
X' \ast Y \ar[d]^-{X'\ast q} &
X \ast Y' \ar[d]^-{v\ast Y'} \\
X^{\prime\prime} \ast Y \ar[d]_-{X^{\prime\prime}\ast q}  &
X' \ast Y'  \ar[d]_-{v'\ast Y'} \ar@{=}[r]  \ar@{}[rd]|-{\Longdownarrow 1} &
X' \ast Y' \ar[d]^-{v'\ast Y'} \\
X^{\prime\prime} \ast Y' \ar@{=}[r] &
X^{\prime\prime} \ast Y' \ar@{=}[r] &
X^{\prime\prime} \ast Y'.}
$$
\item For every 2-cell $\omega$, horizontal 1-cell $h$ and vertical 1-cell $v$, the following naturality conditions hold.
$$
\xymatrix@C=20pt@R=20pt{
X \ast Y \ar[r]^-{X \ast p} \ar[d]_-{X\ast q} 
\ar@{}[rd]|-{\Longdownarrow{X\ast \omega}} &
X \ast Y' \ar[r]^-{h \ast Y'} \ar[d]|-{X \ast r}
\ar@{}[rd]|-{\Longdownarrow{h\ast r}} &
X'\ast Y' \ar[d]^-{X'\ast r}  \\
X\ast Y^{\prime\prime} \ar[r]_-{X \ast s} \ar@{=}[d] 
\ar@{}[rrd]|-{\Longdownarrow{h\ast s}} &
X \ast Y^{\prime\prime\prime} \ar[r]_{h\ast Y^{\prime\prime\prime}} &
X' \ast Y^{\prime\prime\prime} \ar@{=}[d] \\
X \ast Y^{\prime\prime} \ar[r]_-{h \ast Y^{\prime\prime}}  &
X' \ast Y^{\prime\prime}  \ar[r]_-{X' \ast s} &
X' \ast Y^{\prime\prime\prime}}
\raisebox{-38pt}{$=$}
\xymatrix@C=20pt@R=20pt{
X \ast Y \ar[r]^-{X \ast p} \ar@{=}[d]
\ar@{}[rrd]|-{\Longdownarrow{h\ast p
}} &
X \ast Y' \ar[r]^-{h \ast Y'} &
X' \ast Y' \ar@{=}[d] \\
X \ast Y \ar[r]^-{h \ast Y} \ar[d]_-{X\ast q} 
\ar@{}[rd]|-{\Longdownarrow{h\ast q}} &
X' \ast Y \ar[r]^-{X' \ast p} \ar[d]|-{X'\ast q\ }
\ar@{}[rd]|-{\Longdownarrow{X'\ast \omega}} &
X' \ast Y' \ar[d]^-{X' \ast r} \\
X \ast Y^{\prime\prime} \ar[r]_-{h \ast Y^{\prime\prime}}  &
X' \ast Y^{\prime\prime}  \ar[r]_-{X' \ast s} &
X' \ast Y^{\prime\prime\prime} }
$$
$$
\xymatrix@C=20pt@R=20pt{
X \ast Y \ar@{=}[r] \ar[d]_-{v\ast Y} 
\ar@{}[rd]|-{\Longdownarrow{v \ast q}} &
X \ast Y \ar[r]^-{X \ast p} \ar[d]|-{X\ast q} 
\ar@{}[rd]|-{\Longdownarrow{X\ast \omega}} &
X \ast Y' \ar[d]^-{X \ast r} \\
X' \ast Y \ar[d]_-{X'\ast q} &
X \ast Y^{\prime\prime} \ar[r]^-{X \ast s} \ar[d]_-{v \ast Y^{\prime\prime} }
\ar@{}[rd]|-{\Longdownarrow{v \ast s}} &
X \ast Y^{\prime\prime\prime}
\ar[d]^-{v \ast Y^{\prime\prime\prime}} \\
X' \ast Y^{\prime\prime} \ar@{=}[r] &
X' \ast Y^{\prime\prime} \ar[r]_-{X'\ast s} &
X' \ast Y^{\prime\prime\prime}}
\raisebox{-38pt}{$=$}
\xymatrix@C=20pt@R=20pt{
X \ast Y \ar[r]^-{X \ast p} \ar[d]_-{v\ast Y} 
\ar@{}[rd]|-{\Longdownarrow{v \ast p}} &
X \ast Y' \ar@{=}[r]  \ar[d]|-{\ v\ast Y'} 
\ar@{}[rd]|-{\Longdownarrow{v\ast r}} &
X \ast Y' \ar[d]^-{X \ast r} \\
X' \ast Y \ar[r]^-{X' \ast p} \ar[d]_-{X'\ast q} 
\ar@{}[rd]|-{\Longdownarrow{X'\ast \omega}} &
X' \ast Y' \ar[d]^-{X' \ast r}  &
X \ast Y^{\prime\prime\prime} \ar[d]^-{v\ast Y^{\prime\prime\prime}}  \\
X' \ast Y^{\prime\prime} \ar[r]_-{X'\ast s} &
X' \ast Y^{\prime\prime\prime} \ar@{=}[r] &
X' \ast Y^{\prime\prime\prime} }
$$
$$
\xymatrix@C=20pt@R=20pt{
Y \ast X \ar[r]^-{Y\ast h} \ar[d]_-{q\ast X}
\ar@{}[rd]|-{\Longdownarrow{q\ast h}} &
Y \ast X' \ar[r]^-{p\ast X'} \ar[d]|-{q\ast X'}
\ar@{}[rd]|-{\Longdownarrow{\omega \ast X'}} &
Y' \ast X' \ar[d]^-{r\ast X'} \\
Y^{\prime\prime} \ast X \ar[r]_-{Y^{\prime\prime} \ast h} \ar@{=}[d] 
\ar@{}[rrd]|-{\Longdownarrow{s\ast h}} &
Y^{\prime\prime} \ast X' \ar[r]_-{s\ast X'} &
Y^{\prime\prime\prime} \ast X' \ar@{=}[d] \\
Y^{\prime\prime}  \ast X \ar[r]_-{s\ast X} &
Y^{\prime\prime\prime} \ast X \ar[r]_-{Y^{\prime\prime\prime} \ast h} &
Y^{\prime\prime\prime} \ast X'}
\raisebox{-38pt}{$=$}
\xymatrix@C=20pt@R=20pt{
Y \ast X \ar[r]^-{Y\ast h} \ar@{=}[d] 
\ar@{}[rrd]|-{\Longdownarrow{p \ast h}} &
Y \ast X' \ar[r]^-{p\ast X'} &
Y' \ast X' \ar@{=}[d] \\
Y \ast X \ar[r]^-{p\ast X} \ar[d]_-{q\ast X}
\ar@{}[rd]|-{\Longdownarrow{\omega \ast X}} &
Y' \ast X \ar[r]^-{Y'\ast h} \ar[d]|-{r\ast X}
\ar@{}[rd]|-{\Longdownarrow{r\ast h}} &
Y' \ast X' \ar[d]^-{r\ast X'} \\
Y^{\prime\prime}  \ast X \ar[r]_-{s\ast X} &
Y^{\prime\prime\prime} \ast X \ar[r]_-{Y^{\prime\prime\prime} \ast h} &
Y^{\prime\prime\prime} \ast X'}
$$
$$
\xymatrix@C=20pt@R=20pt{
Y \ast X \ar@{=}[r] \ar[d]_-{q\ast X}
\ar@{}[rd]|-{\Longdownarrow{q\ast v}} &
Y \ast X \ar[r]^-{p\ast X} \ar[d]|-{Y \ast v} 
\ar@{}[rd]|-{\Longdownarrow{p\ast v}} &
Y' \ast X \ar[d]^-{Y' \ast v} \\
Y^{\prime\prime}  \ast X \ar[d]_-{Y^{\prime\prime} \ast v} &
Y \ast X'  \ar[r]^-{p\ast X'} \ar[d]_-{q\ast X'}
\ar@{}[rd]|-{\Longdownarrow{\omega \ast X'}} &
Y' \ast X' \ar[d]^-{r\ast X'} \\
Y^{\prime\prime}  \ast X' \ar@{=}[r]  &
Y^{\prime\prime}  \ast X' \ar[r]_-{s\ast X'} &
Y^{\prime\prime\prime} \ast X'}
\raisebox{-38pt}{$=$}
\xymatrix@C=20pt@R=20pt{
Y \ast X \ar[r]^-{p\ast X} \ar[d]_-{q\ast X}
\ar@{}[rd]|-{\Longdownarrow{\omega \ast X}} &
Y' \ast X \ar@{=}[r] \ar[d]|-{r\ast X}
\ar@{}[rd]|-{\Longdownarrow{r\ast v}} &
Y' \ast X \ar[d]^-{Y' \ast v} \\
Y^{\prime\prime}  \ast X \ar[d]_-{Y^{\prime\prime} \ast v} \ar[r]^-{s\ast X}
\ar@{}[rd]|-{\Longdownarrow{s\ast v}} &
Y^{\prime\prime\prime} \ast X \ar[d]^-{Y^{\prime\prime\prime} \ast v} &
Y' \ast X' \ar[d]^-{r\ast X'} \\
Y^{\prime\prime}  \ast X' \ar[r]_-{s\ast X'} &
Y^{\prime\prime\prime} \ast X' \ar@{=}[r]  &
Y^{\prime\prime\prime} \ast X'}
$$
\end{enumerate}

In part (vii) we see the naturality conditions on the horizontal and vertical pseudotransformations obtained as the images of 1-cells under the double functor $M:\mathbb A \to \ldb \mathbb A,\mathbb A\rdb$; and the compatibility conditions on the modifications obtained as the images of 2-cells under $M$.
The conditions in (i-v-vi) come from two sources: from the requirements that $M$, and the image of any 0-cell under it, are double functors; and from the functoriality conditions on the horizontal and vertical pseudotransformations obtained as the images of 1-cells under $M$.
Condition (i) expresses the unitality, and (iii-iv) express the associativity of the monoid. 

There seems to be no evident way to interpret a monoid $\mathbb A$ in $(\mathsf{DblCat},\otimes)$ like above as a monoidal double category in the sense of \cite[Definition 2.9]{Shulman} and \cite[Section 3.1]{GrandisPare:Intercategories_fw}. Using the notation of this section, one can easily introduce multiplication maps $\mathbb A\times \mathbb A \to \mathbb A$ on the various cells. They are associative with the unit $I$:
\begin{itemize}
\item A pair of \underline{0-cells} $X$ and $Y$ is sent to $X\ast Y$.
\item A pair of \underline{horizontal 1-cells} $h$ and $k$ is sent to the horizontal 1-cell on the left; and a pair of \underline{vertical 1-cells} $f$ and $g$ is sent to the vertical 1-cell on the right:
$$
\xymatrix{
X \ast Y \ar[r]^-{h\ast Y} &
X' \ast Y \ar[r]^-{X'\ast k} &
X' \ast Y'}
\qquad \qquad
\raisebox{30pt}{$\xymatrix@R=15pt{
X \ast Y \ar[d]^-{f\ast Y} \\
X' \ast Y \ar[d]^-{X'\ast g} \\
X' \ast Y'.}$}
$$
\item A pair of \underline{2-cells} $\omega$ and $\vartheta$ is sent to
$$
\xymatrix{
X \ast Y \ar[r]^-{h\ast Y} \ar[d]_-{f\ast Y} 
\ar@{}[rd]|-{\Longdownarrow {\omega \ast Y}} &
X' \ast Y \ar[r]^-{X'\ast p} \ar[d]|-{g\ast Y}
\ar@{}[rd]|-{\Longdownarrow {g \ast p}} &
X' \ast Y' \ar[d]^-{g\ast Y'} \\
X^{\prime \prime} \ast Y \ar[r]^-{k\ast Y} \ar[d]_-{X^{\prime \prime} \ast q}
\ar@{}[rd]|-{\Longdownarrow {k \ast q}} &
X^{\prime \prime \prime} \ast Y \ar[d]|-{X^{\prime \prime \prime} \ast q}
\ar[r]^-{X^{\prime \prime \prime} \ast p}
\ar@{}[rd]|-{\Longdownarrow {X^{\prime \prime \prime} \ast \vartheta}} &
X^{\prime \prime \prime} \ast Y' \ar[d]^-{X^{\prime \prime \prime} \ast r} \\
X^{\prime \prime} \ast Y^{\prime \prime} \ar[r]_-{k\ast Y^{\prime \prime} } &
X^{\prime \prime \prime} \ast Y^{\prime \prime} 
\ar[r]_-{X^{\prime \prime \prime} \ast s} &
X^{\prime \prime \prime} \ast Y^{\prime \prime \prime}.}
$$
\end{itemize}
However, these maps do {\em not} constitute a double functor or at least a pseudo double functor in the sense of \cite[Section 2.1]{GrandisPare:dbladj}, \cite[Definition 2.7]{Shulman}. 
Recall that a {\em pseudo double functor} in these references is defined to strictly preserve the composition in one direction; and up-to a coherent natural family of invertible 2-cells in the other direction.
These maps above, however, do not preserve any of the horizontal and vertical compositions in the strict sense, only up-to coherent natural families 
$$
\xymatrix@C=15pt{
X \ast Y \ar[r]^-{h\ast Y} \ar@{=}[d] \ar@{}[rd]|-{\Longdownarrow 1} &
X' \ast Y \ar[r]^-{X' \ast p} \ar@{=}[d] \ar@{}[rrd]|-{\Longdownarrow {k\ast p}} &
X' \ast Y' \ar[r]^-{k\ast Y'} &
X^{\prime \prime} \ast Y' \ar[r]^-{X^{\prime \prime} \ast s}
\ar@{=}[d] \ar@{}[rd]|-{\Longdownarrow 1} &
X^{\prime \prime} \ast Y^{\prime \prime} \ar@{=}[d] \\
X \ast Y \ar[r]_-{h\ast Y} &
X' \ast Y \ar[r]_-{k\ast Y} &
X^{\prime \prime} \ast Y \ar[r]_-{X^{\prime \prime} \ast p} &
X^{\prime \prime} \ast Y' \ar[r]_-{X^{\prime \prime} \ast s} &
X^{\prime \prime} \ast Y^{\prime \prime}}
\qquad
\raisebox{42pt}{$
\xymatrix@R=15pt@C=12pt{
X \ast Y \ar[d]_-{f\ast Y} \ar@{=}[r] \ar@{}[rd]|-{\Longdownarrow 1} &
X \ast Y \ar[d]^-{f\ast Y} \\
X' \ast Y \ar[d]_-{g\ast Y} \ar@{=}[r] \ar@{}[rd]|-{\Longdownarrow {g\ast q}} &
X' \ast Y \ar[d]^-{X' \ast q} \\
X^{\prime \prime} \ast Y \ar[d]_-{X^{\prime \prime} \ast q} &
X' \ast Y' \ar[d]^-{g\ast Y'} \\
X^{\prime \prime} \ast Y' \ar[d]_-{X^{\prime \prime} \ast r} \ar@{=}[r]
\ar@{}[rd]|-{\Longdownarrow 1} &
X^{\prime \prime} \ast Y' \ar[d]^-{X^{\prime \prime} \ast r} \\
X^{\prime \prime} \ast Y^{\prime \prime} \ar@{=}[r] &
X^{\prime \prime} \ast Y^{\prime \prime}}$}
$$
of vertically, respectively, horizontally invertible 2-cells.

Succinctly, monoids in $(\mathsf{DblCat},\otimes)$ determine monoids in the Cartesian monoidal category whose objects are double categories and whose morphisms are pseudo-pseudo double functors (rather than strict--pseudo double functors in \cite{GrandisPare:dbladj}).
By this reason, there seems to be no easy way to regard a monoid in $(\mathsf{DblCat},\otimes)$ in this section as a suitably degenerate {\em intercategory} \cite{GrandisPare:Intercategories,GrandisPare:Intercategories_fw}. 

\subsection{Monoidality of the functor $\Mnd:(\mathsf {DblCat},\otimes) \to (\mathsf {DblCat},\otimes)$ due to Fiore, Gambino and Kock \cite{FioreGambinoKock}}

To any double category $\mathbb A$, the double category $\Mnd(\mathbb A)$ of monads in $\mathbb A$ was associated in \cite{FioreGambinoKock}. This construction can be seen as the object map of the functor in the title, which sends a morphism; that is, a double functor $\mathsf F:\mathbb A \to \mathbb B$ to the double functor $\Mnd(\mathsf F):\Mnd(\mathbb A) \to \Mnd(\mathbb B)$ of `componentwise' action.

As the nullary part of the candidate monoidal structure, we take the evident iso double functor
 $\xymatrix@C=12pt{\mathbbm 1 \ar[r]^-\cong & \Mnd(\mathbbm 1)}$. For any double categories $\mathbb A$ and $\mathbb B$, for the double functors $
\chi_{\mathbb A,\mathbb B}:
\Mnd \ldb \mathbb A,\mathbb B \rdb \to \ldb \Mnd(\mathbb A), \Mnd(\mathbb B) \rdb$ encoding the binary part, we propose the following.

A \underline{0-cell} in $\Mnd \ldb \mathbb A,\mathbb B \rdb$ is by definition a monad 
$(
(\xymatrix@C=12pt{\mathbb A \ar[r]^-{\mathsf T} & \mathbb B},
\xymatrix@C=12pt{{\mathsf T} \ar[r]^-t & {\mathsf T}}),
\theta,\tau)$
in the horizontal 2-category of $\ldb \mathbb A,\mathbb B \rdb$. We have to associate to it a 0-cell in $\ldb \Mnd(\mathbb A), \Mnd(\mathbb B) \rdb$; that is, a double functor $\chi_{\mathbb A,\mathbb B}((\mathsf T,t),\theta,\tau):\Mnd(\mathbb A)\to  \Mnd(\mathbb B)$.

Evaluation at any 0-cell $X$ of $\mathbb A$ gives a 0-cell $((\mathsf T X,t_X),\theta_X,\tau_X)$ in $\Mnd (\mathbb B)$. 
The image of any 0-cell $((X,x),\mu,\eta)$ of $\Mnd(\mathbb A)$ under the double functor $\Mnd(\mathsf T)$ is a 0-cell $((\mathsf TX, \mathsf T x),\mathsf T \mu,\mathsf T \eta)$ in $\Mnd (\mathbb B)$. 
Between these monads in the horizontal 2-category of $\mathbb B$, there is a distributive law $t^x$. It induces a 0-cell in $\Mnd (\mathbb B)$,
$$
\Bigl(\! (\mathsf TX,\!
\xymatrix@C=12pt{\mathsf TX \ar[r]^-{t_X} & \mathsf TX 
\ar[r]^-{\mathsf Tx} & \mathsf TX}\!\!)\! ,\!
\raisebox{38pt}{$\xymatrix@C=12pt{
\mathsf TX \ar[r]^-{t_X} \ar@{=}[d] \ar@{}[rd]|-{\Longdownarrow 1} &
\mathsf TX \ar[r]^-{\mathsf Tx} \ar@{=}[d] \ar@{}[rrd]|-{\Longdownarrow {t^x}} &
\mathsf TX \ar[r]^-{t_X} &
\mathsf TX \ar[r]^-{\mathsf Tx} \ar@{=}[d] \ar@{}[rd]|-{\Longdownarrow 1} &
\mathsf TX  \ar@{=}[d] \\
\mathsf TX \ar[r]^-{t_X} \ar@{=}[d] \ar@{}[rrd]|-{\Longdownarrow {\theta_X}} &
\mathsf TX \ar[r]^-{t_X} &
\mathsf TX \ar[r]^-{\mathsf Tx} \ar@{=}[d] 
\ar@{}[rrd]|-{\Longdownarrow {\mathsf T \mu}} &
\mathsf TX \ar[r]^-{\mathsf Tx} &
\mathsf TX  \ar@{=}[d] \\
\mathsf TX \ar[rr]_-{t_X} &&
\mathsf TX \ar[rr]_-{\mathsf Tx} &&
\mathsf TX}$}\! , \!
\raisebox{38pt}{$\xymatrix@C=12pt@R=63pt{
\mathsf TX  \ar@{=}[d] \ar@{=}[r] \ar@{}[rd]|-{\Longdownarrow {\tau_X}} &
\mathsf TX  \ar@{=}[d] \ar@{=}[r] \ar@{}[rd]|-{\Longdownarrow {\mathsf T\eta}} &
\mathsf TX  \ar@{=}[d] \\
\mathsf TX \ar[r]_-{t_X} &
\mathsf TX \ar[r]_-{\mathsf Tx} &
\mathsf TX}$}\! \Bigr).
$$
It will be the image of the \uuline{0-cell} $((X,x),\mu,\eta)$ of $\Mnd(\mathbb A)$ under the double functor $\chi_{\mathbb A,\mathbb B}((\mathsf T,t),\theta,\tau):\Mnd(\mathbb A)\to  \Mnd(\mathbb B)$.
On the \uuline{horizontal 1-cells} $\chi_{\mathbb A,\mathbb B}((\mathsf T,t),\theta,\tau)$ acts as 
$$
\Bigl(\xymatrix@C=15pt{X \ar[r]^-f & Y},
\raisebox{38pt}{$\xymatrix@C=15pt@R=63pt{
X \ar[r]^-f \ar@{=}[d] \ar@{}[rrd]|-{\Longdownarrow \varphi} &
Y \ar[r]^-y &
Y \ar@{=}[d] \\
X \ar[r]_-x &
X \ar[r]_-f &
Y}$}
\Bigr)
\mapsto
\Bigl(\xymatrix@C=15pt{
\mathsf T X \ar[r]^-{\mathsf T f} & \mathsf T Y},
\raisebox{38pt}{$\xymatrix@C=15pt{
\mathsf T X \ar[r]^-{\mathsf T f} \ar@{=}[d]
\ar@{}[rrd]|-{\Longdownarrow {t^f}} &
\mathsf T Y \ar[r]^-{ t_Y} &
\mathsf T Y \ar[r]^-{\mathsf Ty} \ar@{=}[d] 
 \ar@{}[rd]|-{\Longdownarrow 1} &
\mathsf T Y \ar@{=}[d] \\
\mathsf TX \ar[r]^-{t_X} \ar@{=}[d] \ar@{}[rd]|-{\Longdownarrow 1} &
\mathsf T X \ar[r]^-{\mathsf T f} \ar@{=}[d]
\ar@{}[rrd]|-{\Longdownarrow  {\mathsf T \varphi}} &
\mathsf T Y \ar[r]^-{\mathsf Ty} &
\mathsf T Y \ar@{=}[d] \\
\mathsf TX \ar[r]_-{t_X} &
\mathsf TX \ar[r]_-{\mathsf T x} &
\mathsf TX \ar[r]_-{\mathsf Tf} &
\mathsf T Y}$}\Bigr)
$$
and on the \uuline{vertical 1-cells} it acts as
$$
\Bigl(\raisebox{18pt}{$\xymatrix{X \ar[d]^-g \\ Z}$},
\raisebox{18pt}{$\xymatrix{
X \ar[r]^-x \ar[d]_-g \ar@{}[rd]|-{\Longdownarrow \gamma} &
X \ar[d]^-g \\
Z \ar[r]_-z &
Z}$}\Bigr)
\mapsto
\Bigl(\raisebox{18pt}{$\xymatrix{\mathsf TX \ar[d]^-{\mathsf Tg} \\ \mathsf TZ}$},
\raisebox{18pt}{$\xymatrix{
\mathsf TX \ar[d]_-{\mathsf Tg} \ar[r]^-{t_X} \ar@{}[rd]|-{\Longdownarrow {t_g}} &
\mathsf TX \ar[d]|-{\mathsf Tg} \ar[r]^-{\mathsf Tx} 
\ar@{}[rd]|-{\Longdownarrow {\mathsf T \gamma}} &
\mathsf TX \ar[d]^-{\mathsf Tg} \\
\mathsf TZ \ar[r]_-{t_Z} &
\mathsf TZ \ar[r]_-{\mathsf T z} &
\mathsf TZ}$}\Bigr).
$$
Finally, $\chi_{\mathbb A,\mathbb B}((\mathsf T,t),\theta,\tau)$ sends a \uuline{2-cell} $\omega$ in $\Mnd(\mathbb A)$ to the 2-cell $\mathsf T \omega$ of $\Mnd(\mathbb B)$.

A \underline{horizontal 1-cell} in $\Mnd \ldb \mathbb A,\mathbb B \rdb$ is a monad morphism 
$$
\Bigl(
\xymatrix@C=15pt{\mathsf T \ar[r]^-p & \mathsf T'},
\raisebox{18pt}{$\xymatrix@C=15pt{
\mathsf T \ar[r]^-p \ar@{=}[d] \ar@{}[rrd]|-{\Longdownarrow \pi} &
\mathsf T' \ar[r]^-{t'} &
\mathsf T' \ar@{=}[d]  \\
\mathsf T \ar[r]_-t &
\mathsf T \ar[r]_-p &
\mathsf T'}$}\Bigr)
$$
in the horizontal 2-category of $\ldb \mathbb A,\mathbb B\rdb$ (so that in particular $p$ is a horizontal pseudotransformation and $\pi$ is a modification). The double functor $\chi_{\mathbb A,\mathbb B}$ should send it to a horizontal 1-cell in $\ldb \Mnd(\mathbb A),\Mnd(\mathbb B)\rdb$; that is, the following horizontal pseudotransformation 
$\chi_{\mathbb A,\mathbb B}((T,t),\theta,\tau)\to 
\chi_{\mathbb A,\mathbb B}((T',t'),\theta',\tau')$.

It consists of the horizontal 1-cell in $\Mnd(\mathbb B)$
$$
\Bigl(
\xymatrix@C=15pt{\mathsf T X \ar[r]^-{p_X} & \mathsf {T'} X},
\raisebox{39pt}{$\xymatrix@C=15pt{
\mathsf T X \ar[r]^-{p_X} \ar@{=}[d] 
\ar@{}[rrd]|-{\Longdownarrow {\pi_X}} &
\mathsf {T'} X \ar[r]^-{t'_X} &
\mathsf {T'} X  \ar[r]^-{\mathsf {T'}x} \ar@{=}[d]  \ar@{}[rd]|-{\Longdownarrow 1} &
\mathsf {T'} X \ar@{=}[d] \\
\mathsf TX \ar@{=}[d] \ar[r]^-{t_X} \ar@{}[rd]|-{\Longdownarrow 1} &
\mathsf T X \ar[r]^-{p_X} \ar@{=}[d] 
\ar@{}[rrd]|-{\Longdownarrow {(p^x)^{-1}}} &
\mathsf {T'} X  \ar[r]^-{\mathsf{T'}x} &
\mathsf {T'} X \ar@{=}[d] \\
\mathsf TX \ar[r]_-{t_X} &
\mathsf T X  \ar[r]_-{\mathsf Tx}  &
\mathsf TX \ar[r]_-{p_X} &
\mathsf {T'}X
}$}\Bigr)
$$
for all 0-cells $((X,x),\mu,\eta)$ in $\Mnd(\mathbb A)$ together with the 2-cell $p_g$ in $\Mnd(\mathbb B)$ for all vertical 1-cells $g$ in $\Mnd(\mathbb A)$, and the vertically invertible 2-cell $p^h$ in $\Mnd(\mathbb B)$ for all horizontal 1-cells $h$ in $\Mnd(\mathbb A)$.

A \underline{vertical 1-cell}
$$
\Bigl(
\raisebox{18pt}{$\xymatrix{\mathsf T \ar[d]^-r \\ \mathsf T'}$},
\raisebox{18pt}{$\xymatrix{
\mathsf T \ar[d]_-r \ar[r]^-t \ar@{}[rd]|-{\Longdownarrow \varrho} &
\mathsf T \ar[d]^-r  \\
\mathsf {T'} \ar[r]_-{t'} &
\mathsf {T'}}$}\Bigr)
$$
of $\Mnd \ldb \mathbb A,\mathbb B \rdb$ is sent by $\chi_{\mathbb A,\mathbb B}$ to the vertical pseudotransformation from $\chi_{\mathbb A,\mathbb B}((T,t),\theta,\tau)$ to $\chi_{\mathbb A,\mathbb B}((T',t'),\theta',\tau')$ which consists of the following data.

For all 0-cells $((X,x),\mu,\eta)$ of $\Mnd(\mathbb A)$ the vertical 1-cell in $\Mnd(\mathbb B)$
$$
\Bigl(
\raisebox{18pt}{$\xymatrix{\mathsf TX \ar[d]^-{r_X} \\ \mathsf {T'}X}$},
\raisebox{18pt}{$\xymatrix{
\mathsf TX \ar[d]_-{r_X} \ar[r]^-{t_X} \ar@{}[rd]|-{\Longdownarrow {\varrho_X}} &
\mathsf TX \ar[d]|-{r_X}  \ar[r]^-{\mathsf Tx} \ar@{}[rd]|-{\Longdownarrow {r_x}} &
\mathsf TX \ar[d]^-{r_X} \\
\mathsf {T'}X \ar[r]_-{t'_X} &
\mathsf {T'}X \ar[r]_-{\mathsf T'x} &
\mathsf {T'}X}$}\Bigr);
$$
for all horizontal 1-cells $h$ of $\Mnd(\mathbb A)$ the 2-cell $r_h$ in $\Mnd(\mathbb B)$; and for all vertical 1-cells $g$ of $\Mnd(\mathbb A)$ the horizontally invertible 2-cell $r^g$ in $\Mnd(\mathbb B)$.

The double functor $\chi_{\mathbb A,\mathbb B}$ sends a \underline{2-cell} $\omega$ in $\Mnd \ldb \mathbb A,\mathbb B \rdb$ to the modification whose component at every 0-cell $((X,x),\mu,\eta)$ of $\Mnd(\mathbb A)$ is $\omega_X$.

These double functors $\chi_{\mathbb A,\mathbb B}:
\Mnd \ldb \mathbb A,\mathbb B \rdb \to \ldb \Mnd(\mathbb A), \Mnd(\mathbb B) \rdb$ constitute a natural transformation $\chi:\Mnd \ldb -,-\rdb\to \ldb  \Mnd(-), \Mnd(-)\rdb$.
It is not hard (although a bit long) to see that the double functors of Section \ref{sec:ox}, the trivial isomorphism $\mathbbm 1 \cong \Mnd(\mathbbm 1)$ and the double functors $\chi_{\mathbb A,\mathbb B}$ constructed above, satisfy the conditions of \eqref{eq:unit} and $\eqref{eq:assoc}$. This proves the monoidality of the functor in the title of the section (which sends then monoids as in Section \ref{sec:monoid} to monoids in the same sense).

\subsection{The closed monoidal category $(\mathsf {DblCat},\times)$}
\label{sec:x}

Recall that for any double category $\mathbb A$, the internal hom functor $\langlebar \mathbb A,-\ranglebar : \mathsf{DblCat} \to \mathsf{DblCat}$ of the closed monoidal category in the title sends an object; that is, a double category $\mathbb B$ to the following double subcategory $\langlebar \mathbb A,\mathbb B \ranglebar$ of $\ldb \mathbb A,\mathbb B \rdb$.
\begin{itemize}
\item
The \underline{0-cells} are still the double functors $
\mathbb A \to \mathbb B$.
\item
The \underline{horizontal 1-cells} are the {\em horizontal transformations} of \cite{GrandisPare}. That is, those horizontal pseudotransformations $x$ (see Section \ref{sec:[A,B]}) whose 2-cell parts $x^h$ are vertical identity 2-cells for all horizontal 1-cells $h$ in $\mathbb A$.
\item
Symmetrically, the \underline{vertical 1-cells} are the {\em vertical transformations} of \cite{GrandisPare}. That is, those vertical pseudotransformations $y$
(see Section \ref{sec:[A,B]}) whose 2-cell parts $y^f$ are horizontal identity 2-cells for all vertical 1-cells $f$ in $\mathbb A$.
\item
Finally, the \underline{2-cells} are the modifications of \cite{GrandisPare} (this is the same notion as in Section \ref{sec:[A,B]}).
\end{itemize}
The functor $\langlebar \mathbb A,-\ranglebar : \mathsf{DblCat} \to \mathsf{DblCat}$ sends a morphism; that is, a double functor $\mathsf H:\mathbb B \to \mathbb C$ to the restriction 
$\langlebar \mathbb A,\mathsf H\ranglebar:
\langlebar \mathbb A,\mathbb B \ranglebar \to 
\langlebar \mathbb A,\mathbb C \ranglebar$ 
of the double functor 
$\ldb \mathbb A,\mathsf H\rdb:
\ldb \mathbb A,\mathbb B \rdb \to 
\ldb \mathbb A,\mathbb C \rdb$
in Section \ref{sec:[-,-]}.

For any double category $\mathbb A$, the double functor
$$
\xymatrix{
\mathbbm 1 \ar[r]^-{\eta^{\mathbb A}_\times} &
\langlebar \mathbb A,\mathbbm 1 \times \mathbb A\ranglebar 
\ar[r]^-{\langlebar 1,\lambda_\times\ranglebar} &
\langlebar \mathbb A,\mathbb A\ranglebar}
$$
is $1_{\mathbb A}$, sending the single object of $\mathbbm 1$ to the identity double functor $1_{\mathbb A}:\mathbb A \to \mathbb A$; and the double functor
$$
\xymatrix{
\mathbb A \ar[r]^-{\eta^{\mathbbm 1}_\times} &
\langlebar \mathbbm 1,\mathbb  A \times \mathbbm 1\ranglebar 
\ar[r]^-{\langlebar 1,\varrho_\times \ranglebar} &
\langlebar \mathbbm 1,\mathbb A\ranglebar}
$$
is the canonical isomorphism $\mathbb A\cong \langlebar \mathbbm 1,\mathbb A\ranglebar \cong \ldb \mathbbm 1,\mathbb A\rdb$. 

For any double categories $\mathbb A$, $\mathbb B$ and $\mathbb C$, the double  functor 
$
\mathfrak l^{\mathbb C}_\times:
\langlebar \mathbb A,\mathbb B \ranglebar \to
\langlebar
\langlebar \mathbb C,\mathbb A\ranglebar,\langlebar \mathbb C,\mathbb B \ranglebar
\ranglebar
$
of \eqref{eq:gen_l} is constituted by the following maps.
\begin{itemize}
\item It sends a \underline{0-cell}; that is, a double functor $\mathsf F:\mathbb A \to \mathbb B$ to the double functor 
$\langlebar \mathbb C,\mathsf F \ranglebar:
\langlebar \mathbb C,\mathbb A\ranglebar \to
\langlebar \mathbb C,\mathbb B \ranglebar$.
\item It sends a \underline{horizontal 1-cell}; that is, a horizontal transformation 
$\xymatrix@C=12pt{\mathsf F \ar[r]^-x & \mathsf G}$
to the horizontal transformation $\langlebar \mathbb C,\mathsf F \ranglebar \to \langlebar \mathbb C,\mathsf G \ranglebar$ whose component at any vertical transformation on the left --- between double functors $\mathbb C \to \mathbb A$ --- is the modification on the right:
$$
\xymatrix{
\mathsf H \ar[d]^-q \\ \mathsf {H'}}
\qquad \qquad
\xymatrix{
\mathsf{FH} \ar[r]^-{x_{\mathsf H-}} \ar[d]_-{\mathsf F q_-} 
\ar@{}[rd]|-{\Longdownarrow {x_{q_-}}} &
\mathsf{GH} \ar[d]^-{\mathsf G q_-} \\
\mathsf{FH'} \ar[r]_-{x_{\mathsf {H'}-}} &
\mathsf{GH'}. }
$$
\item Symmetrically, it sends a \underline{vertical 1-cell}; that is, a vertical transformation on the left --- between double functors $\mathbb A \to \mathbb B$ --- to the vertical transformation whose component at any horizontal transformation 
$ \xymatrix@C=12pt{ \mathsf H \ar[r]^-p &  \mathsf {H'}}$ --- between double functors $\mathbb C \to \mathbb A$ --- is the modification on the right:
$$
\xymatrix{\mathsf F \ar[d]^-y \\ \mathsf J}
\qquad \qquad
\xymatrix{
\mathsf{FH} \ar[d]_-{y_{\mathsf H-}} \ar[r]^-{\mathsf F p_-} 
\ar@{}[rd]|-{\Longdownarrow {y_{p_-}}} &
\mathsf{FH'} \ar[d]^-{y_{\mathsf {H'}-}} \\
\mathsf{JH} \ar[r]_-{\mathsf J p_-} &
\mathsf{JH'}. }
$$
\item Finally, it sends a modification on the left to the modification whose component at a double functor $\mathsf H:\mathbb C \to \mathbb A$ is the modification on the right:
$$
\xymatrix{
\mathsf F \ar[r]^-x \ar[d]_-y \ar@{}[rd]|-{\Longdownarrow \Theta} &
\mathsf G \ar[d]^-v \\
\mathsf J \ar[r]_-z &
\mathsf K}
\qquad \qquad
\xymatrix{
\mathsf {FH} \ar[r]^-{x_{\mathsf H -}} \ar[d]_-{y_{\mathsf H -}}
\ar@{}[rd]|-{\Longdownarrow {\Theta_{\mathsf H-}}} &
\mathsf {GH}  \ar[d]^-{v_{\mathsf H -}} \\
\mathsf {JH}  \ar[r]_-{z_{\mathsf H -}} &
\mathsf {KH}. }
$$
\end{itemize}

\subsection{Monoidality of the identity functor $(\mathsf {DblCat},\times) \to (\mathsf {DblCat},\otimes)$}

The evident inclusion double functors
$\langlebar \mathbb A,\mathbb B \ranglebar
\rightarrowtail
\ldb \mathbb A,\mathbb B \rdb 
$,
for all double categories $\mathbb A$ and $\mathbb B$, define a natural transformation 
$\langlebar -,- \ranglebar
\rightarrowtail
\ldb -,- \rdb 
$.
Together with the double functors in Section \ref{sec:ox} and those in Section \ref{sec:x}, and the identity double functor
$\xymatrix@C=12pt{\mathbbm 1 \ar@{=}[r] & \mathbbm 1}$ 
as the nullary part of the stated monoidal structure, they clearly satisfy the conditions in \eqref{eq:unit} and \eqref{eq:assoc}. With this we infer the monoidality of the functor in the title of the section.

In particular, a strict monoidal double category \cite{BruniMeseguerMontanari} -- which is a monoid in $(\mathsf{DblCat},\times)$ -- gives rise to a monoid in $(\mathsf{DblCat},\otimes)$ -- described in Section \ref{sec:monoid}.

\subsection{The closed monoidal category $(\TwoCat,\otimes)$}

In this section we regard the category $\TwoCat$ of 2-categories and 2-functors as a closed monoidal category via the Gray monoidal product $\otimes$ in \cite{Gray}.
Recall that, for any 2-category $\mathcal A$, the internal hom functor $[\mathcal A,-]:\TwoCat \to \TwoCat$ sends an object; that is, a 2-category $\mathcal B$ to the 2-category $[\mathcal A,\mathcal B]$ of 2-functors $\mathcal A \to \mathcal B$, pseudonatural transformations, and modifications. It sends a morphism; that is, a 2-functor $\mathsf F: \mathcal B \to \mathcal C$ to the 2-functor $[\mathcal A,\mathsf F]:[\mathcal A,\mathcal B] \to [\mathcal A,\mathcal C]$ given by postcomposition with $\mathsf F$.

The 2-functor
$$
\xymatrix{
\calone \ar[r]^-{\eta^{\mathcal A}} &
[\mathcal A,\calone \otimes \mathcal A] \ar[r]^-{[1,\lambda]} &
[\mathcal A,\mathcal A]}
$$
is $1_{\mathcal A}$, the 2-functor sending the single object of the terminal 2-category $\calone$ to the identity 2-functor $1_{\mathcal A}:\mathcal A \to \mathcal A$.

The 2-functor
$$
\xymatrix{
\mathcal A \ar[r]^-{\eta^{\scalebox{.6}{$\calone$}}} &
[\calone,\mathcal A\otimes \calone] \ar[r]^-{[1,\varrho]} &
[\calone,\mathcal A]}
$$
is the canonical isomorphism.

For any 2-categories $\mathcal A$, $\mathcal B$ and $\mathcal C$, the 2-functor $\mathfrak l^{\mathcal C}_{\mathcal A,\mathcal B}:[\mathcal A,\mathcal B] \to [[\mathcal C,\mathcal A], [\mathcal C,\mathcal B]]$ in \eqref{eq:gen_l} has the following maps.
\begin{itemize}
\item It sends a \underline{0-cell}; that is, a 2-functor $\mathsf H:\mathcal A \to \mathcal B$ to the 2-functor $[\mathcal C,\mathsf H]:[\mathcal C,\mathcal A] \to [\mathcal C,\mathcal B]$.
\item It sends a \underline{1-cell}; that is, a pseudonatural transformation $\psi:\mathsf H \to \mathsf {H'}$ to the pseudonatural transformation $[\mathcal C,\mathsf H] \to [\mathcal C,\mathsf H']$ whose component at any pseudonatural transformation $\varphi:\mathsf F \to \mathsf {F'}$ between 2-functors $\mathcal C \to \mathcal A$ is
$$
\xymatrix{
\mathsf{HF} \ar[d]_-{\psi_{\mathsf F-}} \ar[r]^-{\mathsf H \varphi_-}
\ar@{}[rd]|-{\longdownarrow {\psi_{\varphi_-}}} &
\mathsf{HF'} \ar[d]^-{\psi_{\mathsf {F'}-}} \\
\mathsf{H'F} \ar[r]_-{\mathsf {H'} \varphi_-}  &
\mathsf{H'F'}.}
$$
\item It sends a \underline{2-cell}; that is, a modification $\omega$ to the modification whose component at any 2-functor $\mathsf F:\mathcal C \to \mathcal A$ is $\omega_{\mathsf F-}$.
\end{itemize}

\subsection{Monoidality of the functors $(\mathsf {DblCat},\otimes) \to (\TwoCat,\otimes)$ sending double categories to their horizontal -- or vertical -- 2-categories} 

Consider the functor $\mathcal H:\mathsf {DblCat} \to \TwoCat$ which sends a double category $\mathbb A$ to its so-called horizontal 2-category. The 0-cells of $\mathcal H \mathbb A$ are the 0-cells of $\mathbb A$, the 1-cells of $\mathcal H \mathbb A$ are the horizontal 1-cells of $\mathbb A$ and the 2-cells of $\mathcal H \mathbb A$ are those 2-cells of $\mathbb A$ which are surrounded by identity vertical 1-cells (and arbitrary horizontal 1-cells). Compositions in $\mathcal H \mathbb A$ are inherited from $\mathbb A$. The functor $\mathcal H$ sends a morphism; that is, a double functor $\mathsf F: \mathbb A \to \mathbb B$ to the 2-functor $\mathcal H \mathsf F:\mathcal H \mathbb A \to \mathcal H \mathbb B$ which acts on the various cells as $\mathsf F$ does.

The horizontal 2-category of the terminal double category $\mathbbm 1$ is the terminal 2-category $\calone$. So we may choose the nullary part $\mathcal H_0$ of the candidate monoidal structure on $\mathcal H$ to be the identity 2-functor $\calone \to \calone$. As the 2-functor $\chi_{\mathbb A,\mathbb B}:\mathcal H \ldb \mathbb A,\mathbb B \rdb\to [\mathcal H \mathbb A, \mathcal H \mathbb B]$ for any double categories $\mathbb A$ and $\mathbb B$, encoding the binary part,  we propose the following.
\begin{itemize}
\item
A \underline{0-cell}; that is, a double functor $\mathsf F:\mathbb A \to \mathbb B$ is sent to the 2-functor $\mathcal H \mathsf F:\mathcal H \mathbb A \to \mathcal H \mathbb B$.
\item
A \underline{1-cell}; that is, a horizontal pseudotransformation $x:\mathsf F \to \mathsf G$ is sent to the pseudonatural transformation $\mathcal H \mathsf F \to \mathcal H \mathsf G$ whose component at any 1-cell of $\mathcal H \mathbb A$ --- that is, horizontal 1-cell $h:A\to C$ of $\mathbb A$ --- is the 2-cell $x^h:x_C.\mathsf F h\to \mathsf G h.x_A$ of $\mathcal H \mathbb B$.
\item
A \underline{2-cell}; that is, a modification of the form
\raisebox{15pt}{$\xymatrix@C=17pt@R=17pt{
\mathsf F \ar[r]^-x \ar@{=}[d] \ar@{}[rd]|-{\Longdownarrow \omega} &
\mathsf G \ar@{=}[d] \\
\mathsf F \ar[r]_-z &
\mathsf G}$},
is sent to the modification whose component at any 0-cell of $\mathcal H \mathbb A$ --- that is, 0-cell $A$ of $\mathbb A$ --- is the 2-cell $\omega_A:x_A \to z_A$ of $\mathcal H \mathbb B$.
\end{itemize}
The so defined 2-functors $\chi_{\mathbb A,\mathbb B}$ constitute a natural transformation $\chi:\mathcal H \ldb -,- \rdb\to [\mathcal H -, \mathcal H -]$ and satisfy the conditions in \eqref{eq:unit} and \eqref{eq:assoc}. Hence there is a corresponding monoidal structure on $\mathcal H$. In particular, applying $\mathcal H$ to a monoid in $(\mathsf {DblCat},\otimes)$ as in Section \ref{sec:monoid}, we obtain a monoid in $(\TwoCat,\otimes)$ (known as a {\em Gray monoid} \cite{GordonPowerStreet}).

Symmetric considerations verify monoidality of the functor $\mathsf V:\mathsf {DblCat} \to \TwoCat$, sending a double category to its vertical 2-category. 

\subsection{Monoidality of the functor $\Sqr:(\TwoCat,\otimes) \to (\mathsf {DblCat}, \otimes)$ due to Ehresmann}

Ehresmann's {\em square-} or {\em quintet construction} \cite{Ehresmann} sends a 2-category $\mathcal A$ to the following double category $\Sqr(\mathcal A)$.
The 0-cells of $\Sqr(\mathcal A)$ are the 0-cells of $\mathcal A$.
Both the horizontal and the vertical 1-cells of $\Sqr(\mathcal A)$ are the 1-cells of $\mathcal A$.
A 2-cell of $\Sqr(\mathcal A)$ with boundaries 
\raisebox{15pt}{$\xymatrix@C=17pt@R=17pt{
A \ar[r]^-t \ar[d]_-l \ar@{}[rd]|-{\Longdownarrow {}} &
C \ar[d]^-r \\
B \ar[r]_-b &
D}$}
is a 2-cell $r.t \to b.l$ of $\mathcal A$.
For any 2-functor $\mathsf F:\mathcal A \to \mathcal B$ there is a double functor $\Sqr(\mathsf F):\Sqr(\mathcal A) \to \Sqr(\mathcal B)$ which acts on the various cells as $\mathsf F$ does.

Applying the so defined functor $\Sqr$ to the terminal 2-category $\calone$, we obtain the terminal double category $\mathbbm 1$. So as the nullary part of the candidate monoidal structure on $\Sqr$, we may choose the identity double functor $\mathbbm 1 \to \mathbbm 1$. For any 2-categories $\mathcal A$ and $\mathcal B$, for the double functor 
$\Sqr[\mathcal A,\mathcal B] \to \ldb \Sqr(\mathcal A),\Sqr(\mathcal B) \rdb$ encoding the binary part, the following choices can be made.
\begin{itemize}
\item A \underline{0-cell}; that is, a 2-functor $\mathsf H:\mathcal A \to \mathcal B$ is sent to the double functor $\Sqr(\mathsf H):\Sqr(\mathcal A) \to \Sqr(\mathcal B)$.
\item A \underline{horizontal 1-cell}; that is, a pseudonatural transformation $p:\mathsf H \to \mathsf K$ is sent to the horizontal pseudotransformation whose components at any 1-cell $f:X\to Y$ of $\mathcal A$ are the 2-cells in $\Sqr(\mathcal B)$
$$
\xymatrix{
\mathsf H X \ar[r]^-{p_X} \ar[d]_-{\mathsf H f} 
\ar@{}[rd]|-{\Longdownarrow {(p_f)^{-1}}} &
\mathsf {H'}X \ar[d]^-{\mathsf {H'} f} \\
\mathsf H Y \ar[r]_-{p_Y} &
\mathsf {H'} Y}
\qquad \qquad
\xymatrix{
\mathsf H X \ar@{=}[d] \ar[r]^-{\mathsf H f} 
\ar@{}[rrd]|-{\Longdownarrow {p_f}} &
\mathsf H Y \ar[r]^-{p_Y} &
\mathsf {H'} Y \ar@{=}[d] \\
\mathsf H X \ar[r]_-{p_X} &
\mathsf {H'}X \ar[r]_-{\mathsf {H'} f}  &
\mathsf {H'} Y.}
$$
\item A \underline{vertical 1-cell}; which is again a pseudonatural transformation $p:\mathsf H \to \mathsf K$ is sent to the vertical pseudotransformation whose components at any 1-cell $f:X\to Y$ of $\mathcal A$ are the 2-cells in $\Sqr(\mathcal B)$
$$
\xymatrix@R=55pt{
\mathsf H X \ar[d]_-{p_X} \ar[r]^-{\mathsf H f} 
\ar@{}[rd]|-{\Longdownarrow {p_f}} &
\mathsf H Y \ar[d]^-{p_Y} \\
\mathsf {H'}X \ar[r]_-{\mathsf {H'} f} &
\mathsf {H'} Y}
\qquad \qquad
\xymatrix@R=20pt{
\mathsf H X \ar@{=}[r] \ar[d]_-{p_X}
\ar@{}[rdd]|-{\Longdownarrow {p_f}} &
\mathsf H X \ar[d]^-{\mathsf H f} \\
\mathsf {H'}X \ar[d]_-{\mathsf {H'} f}  &
\mathsf H Y \ar[d]^-{p_Y} \\
\mathsf {H'} Y \ar@{=}[r] &
\mathsf {H'} Y.}
$$
\item A \underline{2-cell}; that is, a modification of 2-functors on the left, is sent to the modification of double functors whose component at any 0-cell of $\Sqr(\mathcal A)$ --- that is, any 0-cell $X$ of $\mathcal A$ --- is the 2-cell of $\Sqr(\mathcal B)$ on the right:
$$
\xymatrix{
\mathsf H \ar[r]^-p \ar[d]_-q \ar@{}[rd]|-{\longdownarrow \omega} &
\mathsf {H'} \ar[d]^-r \\
\mathsf {H^{\prime\prime}} \ar[r]_-s &
\mathsf {H^{\prime\prime\prime}} }
\qquad \qquad
\xymatrix{
\mathsf HX \ar[r]^-{p_X} \ar[d]_-{q_X} \ar@{}[rd]|-{\Longdownarrow {\omega_X}} &
\mathsf {H'} X\ar[d]^-{r_X} \\
\mathsf {H^{\prime\prime}} X\ar[r]_-{s_X} &
\mathsf {H^{\prime\prime\prime}} X.}
$$
\end{itemize}
The resulting double functors $\Sqr[\mathcal A,\mathcal B] \to \ldb \Sqr(\mathcal A),\Sqr(\mathcal B) \rdb$ constitute a natural transformation  
$\Sqr[-,-] \to \ldb \Sqr(-),\Sqr(-) \rdb$ and satisfy the conditions in \eqref{eq:unit} and \eqref{eq:assoc}. So they render monoidal the functor in the title of the section. In particular, applying the functor $\Sqr$ to a Gray monoid --- that is, a monoid in $(\TwoCat,\otimes)$ --- a monoid in $(\mathsf {DblCat}, \otimes)$ is obtained (which may not be a monoidal bicategory in the sense of \cite[Definition 2.9]{Shulman}).

%%%%%%%%%%%%%%%%%%%%%%%% APPENDIX   %%%%%%%%%%%%%%%%%%%%%

\appendix
\section{Diagrams}

Large size diagrams, used in the earlier sections, are collected in this appendix, after the bibliography.

%%%%%%%%%%%%%%%%%%%%%%%% BIBLIOGRAPHY   %%%%%%%%%%%%%%%%%%%%%

\bibliographystyle{plain}

\begin{amssidewaysfigure}
\thisfloatpagestyle{empty}
{\color{white} .}
\hspace{-2cm}
\centering
\scalebox{.95}{$
\xymatrix@C=4pt@R=55pt{
\mathsf{DblCat}(\mathbb A \! \otimes \! (\mathbb B \! \otimes \! (\mathbb C \! \otimes \! \mathbb D)),\mathbb K)
\ar[rr]^-{\mathsf{DblCat}(\alpha,\mathbb K)}
\ar[dddd]_-{\mathsf{DblCat}(1 \otimes \alpha,\mathbb K)} &&
\mathsf{DblCat}((\mathbb A \! \otimes \! \mathbb B) \! \otimes \! (\mathbb C\! \otimes \! \mathbb D),\mathbb K)
\ar[rr]^-{\mathsf{DblCat}(\alpha,\mathbb K)} &&
\mathsf{DblCat}(((\mathbb A \! \otimes \! \mathbb B) \! \otimes \! \mathbb C)\! \otimes \! \mathbb D,\mathbb K) \\
% 2.3
&& \mathsf{DblCat}(\mathbb A \! \otimes \! \mathbb B , 
\ldb \mathbb C\! \otimes \! \mathbb D,\mathbb K \rdb)
\ar[r]^-{\raisebox{10pt}{${}_
{\mathsf{DblCat}(\mathbb A \otimes \mathbb B , \mathfrak a^{\mathbb K})}$}}
\ar[u]_-\cong &
\mathsf{DblCat}(\mathbb A \! \otimes \! \mathbb B , 
\ldb \mathbb C,\ldb \mathbb D,\mathbb K \rdb \rdb)
\ar[ru]^-\cong \\
% 3.2
& \mathsf{DblCat}(\mathbb A ,
\ldb \mathbb B \! \otimes \! (\mathbb C \! \otimes \! \mathbb D),\mathbb K\rdb)
\ar[d]_-{\mathsf{DblCat}(\mathbb A , \ldb \alpha,1 \rdb)}
\ar[r]^-{\raisebox{10pt}{${}_
{\mathsf{DblCat}(\mathbb A ,\mathfrak a^{\mathbb K})}$}}
\ar[luu]_-\cong 
\ar@{}[rrd]|-{\textrm{Figure}~ \ref{fig:pent_eq}}&
% 3.3
\mathsf{DblCat}(\mathbb A ,
\ldb \mathbb B ,\ldb \mathbb C \! \otimes \! \mathbb D,\mathbb K\rdb \rdb)
\ar[u]_-\cong 
\ar[r]^-{\raisebox{10pt}{${}_
{\mathsf{DblCat}(\mathbb A ,\ldb 1,\mathfrak a^{\mathbb K} \rdb)}$}} &
\mathsf{DblCat}(\mathbb A ,
\ldb \mathbb B ,\ldb \mathbb C  ,\ldb \mathbb D,\mathbb K\rdb \rdb \rdb)
\ar[u]_-\cong  \\
& \mathsf{DblCat}(\mathbb A ,
\ldb (\mathbb B \! \otimes \! \mathbb C) \! \otimes \! \mathbb D,\mathbb K\rdb)
\ar[rr]_-{\mathsf{DblCat}(\mathbb A ,\mathfrak a^{\mathbb K})} \ar[ld]_-\cong &&
\mathsf{DblCat}(\mathbb A ,
\ldb \mathbb B \otimes \mathbb C ,\ldb  \mathbb D,\mathbb K\rdb \rdb)
\ar[u]_-{\mathsf{DblCat}(\mathbb A ,\mathfrak a^{\ldb  \mathbb D,\mathbb K\rdb})}
\ar[rd]^-\cong \\
\mathsf{DblCat}(\mathbb A \! \otimes \! ((\mathbb B \! \otimes \! \mathbb C) \! \otimes \! \mathbb D),\mathbb K)
\ar[rrrr]_-{\mathsf{DblCat}(\alpha,\mathbb K)}&&&&
\mathsf{DblCat}((\mathbb A \! \otimes \! (\mathbb B \! \otimes \! \mathbb C)) \! \otimes \! \mathbb D,\mathbb K)
\ar[uuuu]_-{\mathsf{DblCat}(\alpha \otimes 1,\mathbb K)}}$}
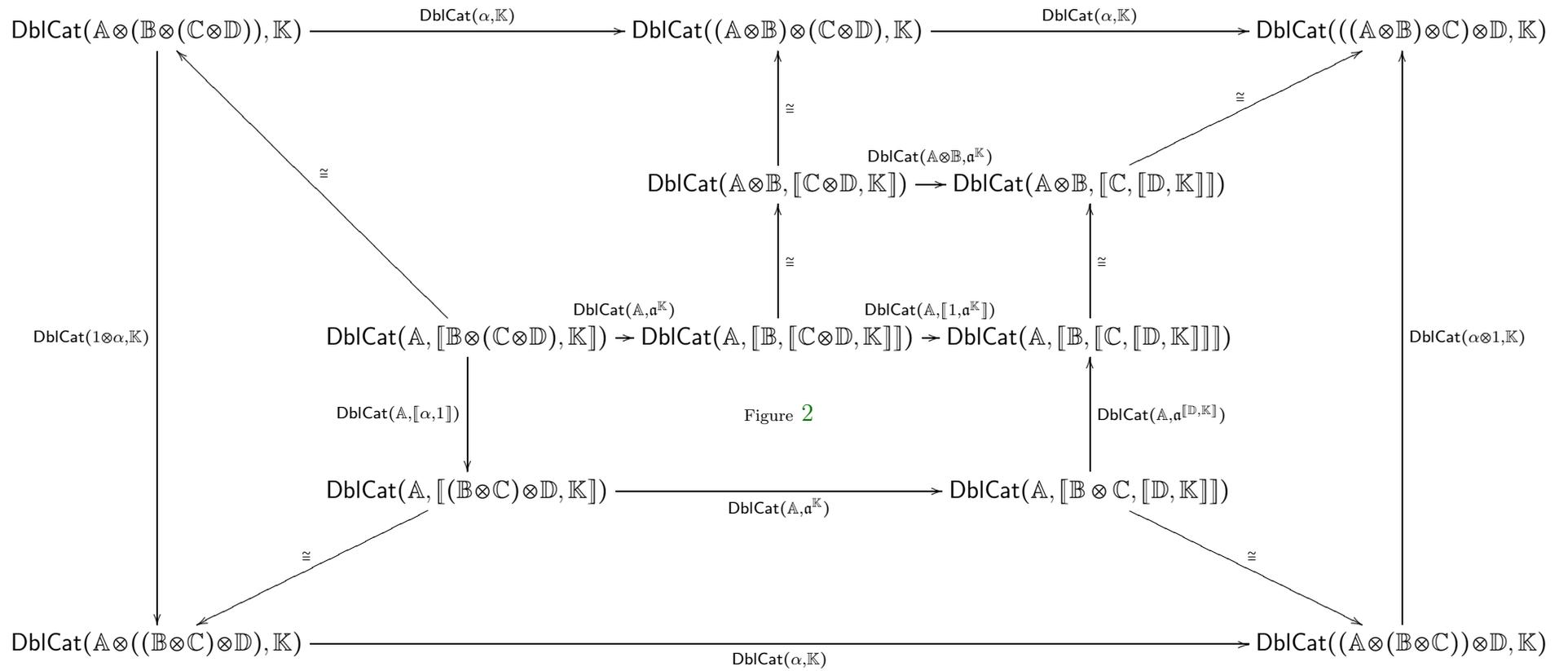
\caption{Equivalent forms of the pentagon condition}
\label{fig:pent}
\end{amssidewaysfigure}

\begin{amssidewaysfigure}
\thisfloatpagestyle{empty}
\centering
{\color{white} .}
\hspace{-3.5cm}
\scalebox{.975}{$
\xymatrix@C=7pt@R=55pt{
% 1.1
\ldb \mathbb A \! \otimes \! (\mathbb B \! \otimes \! \mathbb C)\tp \mathbb K \rdb
\ar@/^3pc/[rrrr]^-{\mathfrak a^{\mathbb K}}
\ar[rrr]_-{\raisebox{-10pt}{${}_{\mathfrak l^{\mathbb B \otimes \mathbb C}}$}}
\ar[rd]^-{\mathfrak l^{\mathbb C}}
\ar[d]_-{\ldb \epsilon^{\mathbb C},1 \rdb}
%\ar@/_9pc/[dddd]_(.25){\ldb \alpha, 1\rdb} 
&&&
% 1.4
\ldb \hspace{-1pt}
   \ldb \mathbb B \! \otimes \! \mathbb C\tp \mathbb A \! \otimes \! (\mathbb B \! \otimes \! \mathbb C) \rdb\tp
   \ldb \mathbb B \! \otimes \! \mathbb C\tp \mathbb K \rdb \hspace{-1pt}
\rdb
\ar[r]_-{\raisebox{-10pt}{${}_{\ldb \eta^{\mathbb B \otimes \mathbb C},1\rdb}$}}
\ar[ddd]_-{\ldb 1, \mathfrak a^{\mathbb K}\rdb} &
% 1.5
\ldb 
   \mathbb A\tp
   \ldb \mathbb B \! \otimes \! \mathbb C\tp \mathbb K \rdb \hspace{-1pt}
\rdb
\ar[dddd]^-{\ldb 1, \mathfrak a^{\mathbb K}\rdb} \\
%2.1
\ldb \hspace{-1pt}
   \ldb \mathbb C\tp \mathbb A \! \otimes \! (\mathbb B \! \otimes \! \mathbb C) \rdb \! \otimes \! \mathbb C\tp
   \mathbb K
\rdb
\ar[r]_-{\mathfrak a^{\mathbb K}}
\ar[d]_-{\ldb \epsilon^{\mathbb B}\otimes 1,1\rdb} 
\ar@{}[ru]|(.35){(\ast)} & 
% 2.2
\ldb \hspace{-1pt}
   \ldb \mathbb C\tp \mathbb A \! \otimes \! (\mathbb B \! \otimes \! \mathbb C) \rdb\tp 
   \ldb \mathbb C\tp  \mathbb K \rdb \hspace{-1pt}
\rdb
\ar[rd]^-{\mathfrak l^{\mathbb B}}
\ar[d]_-{\ldb\epsilon^{\mathbb B},1\rdb} & \\
% 3.1
\ldb \hspace{-1pt}
     ( \hspace{-1pt}
       \ldb
            \mathbb B\tp
            \ldb
                 \mathbb C\tp
                 \mathbb A \! \otimes \! (\mathbb B \! \otimes \!\mathbb C) \hspace{-1pt}
            \rdb \hspace{-1pt}
       \rdb \! \otimes \! \mathbb B
      )\! \otimes \!\mathbb C \tp
      \mathbb K
\rdb
\ar[r]_-{\raisebox{-10pt}{${}_{\mathfrak a^{\mathbb K}}$}}
\ar[d]_-{
\ldb
    (\mathfrak a^{\mathbb A \otimes (\mathbb B \otimes \mathbb C)} \otimes 1) \otimes 1,
    1
\rdb} &
% 3.2
\ldb \hspace{-1pt}
   \ldb
        \mathbb B \tp
        \ldb \mathbb C\tp \mathbb A \! \otimes \!(\mathbb B \! \otimes \! \mathbb C) \rdb \hspace{-1pt}
   \rdb \! \otimes \!\mathbb B\tp
   \ldb \mathbb C \tp \mathbb K \rdb \hspace{-1pt}
\rdb
\ar[r]_-{\raisebox{-10pt}{${}_{\mathfrak a^{\ldb \mathbb C , \mathbb K \rdb}}$}}
\ar[d]^-{
\ldb
    \mathfrak a^{\mathbb A \otimes (\mathbb B  \otimes \mathbb C)}\otimes 1,
    1
\rdb} 
\ar@{}[ru]|(.35){(\ast)} 
\ar@{}[rruu]|(.5){(\ast\ast)} & 
% 3.3
\ldb \hspace{-1pt}
   \ldb
        \mathbb B \tp
        \ldb \mathbb C\tp\mathbb A \! \otimes \! (\mathbb B \! \otimes  \!\mathbb C) \rdb \hspace{-1pt}
   \rdb \tp
   \ldb
         \mathbb B\tp
         \ldb \mathbb C \tp \mathbb K \rdb \hspace{-1pt}
   \rdb\hspace{-1pt}
\rdb
\ar[rd]^-{\quad \ldb \mathfrak a^{\mathbb A \otimes (\mathbb B \otimes \mathbb C)},1\rdb} \\
% 4.1
\ldb \hspace{-1pt}
      ( \hspace{-1pt}
          \ldb
               \mathbb B \!  \otimes  \! \mathbb C \tp
               \mathbb A  \! \otimes \!(\mathbb B \! \otimes \!\mathbb C) \hspace{-1pt}
          \rdb \! \otimes  \!\mathbb B
       ) \! \otimes \! \mathbb C\tp
      \mathbb K
\rdb
\ar[r]^-{\mathfrak a^{\mathbb K}} 
\ar[d]_-{
\ldb
      (\eta^{\mathbb B \otimes \mathbb C} \otimes 1) \otimes 1,
      1
\rdb} &
% 4.2
\ldb \hspace{-1pt}
     \ldb
            \mathbb B \! \otimes \!\mathbb C\tp
            \mathbb A \! \otimes \! (\mathbb B \! \otimes \!\mathbb C)\hspace{-1pt}
     \rdb \! \otimes \! \mathbb B\tp
     \ldb
           \mathbb C\tp
           \mathbb K
     \rdb \hspace{-1pt}
\rdb
\ar[rr]^-{\mathfrak a^{\ldb \mathbb C,\mathbb K\rdb}}
\ar[d]^-{\ldb \eta^{\mathbb B \otimes \mathbb C} \otimes 1,1\rdb } && 
% 4.4
\ldb \hspace{-1pt}
       \ldb
              \mathbb B \! \otimes \! \mathbb C\tp
              \mathbb A \! \otimes \! (\mathbb B \! \otimes \! \mathbb C) \hspace{-1pt}
       \rdb\tp
       \ldb
              \mathbb B \tp
              \ldb \mathbb C\tp\mathbb K\rdb \hspace{-1pt}
       \rdb \hspace{-1pt}
\rdb
\ar[rd]^-{\ldb \eta^{\mathbb B \otimes \mathbb C},1\rdb} \\
% 5.1
\ldb
       (\mathbb A \! \otimes \! \mathbb B) \! \otimes \! \mathbb C \tp
       \mathbb K
\rdb
\ar[r]_-{\mathfrak a^{\mathbb K}} &
% 5.2
\ldb
      \mathbb A \! \otimes \!  \mathbb B \tp
      \ldb \mathbb C\tp \mathbb K\rdb \hspace{-1pt}
\rdb
\ar[rrr]_-{\mathfrak a^{\ldb \mathbb C, \mathbb K\rdb}}
&&& 
% 5.5
\ldb
       \mathbb A\tp
       \ldb
              \mathbb B\tp
              \ldb \mathbb C\tp\mathbb K\rdb \hspace{-1pt}
       \rdb \hspace{-1pt}
\rdb}$}
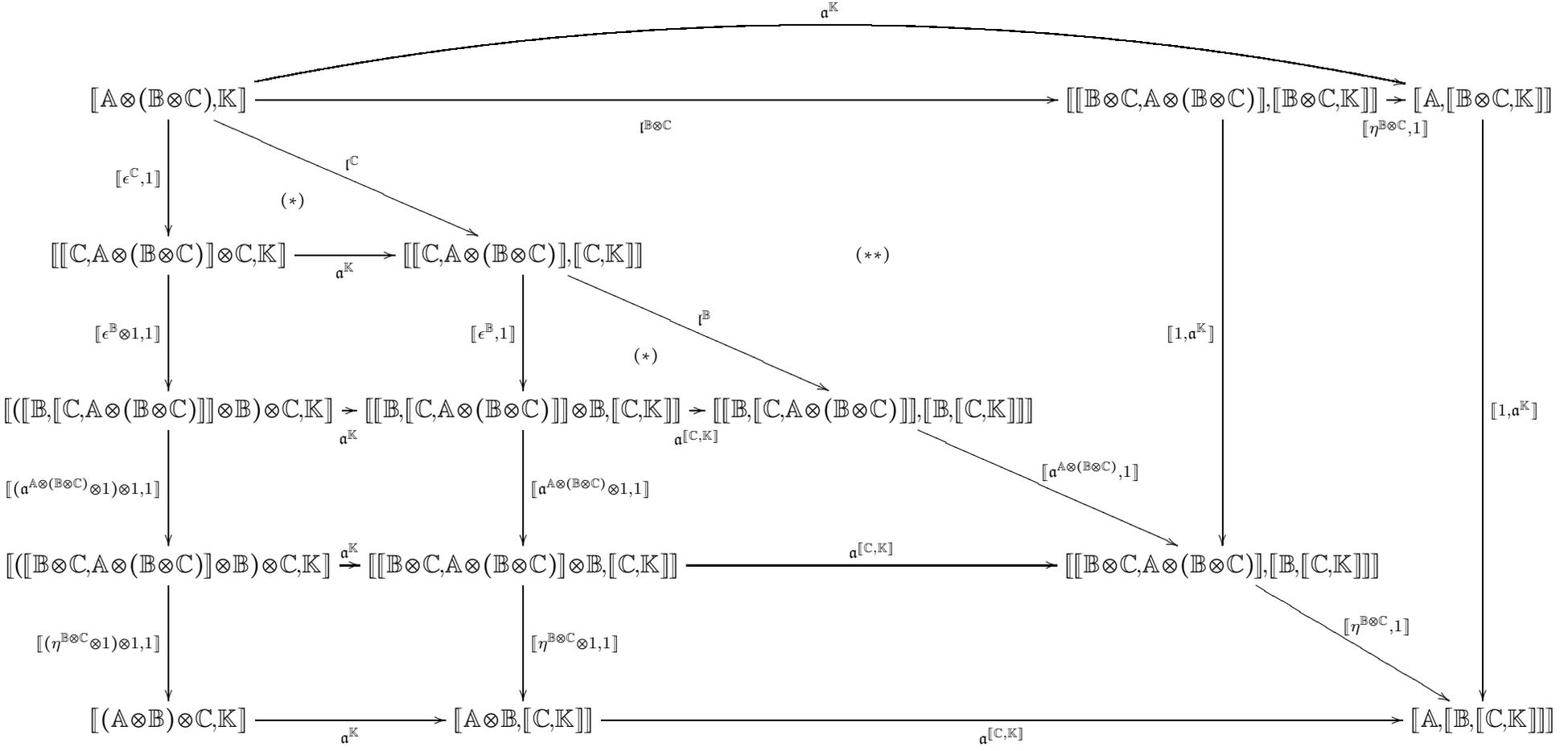
\caption{Proof of the pentagon condition}
\label{fig:pent_eq}
\end{amssidewaysfigure}

\begin{amssidewaysfigure}
\centering
\scalebox{1}{$
\xymatrix@C=15pt@R=55pt{
&& \mathsf{DblCat}(\mathbb A \otimes \mathbb B,\mathbb K)
\ar@/_1.5pc/[llddd]_-{ \mathsf{DblCat}(1 \otimes \lambda,\mathbb K)}
\ar@/^1.5pc/[rrddd]^-{ \mathsf{DblCat}(\varrho \otimes 1,\mathbb K)} \\
&& \mathsf{DblCat}(\mathbb A , \ldb  \mathbb B, \mathbb K\rdb)
\ar[u]_-\cong
\ar[ld]_-{\mathsf{DblCat}(\mathbb A, \ldb \lambda,1 \rdb)\quad}
\ar[rd]^-\cong 
\ar@{}[d]|-{\eqref{eq:triang}} \\
& \mathsf{DblCat}(\mathbb A ,\ldb  \mathbbm 1\otimes  \mathbb B,
\mathbb K\rdb)
\ar[rr]_-{\mathsf{DblCat}(\mathbb A ,\mathfrak a^{\mathbb K})}
\ar[ld]_-\cong &&
\mathsf{DblCat}(\mathbb A , \ldb  \mathbbm 1,\ldb  \mathbb B,\mathbb K\rdb\rdb)
\ar[rd]^-\cong \\
\mathsf{DblCat}(\mathbb A \otimes (\mathbbm 1  \otimes  \mathbb B),
\mathbb K)
\ar[rrrr]_-{\mathsf{DblCat}(\alpha,\mathbb K)} &&&&
\mathsf{DblCat}((\mathbb A \otimes  \mathbbm 1) \otimes 
\mathbb B , \mathbb K)}$}
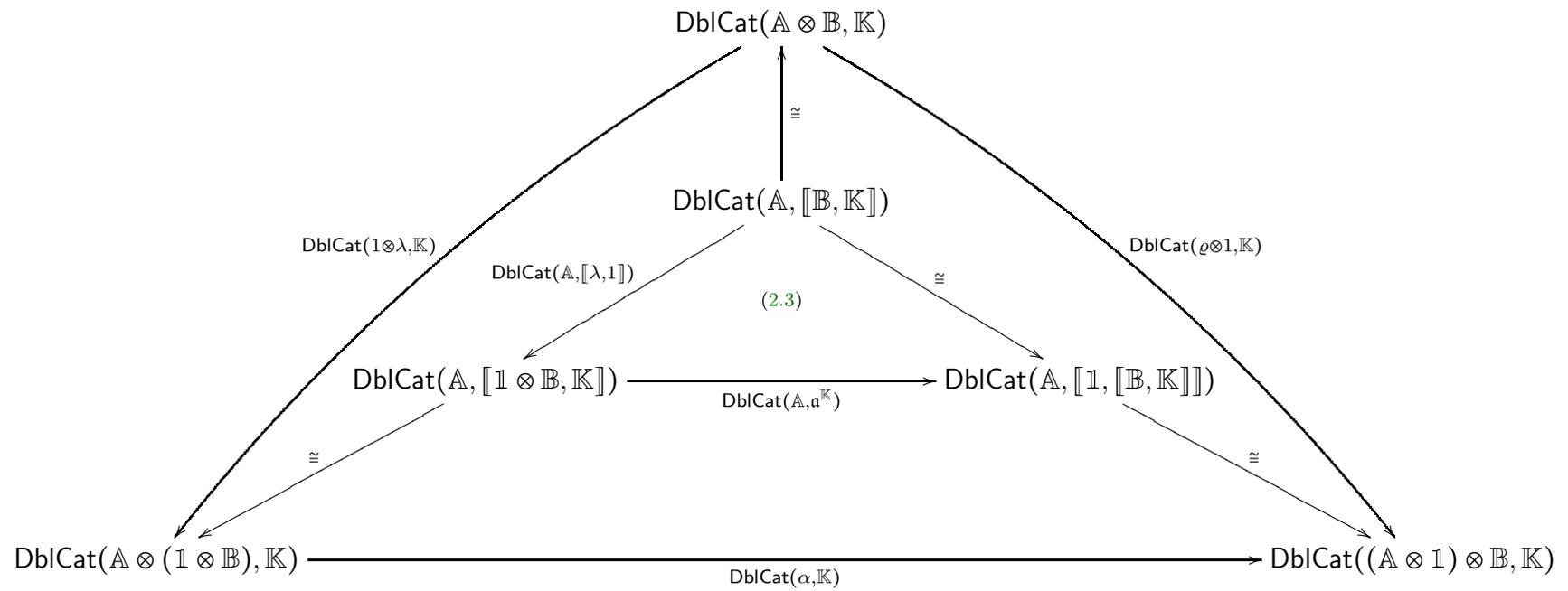
\caption{Equivalent forms of the triangle condition}
\label{fig:triang}
\end{amssidewaysfigure}

\begin{amssidewaysfigure}
\centering
%{\color{white} .}
%\hspace{-2cm}
\scalebox{1}{$
\xymatrix@C=15pt@R=55pt{
% 1.1
\mathsf{DblCat}(\mathbb B \! \otimes \! (\mathbb C \! \otimes \! \mathbb A), \mathbb K)
\ar[rr]^-{\mathsf{DblCat}(1 \otimes \varphi, \mathbb K)}
\ar[dddd]_-{\mathsf{DblCat}(\alpha, \mathbb K)} &&
% 1.3
\mathsf{DblCat}(\mathbb B \! \otimes \! (\mathbb A \! \otimes \! \mathbb C), \mathbb K)
\ar[rr]^-{\mathsf{DblCat}(\alpha, \mathbb K)} &&
\mathsf{DblCat}((\mathbb B \! \otimes \! \mathbb A) \! \otimes \! \mathbb C, \mathbb K)
\ar[dddd]^-{\mathsf{DblCat}(\varphi \otimes 1, \mathbb K)} \\
% 2.2
&  \mathsf{DblCat}(\mathbb B,\ldb \mathbb C \! \otimes \! \mathbb A, \mathbb K \rdb)
\ar[lu]_-\cong
\ar[r]^-{\raisebox{10pt}{${}_{
\mathsf{DblCat}(\mathbb B,\ldb \varphi, \mathbb K\rdb)}$}}
\ar[d]_-{\mathsf{DblCat}(\mathbb B, \mathfrak a^{\mathbb K})} 
\ar@{}[rrd]_(.25){\textrm{Figure}~\ref{fig:hex_eq}}&
% 2.3
\mathsf{DblCat}(\mathbb B,\ldb \mathbb A \! \otimes \! \mathbb C, \mathbb K \rdb)
\ar[u]_-\cong
\ar[r]^-{\raisebox{10pt}{${}_{
\mathsf{DblCat}(\mathbb B,\mathfrak a^{\mathbb K})}$}} &
% 2.4
\mathsf{DblCat}(\mathbb B,\ldb \mathbb A ,\ldb \mathbb C, \mathbb K \rdb\rdb)
\ar[dd]^-{\mathfrak f^{\ldb \mathbb C, \mathbb K \rdb}_0}
\ar[ru]^-\cong \\
% 3.2
& \mathsf{DblCat}(\mathbb B,\ldb \mathbb C ,\ldb \mathbb A, \mathbb K \rdb\rdb)
\ar[d]_-\cong
\ar[rru]_-{\mathsf{DblCat}(\mathbb B,\mathfrak f^{\mathbb K})}
\ar@{}[rr]|-{\eqref{eq:a_inv}} && \\
% 4.2
& \mathsf{DblCat}(\mathbb B\! \otimes \!  \mathbb C ,\ldb \mathbb A, \mathbb K \rdb)
\ar[ld]_-\cong
\ar[r]^-{\raisebox{10pt}{${}_{\mathfrak f^{\mathbb K}_0}$}} &
% 4.3
\mathsf{DblCat}(\mathbb A ,\ldb \mathbb B\! \otimes \!  \mathbb C, \mathbb K \rdb)
\ar[d]^-\cong
\ar[r]^-{\raisebox{10pt}{${}_{
\mathsf{DblCat}(\mathbb A,\mathfrak a^{\mathbb K})}$}} & 
% 4.4
\mathsf{DblCat}(\mathbb A,\ldb \mathbb B ,\ldb \mathbb C, \mathbb K \rdb\rdb)
\ar[rd]^-\cong \\
% 5.1
\mathsf{DblCat}((\mathbb B \! \otimes \! \mathbb C) \! \otimes \! \mathbb A, \mathbb K)
\ar[rr]_-{\mathsf{DblCat}(\varphi, \mathbb K)}&&
% 5.3
\mathsf{DblCat}(\mathbb A \! \otimes \!(\mathbb B \! \otimes \! \mathbb C), \mathbb K)
\ar[rr]_-{\mathsf{DblCat}(\alpha, \mathbb K)} &&
% 5.5
\mathsf{DblCat}((\mathbb A \! \otimes \! \mathbb B) \! \otimes \! \mathbb C, \mathbb K)}$}
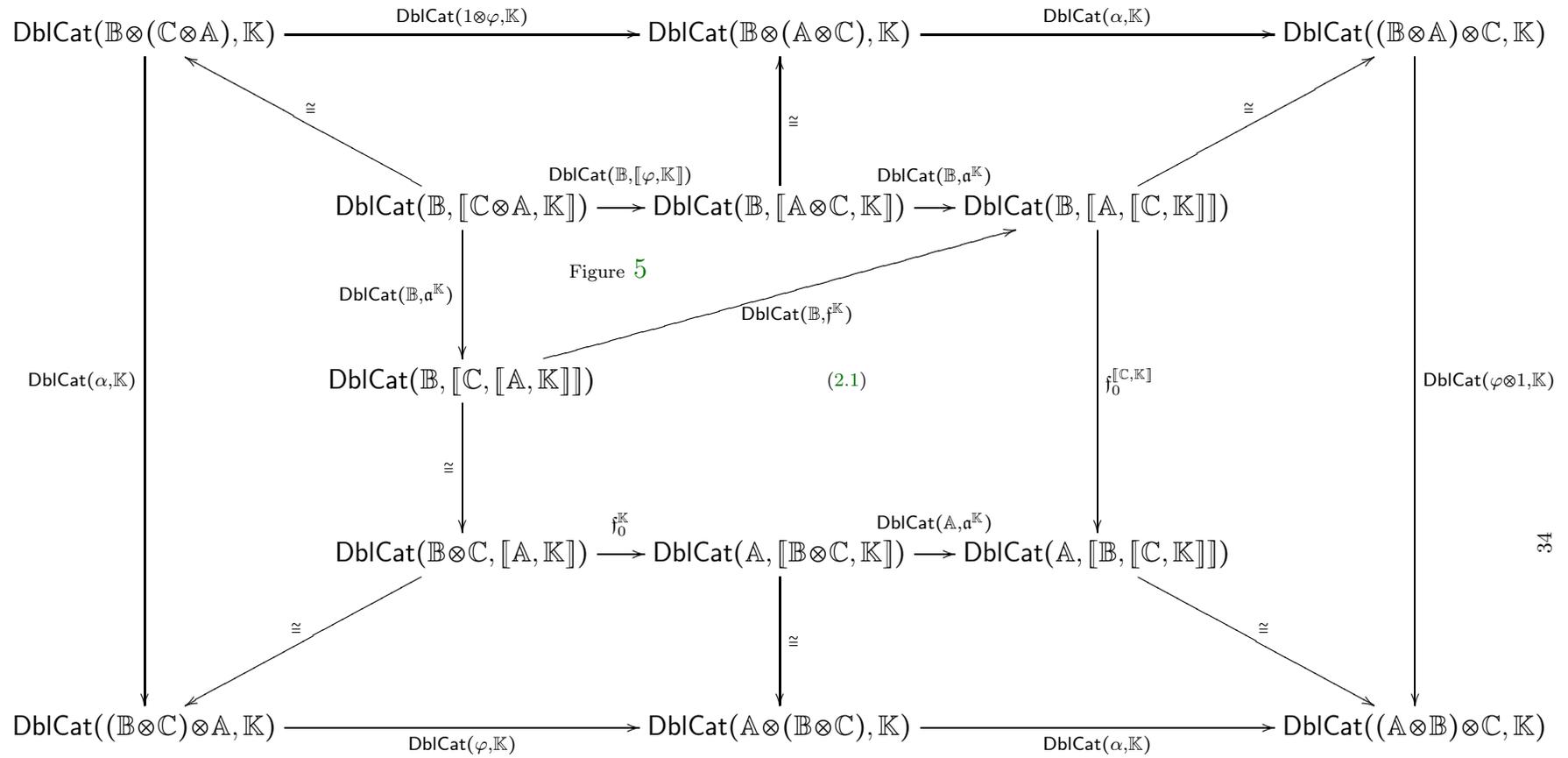
\caption{Equivalent forms of the hexagon condition}
\label{fig:hex}
\end{amssidewaysfigure}

\begin{amssidewaysfigure}
\thisfloatpagestyle{empty}
\centering
%{\color{white} .} 
%\hspace{-2.7cm}
\scalebox{.93}{$
\xymatrix@C=4pt@R=55pt{
% 1.1
\ldb \mathbb C \! \otimes \! \mathbb A \tp  \mathbb K\rdb 
\ar@/^3pc/[rrrr]^-{\ldb\varphi,1\rdb}
\ar[rr]_(.6){\raisebox{-10pt}{${}_{\ldb \epsilon^{\mathbb C},1 \rdb}$}}
\ar[rrd]^-{\mathfrak l^{\mathbb C}}
\ar[rdd]^-{\mathfrak l^{\ldb \mathbb C,\mathbb C \otimes \mathbb A\rdb}}
\ar[ddd]^-{\mathfrak l^{\ldb \ldb\mathbb A,\mathbb C \otimes \mathbb A\rdb,\mathbb C \otimes \mathbb A\rdb}}
\ar@/_11pc/[ddddd]_(.3){\mathfrak l^{\mathbb A}} 
^(.3){\hspace{1.5cm}\eqref{eq:l-r_square}} &&
% 1.3 
\ldb \hspace{-1pt} \ldb \mathbb C\tp \mathbb C \! \otimes \! \mathbb A \rdb\! \otimes \! \mathbb C\tp \mathbb K \rdb 
\ar[r]_-{\raisebox{-10pt}{${}_{\ldb\ldb\eta^{\mathbb A},1\rdb\otimes 1,1\rdb}$}}
\ar[d]^-{\mathfrak a^{\mathbb K}} 
\ar@{}[lld]|(.25){\eqref{eq:eps_a_l}} &
% 1.4
\ldb \hspace{-1pt} \ldb \hspace{-1pt} \ldb \mathbb A\tp \mathbb C \! \otimes \! \mathbb A \rdb \tp
\mathbb C \! \otimes \! \mathbb A \rdb \! \otimes \! \mathbb C\tp \mathbb K \rdb 
\ar[r]_-{\raisebox{-10pt}{${}_{
\ldb\mathfrak r^{\mathbb C \otimes \mathbb A} \otimes 1,1\rdb}$}}
\ar[d]^-{\mathfrak a^{\mathbb K}} &
% 1.5
\ldb \mathbb A \! \otimes \! \mathbb C\tp \mathbb K\rdb 
\ar[d]^-{\mathfrak a^{\mathbb K}} \\
% 2.3
&& \ldb \hspace{-1pt} \ldb \mathbb C\tp \mathbb C \! \otimes \! \mathbb A \rdb\tp \ldb\mathbb C\tp \mathbb K \rdb \hspace{-1pt} \rdb
\ar[r]^-{\raisebox{10pt}{${}_{\ldb\ldb\eta^{\mathbb A},1\rdb,1\rdb}$}}
\ar[d]^-{\mathfrak f^{\mathbb K}} &
% 2.4
\ldb\hspace{-1pt} \ldb \hspace{-1pt} \ldb \mathbb A\tp \mathbb C \! \otimes \! \mathbb A \rdb\tp
\mathbb C \! \otimes \! \mathbb A \rdb\tp\ldb\mathbb C\tp\mathbb K \rdb \hspace{-1pt} \rdb 
\ar[r]^-{\raisebox{10pt}{${}_{
\ldb\mathfrak r^{\mathbb C \otimes \mathbb A} ,1\rdb}$}}
\ar[dd]^-{\mathfrak f^{\mathbb K}} &
% 2.5
\ldb \mathbb A \tp \ldb\mathbb C\tp \mathbb K\rdb \hspace{-1pt} \rdb  
\ar[dddd]^-{\mathfrak f^{\mathbb K}} \\
\ar@{}[r]|-{(\ast)}
% 3.2
& \ldb \hspace{-1pt} \ldb \hspace{-1pt} \ldb \mathbb C\tp \mathbb C \! \otimes \! \mathbb A \rdb \tp
\mathbb C \! \otimes \! \mathbb A \rdb \tp
\ldb \hspace{-1pt} \ldb \mathbb C\tp \mathbb C \! \otimes \! \mathbb A \rdb \tp
\mathbb K \rdb \hspace{-1pt} \rdb 
\ar[r]^-{\raisebox{10pt}{${}_{
\ldb\mathfrak r^{\mathbb C \otimes \mathbb A} ,1\rdb}$}}
\ar[d]^-{\ldb 1,\ldb\ldb\eta^{\mathbb A},1\rdb,1\rdb\rdb} 
\ar@{}[uu]|(.4){\eqref{eq:l-r_pentagon}} &
% 3.3
\ldb \mathbb C \tp \ldb \hspace{-1pt} \ldb \mathbb C \tp \mathbb C \! \otimes \! \mathbb A \rdb \tp
\mathbb K \rdb \hspace{-1pt} \rdb 
\ar[rd]_-{\ldb 1,\ldb\ldb\eta^{\mathbb A},1\rdb,1\rdb\rdb} \\
% 4.1
\ldb \hspace{-1pt} \ldb \hspace{-1pt} \ldb \hspace{-1pt} \ldb \mathbb A\tp \mathbb C \! \otimes \! \mathbb A \rdb \tp
\mathbb C \! \otimes \! \mathbb A \rdb\tp \mathbb C \! \otimes \! \mathbb A \rdb\tp
\ldb \hspace{-1pt} \ldb \hspace{-1pt} \ldb \mathbb A\tp \mathbb C \! \otimes \! \mathbb A \rdb \tp
\mathbb C \! \otimes \! \mathbb A \rdb \tp \mathbb K\rdb \hspace{-1pt} \rdb 
\ar[r]^-{\raisebox{10pt}{${}_{\ldb\ldb\ldb\eta^{\mathbb A},1\rdb,1\rdb,1\rdb}$}}
\ar[d]^-{\ldb\mathfrak r^{\mathbb C \otimes \mathbb A} ,1\rdb} &
% 4.2
\ldb \hspace{-1pt}  \ldb \hspace{-1pt} \ldb \mathbb C \tp \mathbb C \! \otimes \! \mathbb A \rdb \tp
\mathbb C \! \otimes \! \mathbb A \rdb \tp
\ldb \hspace{-1pt} \ldb \hspace{-1pt} \ldb \mathbb A\tp \mathbb C \! \otimes \! \mathbb A \rdb \tp 
\mathbb C \! \otimes \! \mathbb A \rdb \tp \mathbb K\rdb \hspace{-1pt} \rdb
\ar[rr]^-{\raisebox{10pt}{${}_{
\ldb\mathfrak r^{\mathbb C \otimes \mathbb A} ,1\rdb}$}} &&
% 4.4
\ldb \mathbb C \tp
\ldb \hspace{-1pt} \ldb \hspace{-1pt} \ldb \mathbb A\tp \mathbb C \! \otimes \! \mathbb A \rdb \tp
\mathbb C \! \otimes \! \mathbb A \rdb \tp \mathbb K\rdb \hspace{-1pt} \rdb 
\ar[rdd]_-{\ldb1,\ldb\mathfrak r^{\mathbb C \otimes \mathbb A} ,1\rdb\rdb} \\
% 5.1
\quad \ldb \hspace{-1pt} \ldb \mathbb A \tp \mathbb C \! \otimes \! \mathbb A \rdb\tp 
\ldb \hspace{-1pt} \ldb \hspace{-1pt} \ldb \mathbb A\tp \mathbb C \! \otimes \! \mathbb A \rdb \tp
\mathbb C \! \otimes \! \mathbb A \rdb\tp \mathbb K\rdb \hspace{-1pt} \rdb 
\ar[d]^-{\ldb1,\ldb\mathfrak r^{\mathbb C \otimes \mathbb A} ,1\rdb\rdb}\\
% 6.1
\ldb \hspace{-1pt} \ldb \mathbb A\tp  \mathbb C \! \otimes \! \mathbb A \rdb\tp
\ldb \mathbb A\tp \mathbb K \rdb \hspace{-1pt} \rdb 
\ar[rrrr]_-{\ldb \eta^{\mathbb A},1\rdb} &&&&
% 6.5
\ldb \mathbb C \tp \ldb \mathbb A \tp \mathbb K \rdb \hspace{-1pt} \rdb
}$}
\caption{Proof of the hexagon condition}
\label{fig:hex_eq}
\end{amssidewaysfigure}

\end{document}